\documentclass[opew,nonblindrev]{informs2017}
\OneAndAHalfSpacedXI 

\usepackage{endnotes}
\let\footnote=\endnote

\usepackage{amsmath}
\usepackage{amsfonts}
\usepackage{amssymb}
\usepackage{enumerate}
\usepackage{mathrsfs}
\usepackage{graphicx}
\usepackage{makecell}
\usepackage{multirow}
\usepackage{float}
\usepackage{epstopdf}
\usepackage[normalem]{ulem}
\usepackage{bm}
\usepackage{dsfont}
\usepackage{pifont}
\newcommand{\cmark}{\ding{51}}%
\newcommand{\xmark}{\ding{55}}%

\usepackage{hyperref}
\newcommand{\thmref}[2]{\hyperref[#1]{#2 \ref*{#1}}}
\usepackage[dvipsnames]{xcolor}
\newcommand\myshade{100}
\colorlet{mylinkcolor}{MidnightBlue}
\hypersetup{
linkcolor  = mylinkcolor!\myshade!black,
citecolor  = mylinkcolor!\myshade!black,
urlcolor   = mylinkcolor!\myshade!black,
colorlinks = true,
}

\usepackage{graphicx}
\usepackage{subcaption}
\usepackage{caption}
\usepackage{enumitem}
\usepackage{array}
\usepackage{algorithm}
\usepackage{algpseudocode}
\algdef{SE}[SUBALG]{Indent}{EndIndent}{}{\algorithmicend\ }%
\algtext*{Indent}
\algtext*{EndIndent}
\usepackage{soul}

\usepackage{titlesec}
\usepackage[titletoc,toc,title]{appendix}

\usepackage{color}
\usepackage{color-edits}
\addauthor{zc}{cyan}
\addauthor{lf}{magenta}
\addauthor{sm}{brown}

\usepackage{comment}
\definecolor{strcolor}{rgb}{0.6, 0.2, 0.6}
\definecolor{commentcolor}{rgb}{0.3125, 0.5, 0.3125}
\definecolor{keycol}{rgb}{0, 0, 1}

\usepackage{thmtools}
\usepackage{thm-restate}

\usepackage{romannum}
\AtBeginDocument{\pagenumbering{arabic}}

\usepackage{tikz}

\usepackage{adjustbox}
\usepackage[figuresright]{rotating}

\usepackage{mathtools}
\def\sqrtexplained#1{%
  \begingroup
    \sbox0{$#1$}
    \def\underbrace##1_##2{##1}
    \sbox2{$#1$}
    \dimen0=\wd0 \advance\dimen0-\wd2
    \mathrlap{\sqrt{\phantom{\displaystyle#1}\kern\dimen0 }}
    \hphantom{\sqrt{\vphantom{\displaystyle#1}}}
  \endgroup
  #1}

\newcommand{\bmt}[1]{\tilde{\bm{#1}}}

\newcommand{\ubar}[1]{\text{\b{$#1$}}}

\usepackage{tabularx}
\usepackage{pdflscape}
\usepackage{pict2e}
\newcolumntype{L}[1]{>{\raggedright\arraybackslash}p{#1}}
\newcolumntype{C}[1]{>{\centering\arraybackslash}p{#1}}
\newcolumntype{R}[1]{>{\raggedleft\arraybackslash}p{#1}}

\usepackage{natbib}
\bibpunct[, ]{(}{)}{,}{a}{}{,}%
\def\bibfont{\small}%
\TheoremsNumberedThrough   
\ECRepeatTheorems

\EquationsNumberedThrough  

\usepackage{cleveref}
\usepackage{threeparttable}

\begin{document}
\sloppy

\TITLE{Newsvendor under Ambiguity and Misspecification}
\RUNTITLE{Newsvendor under Ambiguity and Misspecification}
\RUNAUTHOR{Liu, Chen, Wang, and Wang}

\ARTICLEAUTHORS{%
\AUTHOR{Feng Liu$^1$, Zhi Chen$^2$, Ruodu Wang$^3$, and Shuming Wang$^1$}
\AFF{$^1$School of Economics and Management, University of Chinese Academy of Sciences \\
Emails: liufeng21@mails.ucas.ac.cn,  wangshuming@ucas.edu.cn \\
$^2$Department of Decisions, Operations and Technology, The Chinese University of Hong Kong \\
Email: zhi.chen@cuhk.edu.hk \\
$^3$Department of Statistics and Actuarial Science, University of Waterloo \\
Email: wang@uwaterloo.ca}
}

\ABSTRACT{\noindent \textbf{Problem definition:} We consider a newsvendor problem with unknown demand distribution, where we distinguish \textit{ambiguity} under which the newsvendor does not differentiate demand distributions of common characteristics (\textit{e.g.}, mean and variance) and \textit{misspecification} under which such characteristics might be misspecified  (due to, \textit{e.g.}, estimation error and/or distribution shift). 

\noindent \textbf{Methodology/results:} The newsvendor hedges against ambiguity and misspecification by maximizing the worst-case expected profit regularized by a distribution's distance to an ambiguity set of distributions with some specified characteristics. Focusing on the popular mean-variance ambiguity set and optimal-transport cost for the misspecification, we show that the decision criterion of misspecification aversion possesses insightful interpretations as distributional transforms. We derive the closed-form optimal order quantity that generalizes the solution of the seminal Scarf model under only ambiguity aversion. We establish finite-sample performance guarantees that consist of two parts: an in-sample optimal value and an out-of-sample effect of misspecification, which can be further decoupled into estimation error and distribution shift. We also extend the framework to multiple products, distributional characteristics specified via optimal transport, and misspecification measured by total variation distance, and derive analytical optimal solutions. 

\noindent \textbf{Managerial implications:} The closed-form solution highlights the impact of misspecification aversion: the optimal order quantity under misspecification aversion can decrease as the price or variance increases, reversing the monotonicity of that under only ambiguity aversion. Hence, ambiguity and misspecification, as different layers of distributional uncertainty, can result in distinct operational consequences. The finite-sample performance guarantee theoretically justifies the need to incorporate misspecification aversion in a non-stationary environment, as demonstrated in our experiments with real-world data.
}

\maketitle

\section{Introduction}

The newsvendor model, as a building block for dealing with uncertain demand in operations management, has found its successful applications in various domains, including inventory management (\citealp{chen2007risk}), revenue management (\citealp{besbes2018dynamic}), capacity planning (\citealp{simchi2015worst}), and healthcare (\citealp{olivares2008structural}), to name a few. When the true demand distribution is fully known, the celebrated \textit{critical fractile} determines an optimal order quantity that maximizes expected profit. In practice, however, the newsvendor often faces incomplete knowledge about the demand distribution. Hence, it can be difficult (if not impossible) to precisely articulate the true demand distribution, leading to demand \textit{ambiguity} for the newsvendors. 

A natural way to mitigate demand ambiguity is to utilize only partial distributional characteristics available for decision-making. In this vein, mean and variance---arguably two of the most widely used and easy-to-estimate statistics that capture the key \textit{location} and \textit{dispersion} characteristics of the underlying distribution, respectively---have been employed. This can be traced back to the seminal work of \cite{scarf1958min} that considers an \textit{ambiguity-averse} newsvendor model maximizing the worst-case expected profit over an ambiguity set of probability distributions with the same mean and variance.\footnotemark\footnotetext{A more general version of the Scarf model allows the mean and variance to be uncertain and to vary within some predetermined intervals, which turns out to be equivalent to the Scarf model characterized by the predetermined lower bound of mean and upper bound of variance (see \citealp{natarajan2011beyond} or Remark~\ref{remark:natarajan}).} In essence, with such a mean-variance ambiguity set, the newsvendor specifies her belief about the (true) demand distribution via an approximation by using mean and variance characteristics. This is also well statistically justified,\footnotemark\footnotetext{It is well known that under mild conditions, a probability distribution can be \textit{uniquely} determined by all integer-order moments (\citealp{durrett2019probability}), among which the first two moments are the most commonly used ones.} especially in the context of the newsvendor problem, by noting that the single-dimensional demand distribution's quantile (\textit{i.e.}, the critical-fractile solution) can be largely characterized by the mean and variance statistics, or even \textit{perfectly} determined under many common distributions (\textit{e.g.}, elliptical, uniform, and exponential, see~\citealp{meyer1987two}). This also underpins the implication that solutions of the ambiguity-averse and nominal ambiguity-neutral models characterized by the same mean and variance can share several key features in operational properties: for instance, the order quantity increases in price (see Sections~\ref{sec:model} and~\ref{sec:sensitivity} for more details). Apart from the newsvendor problem (\citealt{perakis2008regret}, \citealt{han2014risk}),  the mean-variance ambiguity set has been studied in stochastic optimization (\citealp{popescu2007robust}), and used in various applications, including recent literature on decision theory (\citealt{muller2022ranking}), mechanism design (\citealt{carrasco2018optimal, chen2022distribution, chen2024screening}), and risk management (\citealt{li2018worst, nguyen2021mean}).

\begin{figure}[tb]
\includegraphics[width=1. \textwidth]{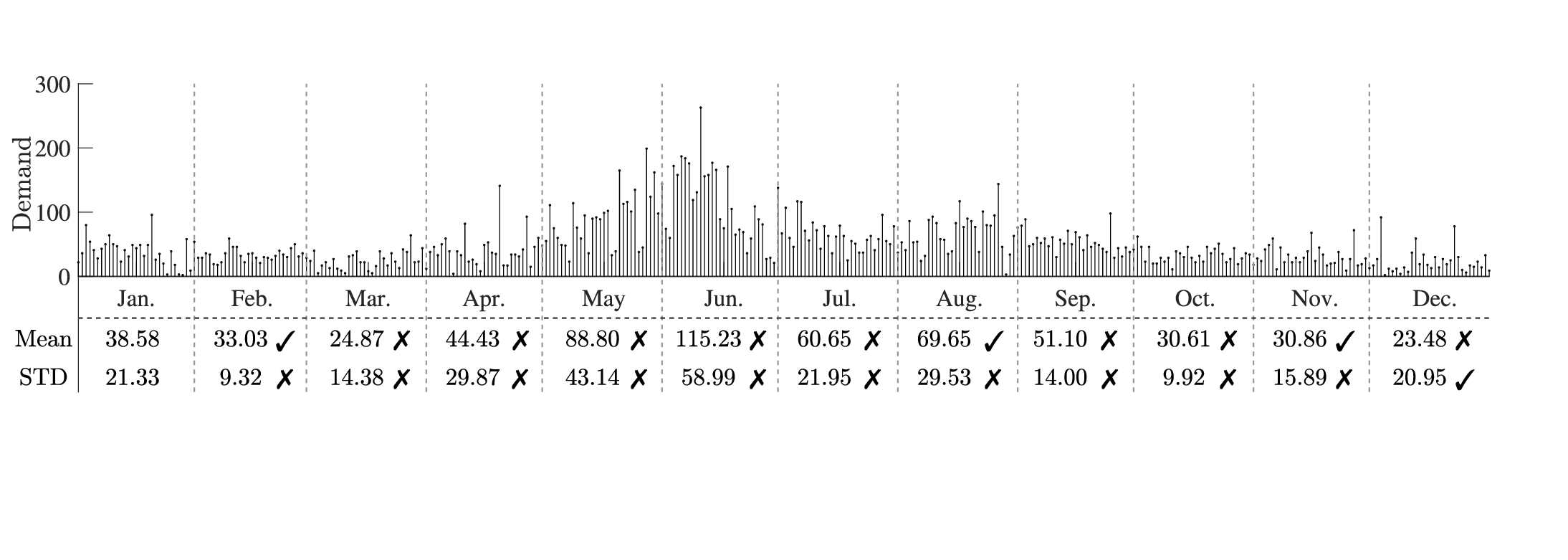}
\vspace{-24mm}
\caption{\textnormal{Daily demand of a product of drinking water over one year with monthly mean and standard deviation (STD). The sign `\cmark' (resp., `\xmark') means that the monthly mean or STD falls in (resp., does not fall in) the 95\% confidence interval estimated with the demand data in the preceding month. Notably, in only $2$ out of $11$ (resp., $1$ out of $11$) instances, the mean (resp., variance) falls within the confidence interval.}}
\footnotesize{{\it Notes.} $\;$ The confidence intervals on the mean and variance are constructed leveraging the $t$-statistic and $\chi^2$-statistic, respectively, without knowing the true values of mean and variance.}
\label{fig:demand data}
\end{figure}

\begin{figure}[tb]
\begin{subfigure}{.39\textwidth}
\centering
\includegraphics[width=0.7\linewidth]{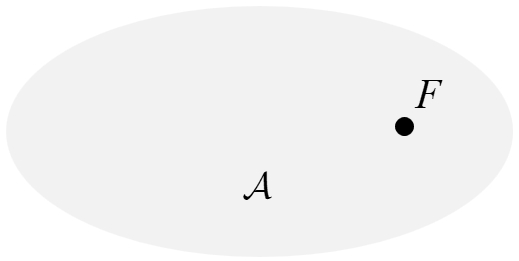}
\end{subfigure}
\begin{subfigure}{.60\textwidth}
\centering
\includegraphics[width=0.7\linewidth]{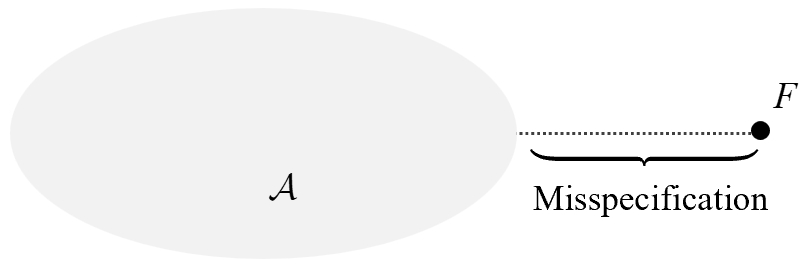}
\end{subfigure}
\caption{\textnormal{\textit{Left}: The true demand distribution $F$ resides in the ambiguity set $\mathcal{A}$, and the newsvendor faces only ambiguity. \textit{Right}: The true demand distribution resides outside the ambiguity set, and hence, the newsvendor faces both ambiguity and misspecification.}}
\label{fig:ambiguity misspecification}
\end{figure}

However, in many practical situations, it is still challenging to accurately estimate the demand mean and variance. For instance, in retailing industries (\textit{e.g.}, E-commerce, grocery, and supermarket), more and more commodities are of ever shorter life cycles (\citealt{calvo2016sourcing}), and many enterprises lack the resources for effective data collection and analysis, especially for new products (\citealt{feiler2022noise}). Insufficient historical data---even when the underlying demand process is stationary---can result in non-negligible {\it estimation error} that challenges the newsvendors' decision under demand ambiguity. On the other hand, the uncertainty is exacerbated if the underlying demand process becomes non-stationary, due to, for instance, the complicated  (time-varying) determinants for the demand~\citep{an2025nonstationary}. This can lead to {\it distribution shift}---the future demand distribution is different from the past\footnotemark\footnotetext{The phenomenon of distribution shift, in statistical learning, also refers to situations where the training and testing samples are governed by different distributions (see, \textit{e.g.}, \citealp{joaquin2008dataset}).}---under which the mean and variance characteristics can constantly change over time, making them difficult to estimate from historical data. In~\Cref{fig:demand data}, we illustrate that even the confidence interval of either mean or variance may depreciate in a non-stationary demand environment. In sum, either estimation error or distribution shift can cause \textit{misspecified} mean and variance parameters, and consequently, the optimal order quantity prescribed under ambiguity may perform inexpertly. This is known as \textit{model misspecification} in statistics and economics (\citealp{hansen2014nobel}). In our context, it refers to the possibility that the true demand distribution may not reside in the ambiguity set with the same distributional characteristics as specified (\textit{e.g.}, mean and variance in the Scarf model or the neighborhood of some reference distribution as in \Cref{sec:ambiguity set}); see~\Cref{fig:ambiguity misspecification} for an illustration.

The above discussion motivates us to caution against the potential issue of model misspecification in the fundamental newsvendor problem. To this end, we introduce a misspecification-averse (and ambiguity-averse) newsvendor model, which is inspired by a recently developed framework of globalized distributionally robust optimization~\citep{liu2023globalized} that can mitigate the performance shortfall due to the possibility that the true distribution falls outside the ambiguity set. Our newsvendor model is also well supported by a recent axiomatic framework~\citep{cerreia2025making} that unifies behavioral decision criteria averse to either ambiguity or misspecification. In particular, we distinguish the ambiguity under which the newsvendor does \textit{not} differentiate demand distributions of common distributional characteristics (\textit{e.g.}, mean and variance) and misspecification under which such characteristics might be misspecified (due to, \textit{e.g.}, estimation error and/or distribution shift as discussed before). The newsvendor hedges against ambiguity and misspecification by maximizing the worst-case expected profit regularized by a distribution's distance to an ambiguity set of distributions with some specified characteristics. We investigate---from decision-criterion, operational, and statistical perspectives---how misspecification aversion may affect the newsvendor's decision and what is the rationale behind it. 

\subsection{Summary of Main Contributions}

\vspace{2mm}
\noindent \textbf{Newsvendor under misspecification}. We investigate the fundamental newsvendor problem under a \textit{structured} decision-under-uncertainty framework that distinguishes the ambiguity and misspecification. We focus on the mean-variance information for characterizing the ambiguity and optimal-transport cost for quantifying the misspecification, respectively, and investigate the rationale behind the misspecification aversion that differentiates from the ambiguity aversion of the seminal Scarf framework. We extend our modeling on the ambiguity of demand distribution via bounding its statistical distance to a reference distribution, misspecification aversion upon ambiguity-aversion is then equivalent to a stronger misspecification aversion to a nominal ambiguity-neutral model under the reference distribution (Section~\ref{sec:ambiguity set}). We also extend to consider other statistical distances (specifically, total variation distance) to measure the misspecification (Section~\ref{sec:phi divergence}). These, together with misspecification aversion to the mean-variance ambiguity set, achieve a comprehensive inspection of the newsvendor under ambiguity and misspecification. 

\vspace{2mm}
\noindent \textbf{Interpretation via distributional transform.}
We investigate and interpret the decision criterion of misspecification aversion by leveraging a decision-analysis vehicle of \textit{distributional transform}~\citep{liu2021distributional}. We show that the decision criterion of the newsvendor under ambiguity and misspecification is essentially a worst-case \textit{transformed} expected profit with each distribution \textit{within} the ambiguity set being transformed---via a {\it transform function} determined by the index of misspecification aversion and operational information (\textit{e.g.}, price $p$ and order quantity $q$) of the newsvendor---to another distribution possibly \textit{beyond} the ambiguity set (\Cref{thm:distributional transform}). 
In this way, the misspecification-aversion effect with valuable information on cost structure (\textit{e.g.}, price) and ordering decision is encapsulated in the transform function along the distributional transform.  

\vspace{2mm}
\noindent \textbf{Analytical optimal solutions.} $\;$
We analytically derive the optimal order quantity of the newsvendor under ambiguity and misspecification (\Cref{thm:newsvendor transport-2}), which generalizes that of the seminal Scarf model under ambiguity alone. The analytical solution enables us to analyze the optimal order quantity's sensitivity to cost structure (\textit{i.e.}, price and cost) and distributional characteristics (\textit{i.e.}, mean and variance) to understand the impact and rationale of misspecification aversion. In particular, while it is \textit{always} optimal to order \textit{more} as price increases under ambiguity aversion, it turns out that ordering \textit{less} instead can be optimal as the price increases under misspecification aversion, \textit{ceteris paribus} (Proposition~\ref{prop:price}). Likewise, in the case of high-profit margin, while it is optimal to order \textit{more} as variance increases under ambiguity aversion, ordering \textit{less} can be optimal as the variance increases under misspecification aversion, \textit{ceteris paribus} (Proposition~\ref{prop:variance}).  These observations reveal that ambiguity and misspecification, as different layers of distributional uncertainty, can result in \textit{distinct} operational consequences, and therefore should be distinguished in the modeling. Furthermore, we extend to derive the analytical optimal order quantities of multiple products under a sum-of-variance constraint, which \textit{unifies} the ambiguity-averse and misspecification-averse single-product model, ambiguity-averse multi-product model, and the Scarf model (\Cref{thm:multi product}). 

\vspace{2mm}
\noindent \textbf{Performance guarantee.} $\;$ 
We establish the finite-sample performance guarantee of the optimal solution that decouples the mixing effect of estimation error and distribution shift in the misspecification statistically (\Cref{thm:guarantee}). In particular, the estimation error is related to the distance between the data-generating distribution and the estimated mean-variance ambiguity set, and it \textit{diminishes} as the sample size approaches infinity; while the distribution shift, captured by the distance between the data-generating distribution and the out-of-sample distribution, is \textit{independent} of the sample size. This theoretically justifies the \textit{rationale} of incorporating misspecification aversion in a non-stationary environment, which is also well demonstrated in our experiments with real-world retailing data. 

\subsection{Related Works}

Our work is related to the following streams of literature. 

\vspace{2mm}
\noindent \textbf{Newsvendor with moment condition.} $\;$
Since the pioneering work of \cite{scarf1958min}, various studies have employed the moment condition to specify the ambiguity-averse newsvendor model with limited demand information. For instance, building on the results of \cite{scarf1958min}, \cite{gallego1993distribution} extend to consider a multi-product newsvendor problem based on \textit{marginal} mean-and-variance information under a budget constraint. \cite{natarajan2011beyond} generalize the Scarf model by allowing the mean and variance to be uncertain and to vary within some predetermined intervals. This formulation turns out to be equivalent to the Scarf model characterized by the lower bound of mean and the upper bound of variance. \cite{zhu2013newsvendor} focus on minimizing the worst-case regret under known mean and variance of the random demand. In addition to ambiguity aversion, \cite{han2014risk} incorporate risk aversion into the newsvendor problem, where the risk is captured by the standard deviation of the newsvendor's profit. Given the mean and variance of the demand, they develop a closed-form solution for the risk-averse (and ambiguity-averse) newsvendor. Apart from mean and variance, other moment information has also been considered in ambiguity-averse newsvendor problems. We refer to \cite{perakis2008regret} for structural information such as median, unimodality, and symmetry, to \cite{ardestani2016robust} for first-order partial moments, to \cite{natarajan2018asymmetry} for asymmetry based on second-order partitioned statistics, to \cite{das2021heavy} for the $t$-th ($t \geq 1$) moment that can capture heavy-tailed demand distributions, and to \cite{padmanabhan2021exploiting} for partial correlations among products. Specifically, \citet{das2021heavy} interpret the Scarf model through a heavy-tailed demand perspective, and reveal that the robust newsvendor problem can be represented under an implied heavy-tailed distribution outside the moment ambiguity set. \cite{govindarajan2021distribution} shift the focus to the inventory pooling problem, where the ambiguity set is specified by mean and covariance (see also \citealt{hanasusanto2015distributionally}), and characterize the closed-form solution of the two-location model. 

\vspace{2mm}
\noindent \textbf{Newsvendor with statistical distance.} $\;$
Another stream of the ambiguity-averse newsvendor model is based on ambiguity sets specified through the closeness to a reference distribution in terms of a certain statistical distance. For instance, \cite{rahimian2019controlling} delve into the total variation distance and obtain an insightful closed-form solution. Based on the Wasserstein distance, \cite{chen2021regret} adopt the minimax regret decision criterion in the presence of both demand and yield randomness, and they show that the optimal order quantity can be determined via an efficient golden section search. \cite{zhang2023optimal} and \cite{fu2024distributionally} further leverage side information from explanatory features, and derive an optimal analytical ordering policy based on the Wasserstein distance and a closed-form solution based on the JW discrepancy measure, respectively. In this work, we also consider the possibility that the misspecification may arise from such ambiguity-averse newsvendor models with statistical distance (\Cref{thm:double misspecification}).

\vspace{2mm}
\noindent \textbf{Model misspecification.} $\;$ 
Statisticians and econometricians have long grappled with the challenge of addressing uncertainty in decision-making, which is categorized as risk, ambiguity, and misspecification (see, \textit{e.g.}, \citealt{hansen2014nobel, hansen2022structured}). Several noteworthy contributions have been made, which, from many angles, shed light on the interplay between the different layers of uncertainty. From the empirical perspective, \cite{aydogan2023three} provide experimental evidence for the role of model misspecification in decision-making under uncertainty. Their work establishes a compelling case for recognizing a distinct behavioral impact among risk, ambiguity, and misspecification, illuminating the nuanced facets of decision theory. From the theoretical perspective, \cite{cerreia2025making} propose an innovative axiomatic framework that unifies the behavioral decision criteria on the aversion to ambiguity and/or misspecification. From the perspective of optimization under uncertainty, the works of \cite{ben2006extending} and \cite{ben2017globalized} propose paradigms named globalized/comprehensive robust counterparts to mitigate the violation of uncertainty-involved constraints when the true parameter does not reside in the pre-specified uncertainty set. Recently, \cite{liu2023globalized} extend the notion of globalized robustness and propose the globalized distributionally robust counterpart to mitigate the violation when the true distribution does not reside in the pre-specified ambiguity set. In addition, \cite{long2023robust} introduce a target-oriented model of robust satisficing that maximizes the robustness to uncertainty of achieving a satisfactory target. Also, the distributionally robust models of \cite{delage2010distributionally} and \cite{wiesemann2014distributionally} can mitigate the possible variations in the moment information and in the probabilities of specific events, respectively. There are also nascent works investigating the issue of model misspecification in revenue management, such as a misspecified demand function (\citealt{nambiar2019dynamic}) and a misspecified choice model (\citealt{chen2024robust}). 


\subsection{Research Gaps}

Although our work is closely related to the studies of \cite{cerreia2025making} and \cite{liu2023globalized}, the gap can be well established.  It is noted that \cite{cerreia2025making} focuses primarily on developing the decision-theoretical axioms for the model misspecification, and \cite{liu2023globalized} develops the framework of globalized distributionally robust optimization and focuses on deriving the associated tractable model counterparts. However, our study is devoted to the analysis of the newsvendor problem under both ambiguity and misspecification. In particular, our developments (for the newsvendor problem)---on the decision-criterion interpretation from the perspective of distributional transform, closed-form and analytical solutions, operational insights, and the performance guarantee that decouples the effects of estimation error and distribution shift---have not been attempted in the above studies. We also point out that although the ambiguity set proposed in \cite{delage2010distributionally} also captures the variation in mean and variance, it is computationally inconvenient for the newsvendor problem. In contrast, our framework provides a systematic and tractable way to account for misspecification while enabling clear operational insights. Moreover, our approach can be readily adapted---by replacing the Scarf ambiguity set with theirs---to further mitigate the misspecification to the bounds on the moment information in their ambiguity set. Although both our work and \cite{das2021heavy} admit an equivalent stochastic representation under an implied distribution outside the prescribed ambiguity set, the interpretations differ fundamentally. \cite{das2021heavy} retroactively interpret Scarf’s ambiguity-averse solution through a heavy-tailed demand distribution, whereas our model proactively incorporates misspecification aversion into the decision criterion, yielding a transformed distribution that depends jointly on the order quantity, price, the index of misspecification aversion, and moment information. This shift in perspective yields distinct optimal solutions, operational insights, and performance guarantees, which constitute the main contribution of our paper.

\subsection{Notation}

We denote by $\mathcal{P}$ (resp., $\mathcal{M}_+$) the set of probability measures (resp., non-negative measures) supported on $\mathbb{R}_+$, and $\mathcal{P}_0$ the set of probability measures supported on $\mathbb{R}$. We use $\tilde{v} \sim F$ to signify a random variable $\tilde{v}$ that follows the probability distribution with a cumulative probability distribution (CDF) $F$, under which the expectation is $\mathbb{E}_{F}[\cdot]$. The Dirac distribution at $v \in \mathbb{R}$ is denoted by $\delta_v$.  We adopt the convention that $1/0 = \infty$, and refer to ``decreasing/increasing'' in the weak sense.

\section{Model}\label{sec:model}

The newsvendor decides the order quantity before demand realization and tries to maximize her expected profit. Given the unit price $p$, unit cost $c$ ($c<p$), and an order quantity $q$, the profit under the materialized demand $v$ amounts to $\pi(q,v) = p \cdot \min\{q, v\} - cq = p \cdot \min\{q, v\} - (1 - \kappa)pq$, where we denote the profit margin by $\kappa = \frac{p-c}{p}$ (giving $c = (1-\kappa)p$), which will be an important parameter for our analysis. Facing stochastic demand, the newsvendor must navigate optimizing her decision-making process to balance the trade-off between lost sales (incurred when $v > q$) and excess inventory (incurred when $v \leq q$). When the precise demand distribution $G$ is fully known, the newsvendor maximizes her expected profit by solving the problem
\begin{equation}\label{prob:nominal}\tag{\sc Nominal}
\max_{q\geq 0}~\mathbb{E}_G[\pi(q,\tilde{v})].
\end{equation}
The optimal order quantity $q^\star_G$ is characterized by the classic critical fractile $q^\star_G = G^{-1}(\kappa)$ (\citealp{arrow1951optimal}). If $G$ is an elliptical distribution with location parameter $\mu$, scale parameter $\sigma$, and some density generator $\xi(\cdot)$,\footnotemark\footnotetext{The density function of the elliptical distribution is $f(x) = \frac{1}{\sigma} \cdot \xi\big( \frac{(x - \mu)^2}{\sigma^2} \big)$ with $\int_{-\infty}^{\infty} \xi(z^2)\; {\rm d}z = 1$.} then $q_G^\star = \mu + \sigma \cdot \Xi^{-1}(\kappa)$, where $\Xi(u) = \int_{-\infty}^u\xi(v^2){\rm d}v$. Such a mean-variance formula holds for many other classes of distributions, including the uniform distributions and the exponential distributions. It is worth mentioning that within each class above, the optimal order quantity is also uniquely determined by only the mean and variance of the demand distribution.

In practice, full information on the demand distribution is often not accessible. To tackle the challenge of partial distributional information, significant advancements have been made in prescribing an ambiguity-averse solution that remains robust against all distributions specified by some common distributional characteristics. The seminal work of \cite{scarf1958min} specifies only the mean and variance of the demand distribution and solves
\begin{equation}\label{prob:ambiguity}\tag{\sc Ambiguity}
\max_{q \geq 0} \; \min_{G \in \mathcal{A}} \; \mathbb{E}_{G}[\pi(q,\tilde{v})]
~~{\rm with}~~
\mathcal{A} = \{G \in \mathcal{P} \mid \mathbb{E}_{G}[\tilde{v}] = \mu,\;\mathbb{E}_{G}[\tilde{v}^2] = \mu^2+\sigma^2\}.
\end{equation}
The \ref{prob:ambiguity} model admits an analytical solution---for the reason that will become clear subsequently, we emphasize it with a subscript `$\infty$'---as follows:
\begin{equation}\label{equ:optimal order ambiguity}
q_\infty^{\star} = \left\{
\begin{array}{ll}
\displaystyle \mu + \sigma f(1 - \kappa) & \displaystyle~~~\kappa \geq \frac{\sigma^2}{\mu^2+\sigma^2} \\ [3mm]
0 & \displaystyle~~~\kappa < \frac{\sigma^2}{\mu^2+\sigma^2}
\end{array}
\right. ~~~{\rm with}~~~ f(x) = \frac{1 - 2x}{2\sqrt{x(1-x)}}~~~\forall 0<x<1.
\end{equation}
As illustrated in \cite{scarf1958min}, the optimal order quantity $q_\infty^{\star}$ of \ref{prob:ambiguity} is comparatively close to the optimal order quantity of \ref{prob:nominal} under a normal approximation of the Poisson distribution for a moderate range of profit margins. In effect, in many commonly used distributions (\textit{e.g.}, elliptical, uniform, and exponential) of a nominal $G$ as mentioned above, the optimal order quantities $q_\infty^\star$ and $q_G^\star$ share important sensitivity features (to be discussed in Section~\ref{sec:sensitivity}). In addition, the \ref{prob:ambiguity} model can also cover the situation where the mean and variance are themselves uncertain and reside in some estimated intervals, as remarked below.   
\begin{remark}[\textnormal{{\sc Generality of \ref{prob:ambiguity}}}]\label{remark:natarajan}
{\rm Consider an ambiguity set parameterized by bounds on the uncertain mean and variance:
\[
\mathcal{V} = 
\big\{
G \in \mathcal{P} \mid \mathbb{E}_{G}[\tilde{v}] = \mu,\;\mathbb{E}_{G}[\tilde{v}^2] = \mu^2+\sigma^2 ~~\mbox{for some}~~ \mu \in [\ubar{\mu},\bar{\mu}] ~~\mbox{and}~~ \sigma^2 \in [\ubar{\sigma}^2,\bar{\sigma}^2]
\big\}.
\]
By corollary~5.1 of \cite{natarajan2011beyond}, it holds that
\[
\max_{q \geq 0} \; \min_{G \in \mathcal{V}} \; \mathbb{E}_{G}[\pi(q,\tilde{v})] = \max_{q \geq 0} \; \min_{G \in \mathcal{A}} \; \mathbb{E}_{G}[\pi(q,\tilde{v})] ~~{\rm with}~~ \mathcal{A} = 
\{
G \in \mathcal{P} \mid \mathbb{E}_{G}[\tilde{v}] = \ubar{\mu},\;\mathbb{E}_{G}[\tilde{v}^2] = \ubar{\mu}^2+\bar{\sigma}^2
\}.
\]
The right-hand side equivalence is indifferent to~\ref{prob:ambiguity} with mean $\ubar{\mu}$ and variance $\bar{\sigma}^2$. Similarly, if the moment ambiguity set developed in \cite{delage2010distributionally} is employed, the bounds on the mean and variance could also be misspecified.
}
\end{remark}

An underlying assumption of~\ref{prob:ambiguity} is that the true distribution $F_{\rm true}$ falls in the specified ambiguity set $\mathcal{A}$. However, the ambiguity set $\mathcal{A}$ is constructed using the mean and variance estimates for those of the true $F_{\rm true}$, and is, therefore, may be misspecified, that is, $F_{\rm true} \notin \mathcal{A}$ (see also the empirical evidence in \Cref{fig:demand data}). This leads to \textit{misspecification}---an issue we address next.

To capture misspecification formally and in a computationally appealing fashion,\footnotemark\footnotetext{When considering misspecification aversion, \cite{cerreia2025making} focus on $\phi$-divergence, which is generally computationally difficult in the newsvendor problem, except for the total variation distance as discussed in Section~\ref{sec:phi divergence}. Another alternative is the Gelbrich distance (\citealp{gelbrich1990formula}), which, however, does not lead to an analytical solution.} for a demand distribution $F \in \mathcal{P}$, we measure its closeness to the ambiguity set $\mathcal{A}$ by
\[
d(F,\mathcal{A}) = \min_{G \in \mathcal{A}} \; d_{\rm OT}(F,G)
\]
with $d_{\rm OT}(F, G)$ being the optimal-transport cost (\citealp{villani2009optimal}) between two distributions $F$ and $G$ with a quadratic cost function $|\cdot|^2$ defined as
\begin{equation}\label{equ:optimal transport 2}
d_{\rm OT}(F, G) =\min_{\Gamma \in \mathcal{W}(F,G)} \; \int_{\mathbb{R}_+\times\mathbb{R}_+}|u-v|^2{\rm d}\Gamma(u,v),
\end{equation}
where $\mathcal{W}(F,G)$ is the set of joint probability distributions on $\mathbb{R}_+ \times \mathbb{R}_+$ with marginals $F$ and $G$. The quantity $\sqrt{d_{\rm OT}(F, G)}$ is also known as the type-2 Wasserstein distance between $F$ and $G$. Note that $d(F,\mathcal{A}) > 0$ if and only if $F \notin \mathcal{A}$. In the main content, we focus on the optimal-transport cost in characterizing misspecification, for the sake of tractability and statistical convenience (see more details in \Cref{sec:optimal solution} and \Cref{sec:shadow value}). Following the spirit of~\cite{cerreia2025making}, given the ambiguity set $\mathcal{A}$, we incorporate misspecification into the newsvendor's decision criterion so that a misspecification-averse (and ambiguity-averse) newsvendor solves
\begin{equation}\label{prob:misspecification}\tag{\sc Misspecification}
\Pi_\alpha^\star = \max_{q \geq 0} \; \min_{F \in \mathcal{P}} \; \Big\{\mathbb{E}_{F}[\pi(q,\tilde{u})] + \alpha \cdot d(F,\mathcal{A})\Big\}
\end{equation}
for some $\alpha \geq 0$ that represents the index of misspecification aversion: the lower the index, the stronger the aversion to misspecification. Intuitively speaking, a larger value of $\alpha$ places a larger penalty on deviation from the ambiguity set $\mathcal{A}$ (as measured by $d(F,\mathcal{A})$) and corresponds to higher confidence in $\mathcal{A}$ (or equivalently, the \ref{prob:ambiguity} model). On the one end, when $\alpha \to \infty$, misspecification aversion is absent (see section~4.1 in \citealt{cerreia2025making}), and \ref{prob:misspecification} reduces to \ref{prob:ambiguity} as the newsvendor is absolutely confident with $\mathcal{A}$. On the other end, when $\alpha \to 0$, misspecification aversion is strongest, and \ref{prob:misspecification} reduces to the robust model
\[
\max_{q \geq 0} \; \min_{F \in \mathcal{P}} \; \mathbb{E}_{F}[\pi(q,\tilde{u})] = \max_{q \geq 0} \; \min_{u \in \mathbb{R}_+} \; \pi(q,u),
\]
wherein the newsvendor is so unconfident that she disregards the distributional characteristics specified in $\mathcal{A}$. An analogous observation to Remark~\ref{remark:natarajan} for \ref{prob:ambiguity} can be established for \ref{prob:misspecification}: misspecification over the uncertain mean and variance is still equivalent to misspecification over the lowest mean and highest variance; see \Cref{appendix:equivalence}.

It is worth noting that the \ref{prob:misspecification} problem can also be regarded as the Lagrangian relaxation of an alternate formulation of the newsvendor model under misspecification as follows:
\begin{equation}\label{prob:misspecification constrained}
\max_{q \geq 0} \; \min_{d(F,\mathcal{A}) \leq \varepsilon} \; \mathbb{E}_F[\pi(q,\tilde{u})]
\end{equation}
for some $\varepsilon \geq 0$. When $\varepsilon = 0$, problem~\eqref{prob:misspecification constrained} reduces to \ref{prob:ambiguity}; and when $\varepsilon > 0$, the optimal solution to problem~\eqref{prob:misspecification constrained} can be constructed from that of~\ref{prob:misspecification} (\Cref{lemma:misspecification constrained}, \Cref{sec:techniques}). In other words, the decision of~\ref{prob:misspecification} essentially hedges against some worst-case distribution that can stay \textit{outside} the ambiguity set $\mathcal{A}$ characterized by mean and variance. This makes~\ref{prob:misspecification} distinct from~\ref{prob:ambiguity} (with the underlying worst-case distribution within $\mathcal{A}$) in many aspects, which we further explore in the forthcoming sections. 
In Section~\ref{extension}, we extend  \ref{prob:misspecification} to involve multiple products, distance-based ambiguity set defined via optimal transport, and misspecification measured by total variation distance.

\section{Decision Criterion via Distributional Transform}\label{sec:decision criterion}

To understand the decision criterion of misspecification aversion in our newsvendor context, we investigate the objective function of \ref{prob:misspecification} from a perspective of {\it distributional transform}. We show that the objective function essentially transforms distributions in $\mathcal{A}$---via a {\it transform function} determined by the index $\alpha$ of misspecification aversion and the newsvendor's profit function---to new ones that possibly violate the mean and variance constraints specified by $\mathcal{A}$. 

Exploring the definition of $d(F,\mathcal{A})$ and interchanging the minimization over $F$ and $G$, we can represent equivalently the \ref{prob:misspecification} problem as follows:
\[
\max_{q \geq 0} \; \min_{G \in \mathcal{A}} \; \min_{F \in \mathcal{P}} \; \Big\{\mathbb{E}_{F}[\pi(q,\tilde{u})]+\alpha\cdot d_{\rm OT}(F,G)\Big\}.
\]
Here, the inner minimization features the potential misspecification of a fixed distribution $G \in \mathcal{A}$, which is then robustified over the specified ambiguity set $\mathcal{A}$ via the outer minimization. Recall that a distributional transform $T_{\varphi}[\cdot]: \mathcal{P} \mapsto \mathcal{P}_0$ maps distributions in $\mathcal{P}$ to $\mathcal{P}_0$ via a \textit{transform function} $\varphi$ (\citealt{liu2021distributional}). We can then represent the innermost minimization above through a distributional transform, leading to the following result that plays a key role in characterizing the decision criterion of \ref{prob:misspecification}. Different from \cite{liu2021distributional}, the transform function we identify further encodes operational information such as price parameter~$p$ and order quantity~$q$.

\begin{theorem}[{\sc Distributional Transform}]\label{thm:distributional transform}
Given $\alpha \geq 0$ and $q \geq 0$, it holds that
\[
\min_{F \in \mathcal{P}} \; \Big\{\mathbb{E}_{F}[\pi(q,\tilde{u})] + \alpha \cdot d(F,\mathcal{A})\Big\} = \min_{G \in \mathcal{A}} \; \int_{\mathbb{R}_+} \pi(q,v) \; {\rm d} T_{\varphi_\alpha}[G](v),
\]
where $T_{\varphi_\alpha}[G](v)=G\circ \varphi_\alpha^{-1}(v) ~\forall v \in \mathbb{R}_+$, and $T_{\varphi_\alpha}[\cdot]$ is a distributional transform of $G\in\mathcal{A}$ with an increasing and continuous transform function $\varphi_\alpha: \mathbb{R}_+ \mapsto \mathbb{R}_+$ defined as follows.
\begin{enumerate}
\item[(i)] If $\alpha < \frac{p}{4q}$, then $\varphi_\alpha(v) = \frac{\alpha}{p} \cdot v^2$. 

\item[(ii)] If $\alpha \geq \frac{p}{4q}$, then $\varphi_\alpha(v) =  \frac{\alpha}{p} \cdot v^2$ when $v < \frac{p}{2\alpha}$ and $\varphi_\alpha(v) = v-\frac{p}{4\alpha}$ otherwise.
\end{enumerate}
\end{theorem}

 \begin{figure}[t!]
\centering
\includegraphics[scale=0.2]{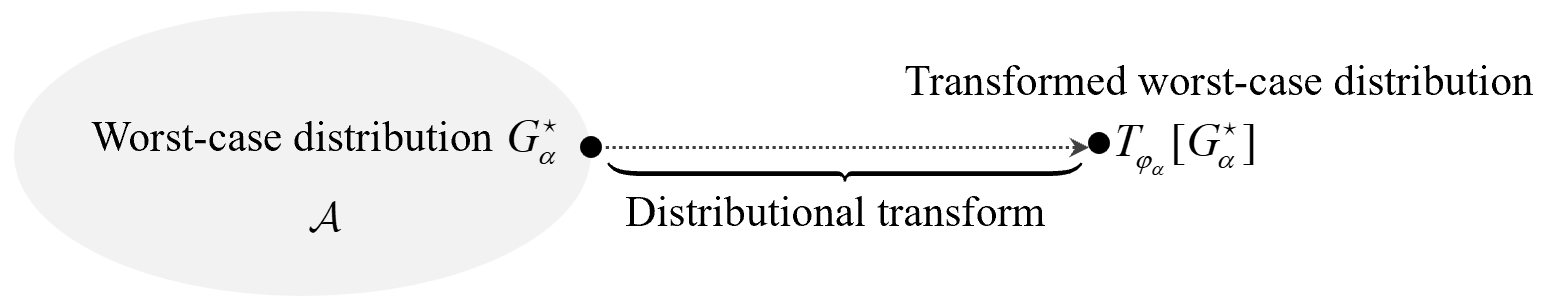}
\caption{\textnormal{The decision criterion of~\ref{prob:misspecification}  transforms the worst-case distribution $G_\alpha^\star$ in the ambiguity set $\mathcal{A}$ into the transformed worst-case distribution $T_{\varphi_\alpha}[G_\alpha^\star]$ that can be outside of $\mathcal{A}$.}}\label{fig:transform}
\end{figure}

By \Cref{thm:distributional transform}, the decision criterion of~\ref{prob:misspecification} serves as an expectation under a transformed worst-case distribution (Figure~\ref{fig:transform}), which establishes the equivalence between \ref{prob:misspecification} and the following problem
\begin{equation}\label{prob:transform}\tag{\sc Transform}
\max_{q \geq 0} \; \min_{G \in \mathcal{A}} \; \mathbb{E}_{T_{\varphi_\alpha}[G]}[\pi(q,\tilde{v})].
\end{equation}
For the remainder of this section, we may use \ref{prob:transform} and \ref{prob:misspecification} interchangeably. 

\begin{figure}[tb]
\begin{subfigure}{.49\textwidth}
\centering
\includegraphics[width=0.85\linewidth]{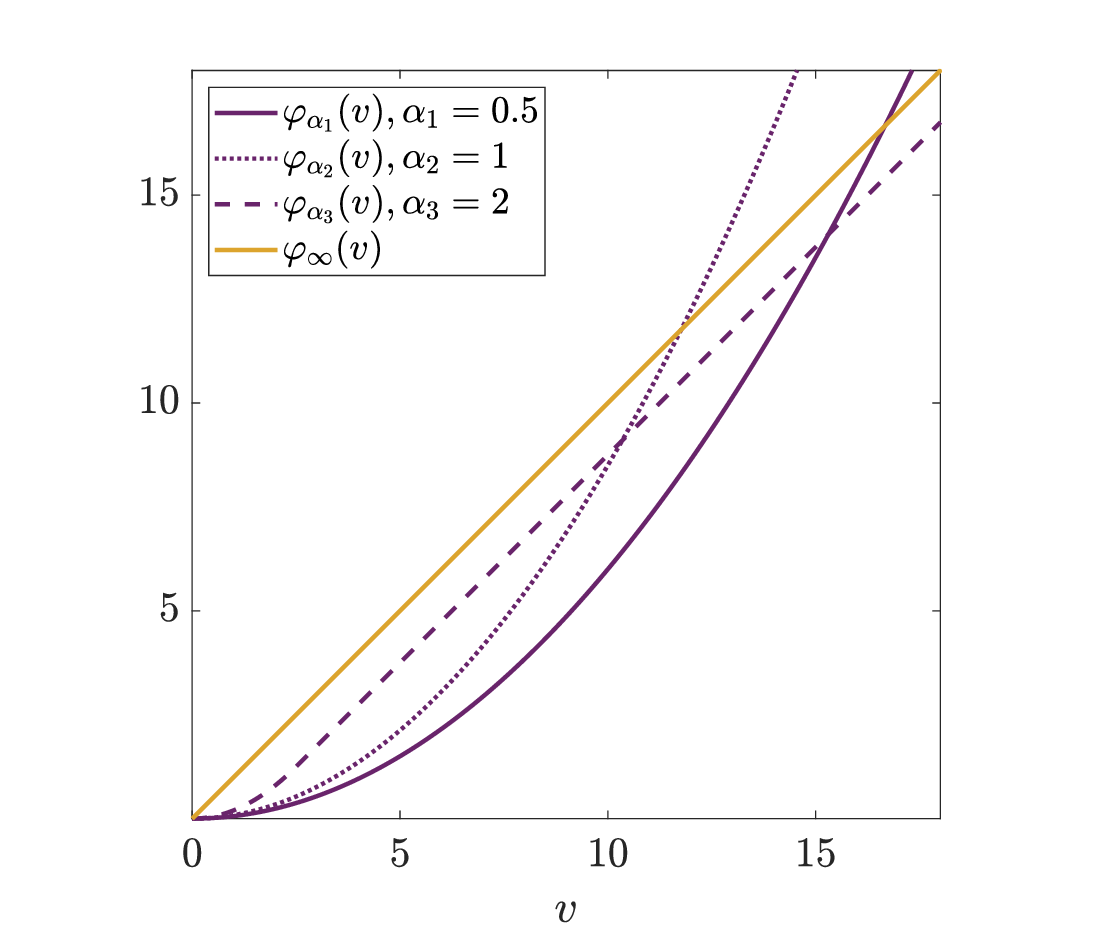}
\end{subfigure}
\hspace{-5mm}
\begin{subfigure}{.49\textwidth}
\centering
\includegraphics[width=0.85\linewidth]{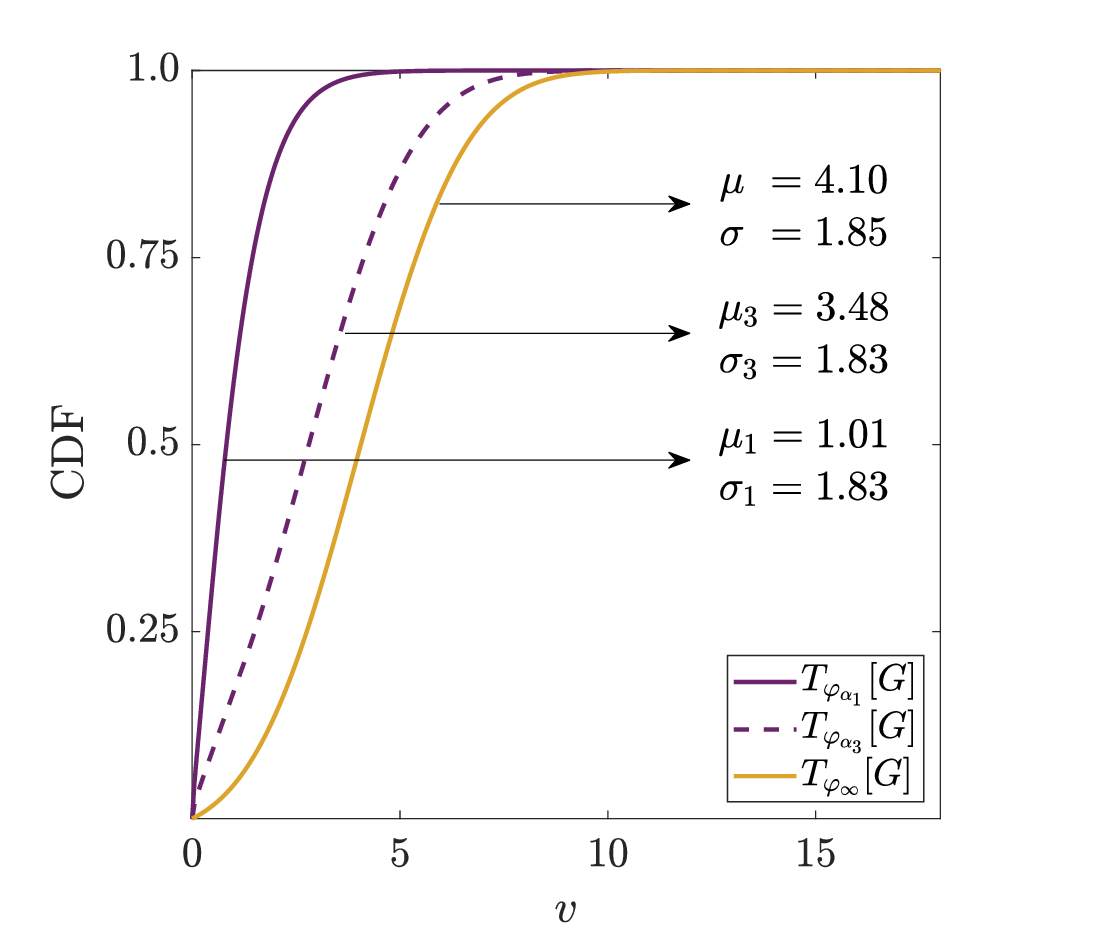}
\end{subfigure}
\caption{\textnormal{{\it Left}: Transform function $\varphi_\alpha(v)$. {\it Right}: CDF of a truncated normal distribution $G$ and the transformed distribution $T_{\varphi_{\alpha}}[G]$. On both panels,  $\alpha_1 < \alpha_2 < \frac{p}{4q} < \alpha_3$. For the right panel, $G \in \mathcal{A}$ is a normal distribution truncated to $\mathbb{R}_+$ with mean $\mu$ and standard deviation $\sigma$. Here, $\mu_1$ and $\sigma_1$ (resp., $\mu_3$ and $\sigma_3$) are the mean and  standard deviation of the transformed distribution $T_{\varphi_{\alpha_1}}[G]$ (resp., $T_{\varphi_{\alpha_3}}[G]$).}}
\label{fig:distortion}
\vspace{-4mm}
\end{figure}

Recall that \ref{prob:ambiguity} evaluates the performance of an order quantity $q$ via the decision criterion $\min_{G \in \mathcal{A}} \; \mathbb{E}_{G}[\pi(q,\tilde{v})]$.
The above equivalence reveals that the key {\it difference} between \ref{prob:ambiguity} and \ref{prob:misspecification} lies in the transform function $\varphi_\alpha$ applied to the probability distributions in the mean-variance ambiguity set $\mathcal{A}$. In particular, for any \textit{physical}\footnotemark\footnotetext{All distributions in $\mathcal{A}$ sharing the same \textit{physically} observable mean-variance information are treated indifferently.} probability distribution $G \in \mathcal{A}$, $\varphi_\alpha$ transforms it to a probability distribution $T_{\varphi_\alpha}[G]$ that can be {\it outside} $\mathcal{A}$, 
so that the resulting transformed expectation $\mathbb{E}_{T_{\varphi_\alpha}[G]}[\cdot]$ reflects the newsvendor’s aversion to misspecification. On the one hand, for a small value of $\alpha$ such that $\alpha < \frac{p}{4q}$, the transform function $\varphi_\alpha$ compresses (resp., amplifies) low (resp., high) demand realizations of $G$,\footnotemark\footnotetext{When $\alpha < \frac{p}{4q}$, $T_{\varphi_\alpha}[G] (v) = G(\sqrt{\frac{p}{\alpha}v})$ for every $G \in \mathcal{A}$. Since $G(v)$ is increasing in $v$, it can be seen that $T_{\varphi_\alpha}[G] (v) > G(v)$ when $v<\frac{p}{\alpha}$ and $T_{\varphi_\alpha}[G] (v) \leq G(v)$ when $v \geq \frac{p}{\alpha}$.} and a smaller $\alpha$ results in more demand realizations being compressed, see the cases of $\alpha_1, \alpha_2$ on the left panel of~\Cref{fig:distortion}.  On the other hand, for a large value of $\alpha$ such that $\alpha \geq \frac{p}{4q}$, $\varphi_\alpha$ compresses all demand realizations of $G$;\footnotemark\footnotetext{When $\alpha \geq \frac{p}{4q}$, $T_{\varphi_\alpha}[G] (v) = G(\sqrt{\frac{p}{\alpha}v}) \geq G(v)$ if $v < \frac{p}{2\alpha}$ and $T_{\varphi_\alpha}[G] (v) = G(v + \frac{p}{4\alpha}) \geq G(v)$ otherwise. Hence, for every $G \in \mathcal{A}$, we always have $T_{\varphi_\alpha}[G] (v) \geq G(v)$ for all $v \geq 0$.} see the case of $\alpha_3$ on the left panel of~\Cref{fig:distortion}. In this case, misspecification aversion compresses all probability distributions in $\mathcal{A}$, and a smaller $\alpha$ also leads to a stronger compression. Importantly, due to the transform function $\varphi_\alpha$, the mean or variance of the transformed distribution $T_{\varphi_\alpha}[G]$ can be different from that of the original distribution $G$ in the ambiguity set $\mathcal{A}$, namely $T_{\varphi_\alpha}[G] \notin \mathcal{A}$; see the right panel of \Cref{fig:distortion} for a visualization. As $\alpha \to \infty$, the transform function becomes $\varphi_\alpha(v)=v$,  implying $T_{\varphi_\alpha}[G]=G$ and that \ref{prob:misspecification} reduces to \ref{prob:ambiguity}. That is, the index $\alpha$ of misspecification aversion, or equivalently, the newsvendor's aversion against misspecification of $\mathcal{A}$, is fully encoded in the transform function $\varphi_\alpha$ of \ref{prob:transform}. The transformed distribution $T_{\varphi_\alpha}[G]$ could be far away from the ambiguity set $\mathcal{A}$, especially when $\alpha$ is small such that the mean of $T_{\varphi_\alpha}[G]$ is significantly lower than that of distributions in $\mathcal{A}$ (see the case of $\alpha_1$ on the right panel of~\Cref{fig:distortion}).

Interestingly, we emphasize that given $\alpha$ and an order quantity $q$, the distributional transform $T_{\varphi_\alpha}$ is determined by the price $p$ (see the transform function $\varphi_\alpha$ derived in \Cref{thm:distributional transform}), leading to the transformed worst-case distribution $T_{\varphi_\alpha}[G_\alpha^\star]$ with $G_\alpha^\star\in\argmin_{G\in\mathcal{A}} \;  \mathbb{E}_{T_{\varphi_\alpha}[G]}[\pi(q,\tilde{v})]$ being dependent on the price (see~\Cref{prop:worst-case transformed distribution} in~\Cref{appendix:worst-case transformed distribution}). In contrast, we recall that the worst-case distribution implied by \ref{prob:ambiguity}, \textit{i.e.}, $G_\infty^\star\in\argmin_{G\in\mathcal{A}} \;  \mathbb{E}_{G}[\pi(q,\tilde{v})]$, is however independent of the cost structure.\footnotemark\footnotetext{Given an order quantity $q$, it can be shown that 
$G^\star_\infty =(\frac{\sigma^2}{\mu^2 + \sigma^2}) \cdot \delta_0 + (\frac{\mu^2}{\mu^2 + \sigma^2}) \cdot \delta_{\frac{\mu^2 + \sigma^2}{\mu}}$ when $q < \frac{\mu^2 + \sigma^2}{2\mu}$, and $G^\star_\infty =(\frac{1}{2} + \frac{q-\mu}{2w}) \cdot \delta_{q - w} + (\frac{1}{2} - \frac{q-\mu}{2w}) \cdot \delta_{q + w}$ otherwise, where $w=\sqrt{(q-\mu)^2 + \sigma^2}$, which is independent of $p$.} In other words, the \textit{price-independent} worst-case distribution  \textit{inside} the ambiguity set $\mathcal{A}$ in the decision criterion  of~\ref{prob:ambiguity} now becomes a \textit{price-dependent} transformed worst-case distribution \textit{outside} of $\mathcal{A}$ in the decision criterion of~\ref{prob:misspecification}. We emphasize that such price-dependency effect indeed leads to operational consequences of~\ref{prob:misspecification} being distinct from that of~\ref{prob:ambiguity}, which will be further explored in Section~\ref{sec:optimal solution}.

\section{Optimal Solution and Sensitivity Analysis}\label{sec:optimal solution}

In this section, we first derive the optimal order quantity of~\ref{prob:misspecification} in closed form. We then investigate the impact of misspecification aversion via the optimal order quantity's sensitivity to the cost-structure information and distributional information, revealing important operational implications of \ref{prob:misspecification} distinct from~\ref{prob:ambiguity}.

\subsection{Analytical Solution}

To proceed, we recall that \ref{prob:misspecification} can be reformulated as
\[
\max_{q \geq 0} \; \min_{F \in \mathcal{P}} \; \{\mathbb{E}_{F}[\pi(q,\tilde{u})] + \alpha \cdot d(F,\mathcal{A})\} 
= \max_{q \geq 0} \; \min_{G \in \mathcal{A}} \; \mathbb{E}_{T_{\varphi_\alpha}[G]}[\pi(q,\tilde{v})] 
= \max_{q \geq 0} \; \min_{G \in \mathcal{A}} \; \mathbb{E}_{G}[\pi(q,\varphi_\alpha(\tilde{v}))],
\]
where the second equality follows from the definition of distributional transform $T_{\varphi_\alpha}$ as identified in~\Cref{thm:distributional transform}. Hence, one can tackle~\ref{prob:misspecification} by adapting the primal-dual machinery for solving a maximin problem similar to \ref{prob:ambiguity} but with a new ``profit" function $\Psi(\alpha,q,v) = \pi(q,\varphi_\alpha(v))$. Specifically, given $q \geq 0$, the primal and dual reformulation of the inner worst-case expectation can be written as
\begin{equation}\label{prob:primal misspecification}
\begin{array}{cll}
\displaystyle \min_{G\in \mathcal{M}_+} & \displaystyle \int_{\mathbb{R}_+} \; \Psi(\alpha, q, v) \; {\rm d}G(v) \\[3mm]
{\rm s.t.} & \displaystyle \int_{\mathbb{R}_+} \; v \; {\rm d}G(v) = \mu &~~\cdots~s_\alpha \\[3mm]
& \displaystyle \int_{\mathbb{R}_+}\; v^2 \; {\rm d}G(v) = \mu^2 + \sigma^2 &~~\cdots~r_\alpha  \\[3mm]
& \displaystyle \int_{\mathbb{R}_+}\;{\rm d}G(v) = 1 &~~\cdots~t_\alpha,
\end{array}
~~~\Longleftrightarrow~~~
\begin{array}{cll}
\displaystyle \max_{s_\alpha, r_\alpha, t_\alpha} & \mu s_\alpha - (\mu^2+\sigma^2) r_\alpha - t_\alpha \\[1mm]
{\rm s.t.} & v s_\alpha - v^2 r_\alpha - t_\alpha \leq \Psi(\alpha, q, v) &~\forall v \geq 0 \\
& s_\alpha \in \mathbb{R}, \; r_\alpha \in \mathbb{R}, \; t_\alpha \in \mathbb{R}.
\end{array}
\end{equation}
where $s_\alpha$, $r_\alpha$, and $t_\alpha$ are dual variables associated with the mean, variance, and support of the ambiguity set $\mathcal{A}$, respectively.
The key to the primal-dual machinery is to construct a pair of primal and dual solutions that share identical objective values. In particular, the primal solution is a worst-case distribution constructed by identifying tangent points between the function $s_\alpha v - r_\alpha v^2 - t_\alpha$ and the function $\Psi(\alpha, q, v)$ in the dual reformulation. For \ref{prob:misspecification}, however, the new ``profit'' function $\Psi$---neither convex nor concave in $v$---is less structured than the concave piecewise affine function $\pi$ of \ref{prob:ambiguity}, making the primal-dual procedure more involved to analyze. Fortunately, leveraging the closed form of $\Psi(\alpha,q,v) = \pi(q,\varphi_\alpha(v))$ given by~\Cref{thm:distributional transform}, we can derive an analytical reformulation of the objective function in \ref{prob:misspecification}.

\begin{proposition}[{\sc Worst-Case Transformed Expectation}]\label{prop:transformed expectation}
Given $\alpha \geq 0$ and $q \geq 0$, 
\begin{equation}\label{equ:transformed-expectation}
\mspace{-20mu}
\min_{G \in \mathcal{A}} \; \mathbb{E}_{T_{\varphi_\alpha}[G]}[\pi(q,\tilde{v})]= \begin{cases}
\displaystyle \frac{p}{2}\bigg(q+\mu - \frac{p}{4\alpha} - \sqrt{\bigg(q-\mu + \frac{p}{4\alpha}\bigg)^2 + \sigma^2}\bigg) - cq & \displaystyle {\rm if}~ q \in \mathcal{Q}\\[5mm]
\displaystyle \frac{\alpha}{2}\bigg(\frac{pq}{\alpha}+\mu^2+\sigma^2-\sqrt{\bigg(\frac{pq}{\alpha}+\mu^2+\sigma^2\bigg)^2-4\mu^2\frac{pq}{\alpha}}\bigg) - cq & {\rm otherwise},
\end{cases}    
\end{equation}
where $\mathcal{Q}= \{q \in \mathbb{R}_+ \mid q \geq \frac{p}{4\alpha}, \;(2\mu-\frac{p}{\alpha})q \geq \mu^2+\sigma^2 - \frac{p\mu}{2\alpha}\}$.
\end{proposition}

The analytical form of the worst-case transformed expectation is non-trivial and generalizes the worst-case expected cost of \ref{prob:ambiguity} as $\alpha \to \infty$. Indeed, as $\alpha \to \infty$, $\mathcal{Q}$ becomes $\{q \in \mathbb{R}_+ \mid q \geq \frac{\mu^2 + \sigma^2}{2\mu}\}$ and~\eqref{equ:transformed-expectation} recovers the worst-case expected cost of~\ref{prob:ambiguity} such that $\min_{G \in \mathcal{A}} \mathbb{E}_{G}[\pi(q,\tilde{v})] = \frac{\mu^2 pq}{\mu^2 + \sigma^2} - cq$ if $q \ge \frac{\mu^2 + \sigma^2}{2\mu}$ and 
$\min_{G \in \mathcal{A}} \mathbb{E}_{G}[\pi(q,\tilde{v})] = \frac{p}{2}(q + \mu - \sqrt{(q - \mu)^2 + \sigma^2}) - cq$ otherwise. Equipped with the analytical form, we can then derive the optimal solution of~\ref{prob:misspecification} as follows. 

\begin{theorem}[{\sc Optimal Solution}]\label{thm:newsvendor transport-2}
Given $\alpha\geq 0$, the optimal order quantity $q^\star_\alpha$ of \ref{prob:misspecification} is
\begin{equation}\label{equ:quantity misspecification}
q_\alpha^\star=\left\{
\begin{array}{ll}
\displaystyle \mu + \sigma f(1 - \kappa) - \frac{p}{4\alpha} &~~~\displaystyle \kappa \geq \frac{\sigma^2}{\mu^2+\sigma^2}, \; \displaystyle \alpha \geq \frac{p}{2(\mu - \sigma \sqrt{(1-\kappa)/\kappa})} \\ [3.5mm]
\displaystyle \big(\mu^2-\sigma^2 + 2\mu\sigma f(1 - \kappa)\big) \cdot \frac{\alpha}{p} &~~~\displaystyle \kappa \geq \frac{\sigma^2}{\mu^2+\sigma^2}, \; \displaystyle \alpha < \frac{p}{2(\mu - \sigma \sqrt{(1-\kappa)/\kappa})}\\ [3.5mm]
0 &~~~\displaystyle \kappa < \frac{\sigma^2}{\mu^2+\sigma^2},
\end{array}
\right.
\end{equation}
where $f(\cdot)$ is defined in \eqref{equ:optimal order ambiguity}. The optimal order quantity $q^\star_\alpha$ is increasing in $\alpha$.
\end{theorem}

\begin{figure}[tb]
\begin{subfigure}{0.33\textwidth}
\centering
\includegraphics[width=1.11\linewidth]{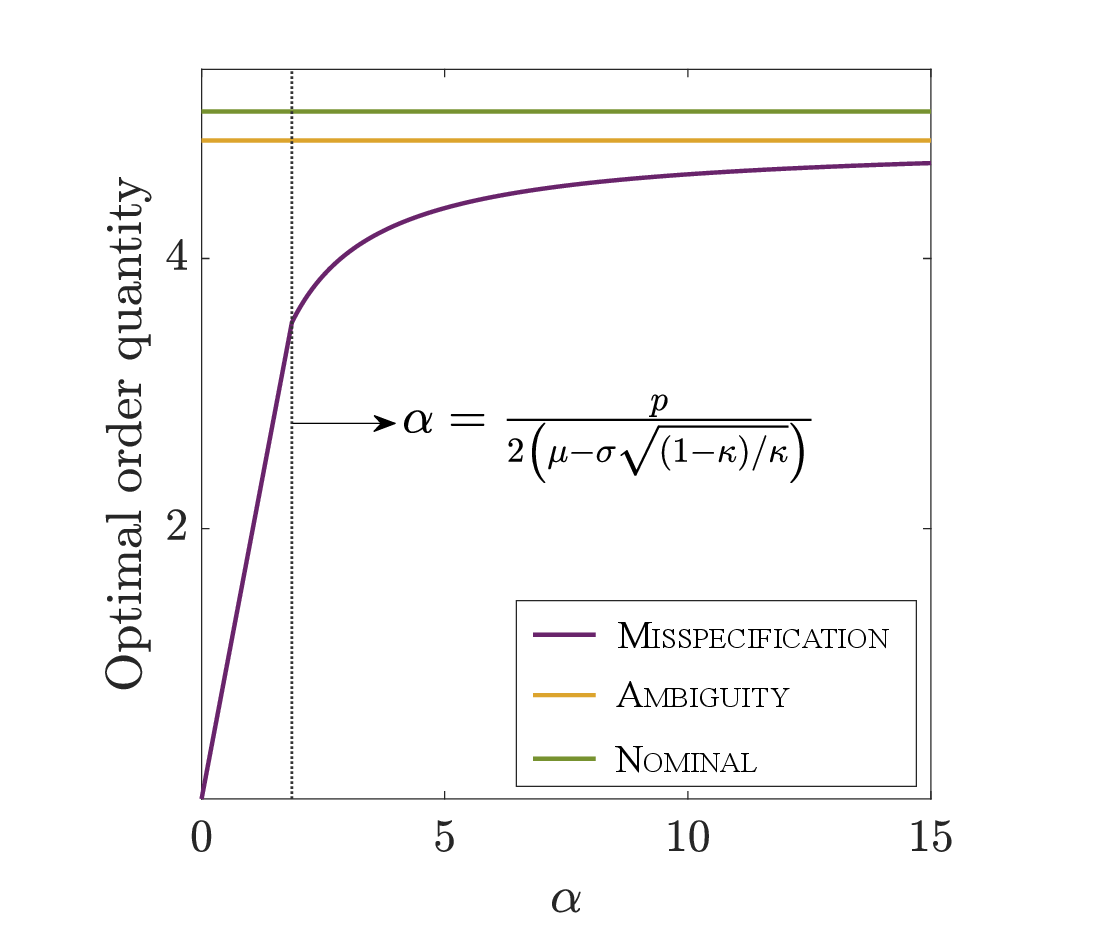}
\end{subfigure}
\begin{subfigure}{0.33\textwidth}
\centering
\includegraphics[width=1.11\linewidth]{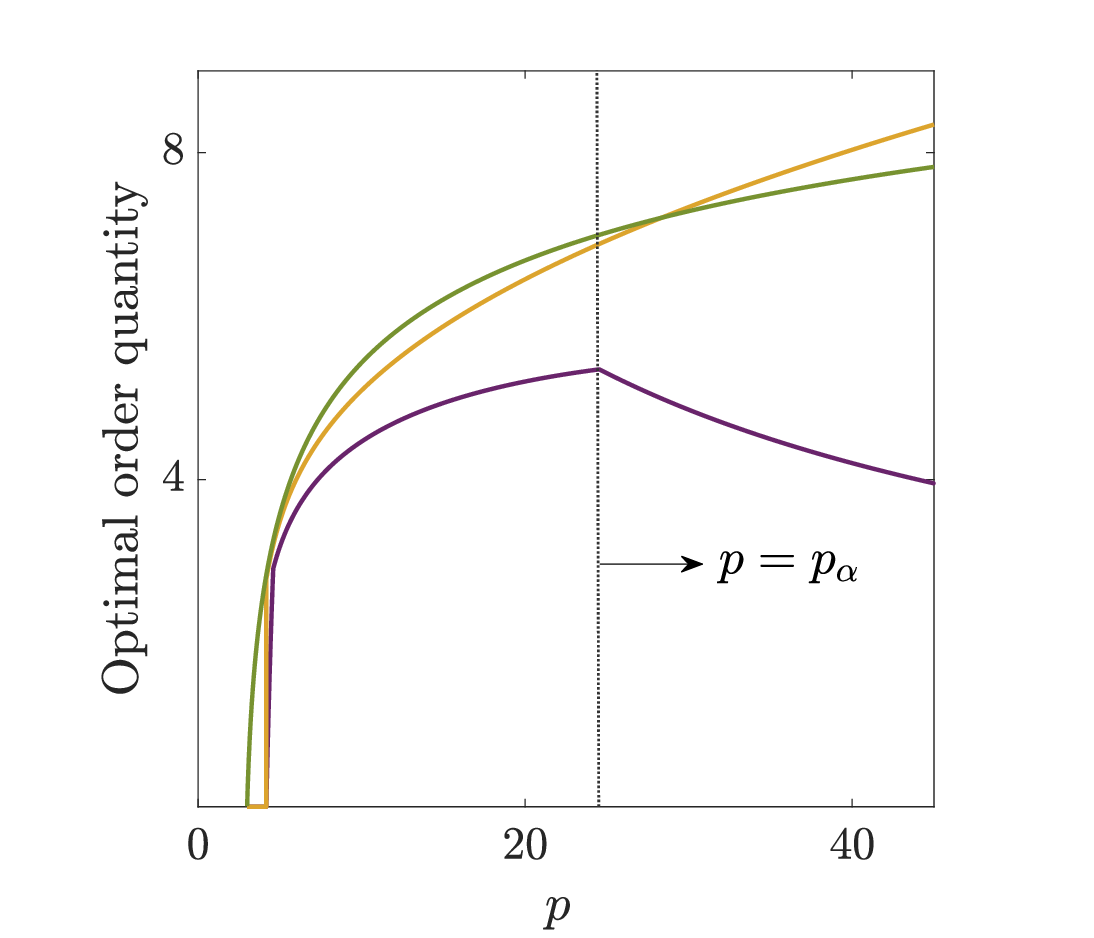}
\end{subfigure}
\begin{subfigure}{0.33\textwidth}
\centering
\includegraphics[width=1.11\linewidth]{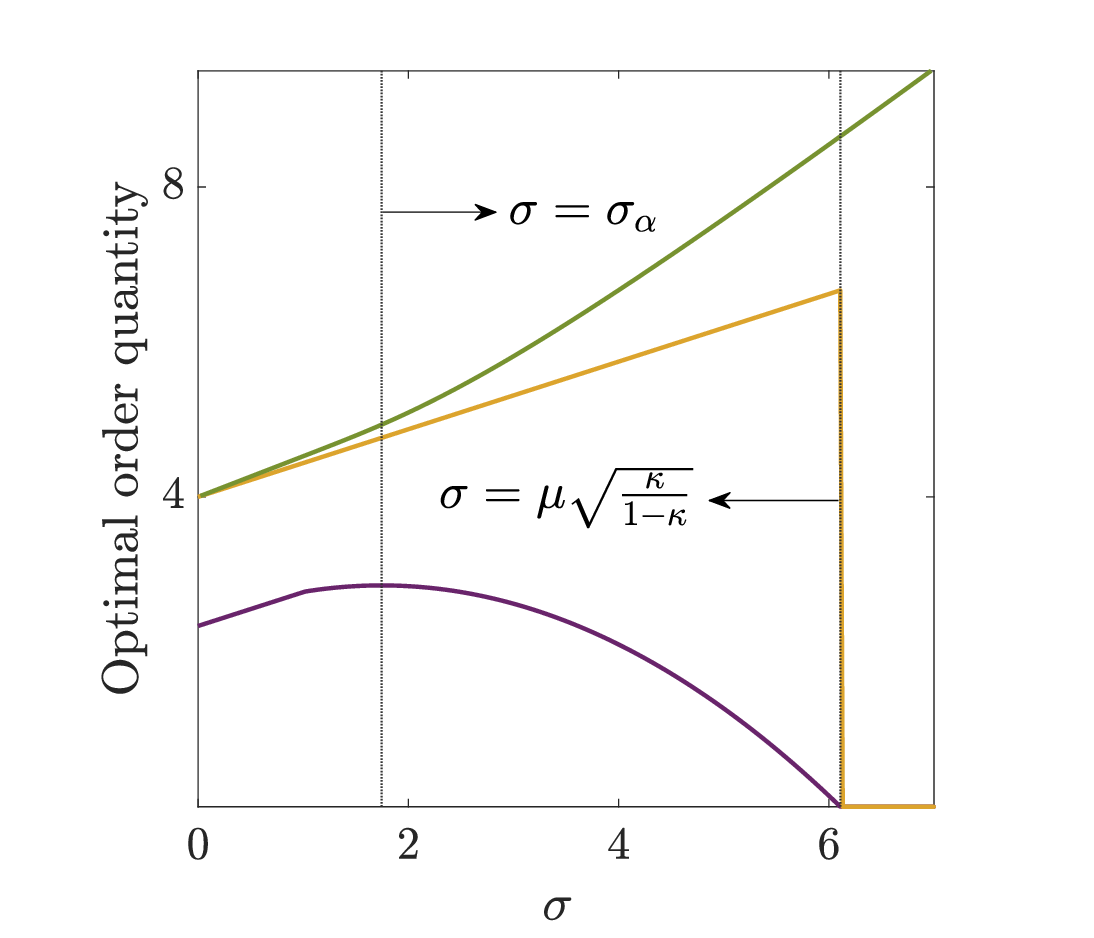}
\end{subfigure}
\caption{\textnormal{Optimal order quantity as a function of $\alpha$ (left), $p$ (middle), and $\sigma$ (right), respectively.}}
\footnotesize{{\it Notes.} $\;$ The optimal order quantity of \ref{prob:nominal} is obtained under a normal distribution truncated to $\mathbb{R}_+$ with mean $\mu$ and standard deviation $\sigma$. For all three panels, we set $c = 3$ and $\mu = 4$. On the left panel, $p=10$ and $\sigma = 2$. On the middle panel, $\sigma=2.5$ and $\alpha=4$, and we identify $p_\alpha = 22.5$. On the right panel, $p=10$ and $\alpha = 1.5$, for which $\kappa=0.7$ and the non-degenerate region is $\textstyle\sigma \in \big[0, \mu\sqrt{\frac{\kappa}{1-\kappa}}\big]=\big[0,4\sqrt{\frac{7}{3}}\big]$, and we identify $\sigma_\alpha = \frac{8}{\sqrt{21}}$.}
\label{fig:optimal order quantity type-2}
\end{figure}

Focusing on the non-degenerate case that $\kappa \geq \frac{\sigma^2}{\mu^2+\sigma^2}$, for $0\le \alpha_1 < \frac{p}{2(\mu - \sigma \sqrt{(1-\kappa)/\kappa})}\le \alpha_2$, we have
\[
q_{\alpha_1}^\star=\big(\mu^2-\sigma^2 + 2\mu\sigma f(1 - \kappa)\big) \cdot \frac{\alpha_1}{p} \le q_{\alpha_2}^\star = \mu + \sigma f(1 - \kappa) - \frac{p}{4\alpha_2} \le \mu + \sigma f(1 - \kappa)=q_{\infty}^\star;
\]
see left panel of \Cref{fig:optimal order quantity type-2}. This implies that the optimal order quantity $q^\star_\alpha$ of \ref{prob:misspecification} is no larger than that of~\ref{prob:ambiguity} (\textit{i.e.}, $q^\star_\infty$)---an intuitive result due to the additional aversion to misspecification---and $q^\star_\alpha \to q^\star_\infty$ as $\alpha \to \infty$. It is also notable that the optimal order quantity $q^\star_\alpha$ is affected by misspecification aversion and ambiguity aversion {\it separately}. When $\alpha < \frac{p}{2(\mu - \sigma \sqrt{(1-\kappa)/\kappa})}$, $q_\alpha^\star$ is a product of $\mu^2-\sigma^2+2\mu\sigma f(1-\kappa)$ and $\alpha/p$ that are purely determined by the mean and variance information specified in $\mathcal{A}$ and purely determined by misspecification, respectively. When $\alpha \geq \frac{p}{2(\mu - \sigma \sqrt{(1-\kappa)/\kappa})}$, $q_\alpha^\star$ is obtained by $\mu + \sigma f(1-\kappa)$ that is exactly the optimal order quantity $q^\star_\infty$ of \ref{prob:ambiguity} minus $\frac{p}{4\alpha}$ that is purely determined by misspecification.

\subsection{Sensitivity and its Implications}\label{sec:sensitivity}

We next look at the optimal order quantity $q_\alpha^\star$'s sensitivity to the cost-structure information (\textit{i.e.}, $c$ and $p$) and distributional information (\textit{i.e.}, $\mu$ and $\sigma^2$). Since it is straightforward that $q_\alpha^\star$ is decreasing (resp., increasing) in $c$ (resp., $\mu$), being consistent with that of~\ref{prob:ambiguity}, we focus on its sensitivity to price $p$ and variance $\sigma^2$, which exhibits a different pattern from that of~\ref{prob:ambiguity}. In particular, we also focus on the non-degenerate case in \Cref{thm:newsvendor transport-2} that $\kappa \geq \frac{\sigma^2}{\mu^2+\sigma^2}$. In this case, recall from~\eqref{equ:optimal order ambiguity} that the optimal order quantity $q_\infty^\star$ of \ref{prob:ambiguity} is
\begin{equation}\label{equ:optimal order ambiguity non-degenerate}
q_\infty^\star = \mu + \sigma\frac{2\kappa-1}{2\sqrt{\kappa(1-\kappa)}} = \mu + \frac{\sigma}{2}\bigg(\sqrt{\frac{\kappa}{1-\kappa}} - \sqrt{\frac{1-\kappa}{\kappa}}\bigg).
\end{equation}
By~\eqref{equ:optimal order ambiguity non-degenerate}, $q^\star_\infty$ of~\ref{prob:ambiguity} is always \textit{increasing} in $p$, so is $q_G^\star = G^{-1}(\kappa)$ of~\ref{prob:nominal} under any distribution $G$, by noting that the profit margin $\kappa = \frac{p-c}{p}$ is increasing in $p$. However, the monotonicity of $q_\alpha^\star$ to $p$ can reverse that of $q_\infty^\star$ ($q_G^\star$) as shown below.

\begin{proposition}[{\sc Sensitivity to Price}]\label{prop:price}
There exists some $p_\alpha \geq \max\big\{\frac{\mu^2+\sigma^2}{\mu^2}c, 2\alpha\mu\big\}$ such that $q_\alpha^\star$ is decreasing in $p \in (p_\alpha, \infty)$.
\end{proposition}

Notably, \Cref{prop:price} points out that when price $p$ is sufficiently large, the optimal order quantity $q_\alpha^\star$ of \ref{prob:misspecification}, in stark contrast to $q_\infty^\star$ of \ref{prob:ambiguity}, is \textit{decreasing} in $p$;  see the middle panel of \Cref{fig:optimal order quantity type-2}. We emphasize that such distinct sensitivity is rooted from the induced \textit{price dependency} of the transformed worst-case demand distribution $T_{\varphi_\alpha}[G_\alpha^\star]$ in \ref{prob:misspecification}, which stochastically reduces in the price $p$ (see transform function $\varphi_\alpha$ in \Cref{thm:distributional transform}), while the worst-case demand distribution of \ref{prob:ambiguity} is \textit{independent} of the price. As the price becomes sufficiently large, the effect of reduced ``demand" ($T_{\varphi_\alpha}[G_\alpha^\star]$) outweighs that of the increased profit margin, leading to the $q_\alpha^\star$ of \ref{prob:misspecification} being decreased.

Likewise, by~\eqref{equ:optimal order ambiguity non-degenerate}, the optimal order quantity $q^\star_\infty$ of~\ref{prob:ambiguity} is \textit{increasing} (resp., decreasing) in $\sigma$ when $\kappa \geq \frac{1}{2}$ (resp., when $\kappa < \frac{1}{2}$), under the non-degenerate condition $\kappa \geq \frac{\sigma^2}{\mu^2+\sigma^2}$, \textit{i.e.}, $\textstyle\sigma \in \big[0, \mu\sqrt{\frac{\kappa}{1-\kappa}}\big]$. Such monotonicity also holds for the optimal order quantity $q_G^\star = G^{-1}(\kappa)$ of~\ref{prob:nominal} under many commonly used distributions (\textit{e.g.}, elliptical, uniform, and exponential). However, the situation becomes different when considering the misspecification aversion: the optimal order quantity $q^\star_\alpha$ of~\ref{prob:misspecification}, can be \textit{decreasing} in $\sigma$ when $\kappa \geq \frac{1}{2}$ in the situation of non-degeneracy; see the right panel of \Cref{fig:optimal order quantity type-2}. The different sensitivity pattern also uncovers an advantage of \ref{prob:misspecification} in the solution's smoothness to parameters, by noting that the solution to \ref{prob:ambiguity} is overly sensitive to the parameter $\sigma$ as it can jump as $\sigma$ changes slightly; see, \textit{e.g.}, \cite{embrechts2022robustness}, for various smoothness issues in optimizing risk measures.

\begin{proposition}[{\sc Sensitivity to Variance}]\label{prop:variance}
Given $\kappa \geq \frac{1}{2}$, there exists some $\sigma_\alpha \leq \mu\sqrt{\frac{\kappa}{1-\kappa}}$ such that $q_\alpha^\star$ is decreasing in $\textstyle\sigma \in \big[\sigma_\alpha, \mu\sqrt{\frac{\kappa}{1-\kappa}}\big]$. 
\end{proposition}

To summarize, sharing the mean and variance characteristics, the optimal order quantity of~\ref{prob:ambiguity}---hedging against the distributional uncertainty \textit{within} the ambiguity set $\mathcal{A}$---may exhibit an \textit{identical} sensitivity pattern to the cost-structure parameters with that of~\ref{prob:nominal}, and, in the situation of non-degeneracy, to the distributional characteristics with that of~\ref{prob:nominal} under many common distributions inside $\mathcal{A}$. For instance, both solutions of~\ref{prob:ambiguity} and~\ref{prob:nominal} share the same monotonicity to the cost $c$ and to the mean value $\mu$. The optimal order quantity of~\ref{prob:misspecification}, however, hedges against another layer of distributional uncertainty \textit{beyond} the ambiguity set $\mathcal{A}$, which therefore can \textit{break} the sensitivity pattern of ordering characterized by the mean and variance information. This suggests that the ambiguity and misspecification, as different layers of distributional uncertainty, could result in \textit{distinct} operational consequences, and therefore should be distinguished in the modeling.

\section{Performance Guarantee}\label{sec:shadow value}

In this section, we investigate the out-of-sample performance guarantee of the optimal order quantity $q_\alpha^\star$ of~\ref{prob:misspecification}. As we have mentioned, in many practical situations, the newsvendor has only access to incomplete knowledge on demand, and the misspecification can arise from a mixing effect of estimation error (\textit{e.g.}, due to data limitation) and distribution shift (\textit{e.g.}, due to non-stationarity). We consider a data-driven setting where the mean-variance ambiguity set $\mathcal{A}$ is estimated as $\mathcal{A}_N$ by using demand samples drawn from a data-generating distribution $D$ with finite mean $\mu$ and standard deviation $\sigma$. 

\begin{assumption}\label{assump:iid}
Random samples $\hat{v}_1,\ldots,\hat{v}_N$ are independently drawn from the data-generating distribution $D$, and the mean-variance ambiguity set $
\mathcal{A}_N = \{G \in \mathcal{P} \mid \mathbb{E}_{G}[\tilde{v}] = \hat{\mu},\;\mathbb{E}_{G}[\tilde{v}^2] = \hat{\mu}^2+\hat{\sigma}^2\}
$
is constructed from sample mean $\hat{\mu} = \frac{1}{N}\sum_{i = 1}^N \hat{v}_i > 0$ and variance $\hat{\sigma}^2 = \frac{1}{N}\sum_{i = 1}^N \hat{v}_i^2 - \big(\frac{1}{N}\sum_{i =1}^N \hat{v}_i\big)^2 > 0$.
\end{assumption}

We consider the possibility that the out-of-sample distribution $F$ can be {\it different} from $D$---a phenomenon of distribution shift. In particular, we examine the finite-sample performance guarantee through a statistical approach that provides an insightful interpretation of the performance guarantee of \ref{prob:misspecification} by decoupling the effects of estimation error and distribution shift. 

Our analysis relies on the concentration of the estimated mean-variance ambiguity set $\mathcal{A}_N$, for which we need to investigate the optimal-transport cost $d(D, \mathcal{A}_N)$ that is closely related to the \textit{Gelbrich distance} \citep{gelbrich1990formula}. For any $G \in \mathcal{A}_N$, the Gelbrich distance between $G$ and the data-generating $D$ is $\sqrt{(\hat{\mu} - \mu)^2 + (\hat{\sigma} - \sigma)^2}$, and  $d_{\rm OT}(D,G) \geq (\hat{\mu} - \mu)^2 + (\hat{\sigma} - \sigma)^2$ with the inequality being tight whenever $G$ is an affine transformation of $D$ (as shown in \citealt{gelbrich1990formula, nguyen2021mean}).
That is to say, if $\mathcal{A}_N$ is supported on the whole space $\mathbb{R}$, then $d(D, \mathcal{A}_N) = (\hat{\mu} - \mu)^2 + (\hat{\sigma} - \sigma)^2$---the optimal-transport cost amounts to the Gelbrich distance squared. This is \textit{not} necessarily true in our newsvendor context as the demand is non-negative and $\mathcal{A}_N$ should be supported on $\mathbb{R}_+$. Quite notably, we show that the optimal-transport cost still coincides with the Gelbrich distance squared if $\frac{\hat{\mu}}{\hat{\sigma}} \geq \frac{\mu}{\sigma}$ and otherwise, is bounded from above by the Gelbrich distance squared plus a term related to the true and estimated mean and variance---a result may be of independent interest.

\begin{lemma}\label{lemma:distance}
Under Assumption~\ref{assump:iid},  the optimal-transport cost of moving the data-generating distribution $D$ to the mean-variance ambiguity set $\mathcal{A}_N$ can be characterized as follows.
\begin{enumerate}
\item[(i)] If $\frac{\hat{\mu}}{\hat{\sigma}} \geq \frac{\mu}{\sigma}$, then $d(D, \mathcal{A}_N) = (\hat{\mu} - \mu)^2 + (\hat{\sigma} - \sigma)^2$.

\item[(ii)] If $\frac{\hat{\mu}}{\hat{\sigma}} < \frac{\mu}{\sigma}$, then for sufficiently large $N$, it holds that
\[
(\hat{\mu} - \mu)^2 + (\hat{\sigma} - \sigma)^2 \leq d(D,\mathcal{A}_N) \leq (\hat{\mu} - \mu)^2 + (\hat{\sigma} - \sigma)^2 + \frac{\mu^2\hat{\sigma}^2 - \hat{\mu}^2\sigma^2}{\sigma\hat{\sigma}}.
\]
\end{enumerate}
\end{lemma}

We also assume the following regularity condition on the data-generating distribution $D$.

\begin{assumption}\label{assump:sub-Gaussian}
The data-generating distribution $D$ is sub-Gaussian with a variance proxy $\nu^2$, \textit{i.e.}, $\mathbb{E}_D[{\rm exp}(x(\tilde{v}-\mu))] \leq {\rm exp}(\frac{x^2\nu^2}{2}),~\forall x \in \mathbb{R}$.
\end{assumption}

Sub-Gaussianity, as a common type of light-tailed characteristics,\footnotemark\footnotetext{The light-tailed assumption is typically necessary for establishing the large-deviation properties for statistics of the sample mean and sample variance. As for the heavy-tailed distributions, more complex estimation procedures (for the mean and variance) are needed to achieve acceptable convergence rates (\citealp{cai2010optimal}).} captures a wide range of distributions, including, among many others, the Gaussian distribution, the Bernoulli distribution, the uniform distribution on a convex set, and \textit{any} bounded distributions (\citealp{vershynin2010introduction}). With the characterization of $d(D, \mathcal{A}_N)$ in \Cref{lemma:distance}, we derive the following concentration inequality.

\begin{proposition}[\sc{Concentration of Mean-Variance Ambiguity Set}]\label{prop:concentration}
Under Assumptions~\ref{assump:iid} and~\ref{assump:sub-Gaussian}, for a given confidence level $\eta\in(0,1]$, it holds for sufficiently large $N$ that 
\[
\mathbb{P}_{D^N}\bigg[d(D, \mathcal{A}_N) \leq \frac{(c_1 + c_2 \log(1/\eta))^2}{\sqrt{N}}\bigg] \geq 1 - \eta,
\]
where $c_1, c_2>0$ are constants that only depend on $\mu$, $\sigma$, and $\nu$.
\end{proposition}

Leveraging the concentration of ambiguity set $\mathcal{A}_N$ and the closed-form expression of the worst-case transformed expectation (\Cref{prop:transformed expectation}), we can then establish a finite-sample performance guarantee of the optimal solution to \ref{prob:misspecification}. 

\begin{theorem}[\sc{Finite-Sample Performance Guarantee}]\label{thm:guarantee}
Under Assumptions~\ref{assump:iid} and~\ref{assump:sub-Gaussian}, for a given confidence level $\eta\in(0,1]$, let $\varepsilon_N = \frac{(c_1 + c_2 \log(1/\eta))^2}{\sqrt{N}}$ with $c_1, c_2>0$ being constants that only depend on $\mu$, $\sigma$, and $\nu$. Let
\begin{equation}\label{equ:alpha-N}
\textstyle \alpha_N =
\frac{1}{2}\sqrt{\frac{p(p-c)}{\varepsilon_N + d_{\rm OT}(F,D)}},
\end{equation}
if $\varepsilon_N + d_{\rm OT}(F,D) < \kappa(\hat{\mu} - \hat{\sigma}\sqrt{(1 - \kappa)/\kappa})^2$; and $\alpha_N = 0$, otherwise. Consider the optimal solution $q^\star_{\alpha_N}$ and the optimal value $\Pi_{\alpha_N}^\star$ of \ref{prob:misspecification} with $\alpha=\alpha_N$ and $\mathcal{A} = \mathcal{A}_N$. For sufficiently large $N$, it holds that 
\begin{equation*}
\displaystyle \mathbb{P}_{D^N}\bigg[\mathbb{E}_{F}[\pi(q^\star_{\alpha_N},\tilde{u})] \geq \Big(\hspace{-2mm}\underbrace{\Pi_{\alpha_N}^\star}_{\substack{\textnormal{in-sample}\\ \textnormal{optimal value}}} \hspace{-2mm} - \hspace{3mm} \frac{1}{2}\sqrtexplained{\hspace{-2mm}\underbrace{p(p-c)\varepsilon_N\vphantom{\Pi_{\alpha_N}^\star}}_{\substack{\textnormal{effect of}\\[0.6mm] \textnormal{estimation error}}} \hspace{-1mm}+\hspace{0.5mm}\underbrace{p(p-c)d_{\rm OT}(F,D)\vphantom{\Pi_{\alpha_N}^\star}}_{\substack{\textnormal{effect of}\\[0.6mm] \textnormal{distribution shift}}}}\Big)^+\bigg] \geq 1 - \eta,    
\end{equation*}
where $D^N$ is the $N$-fold product of $D$. 
\end{theorem}

The guarantee derived in \Cref{thm:guarantee} is the in-sample optimal value of \ref{prob:misspecification} subtracting the out-of-sample misspecification effect described by the estimation error in mean and variance and the distribution shift. For the estimation error, it is related to the upper bound $\varepsilon_N$ on the distance $d(D,\mathcal{A}_N)$ between the data-generating distribution and the estimated mean-variance ambiguity set, which diminishes as $N\to\infty$. For the distribution shift captured by $d_{\rm OT}(F,D)$, it is \textit{independent} of the sample size $N$. A key statistical implication is that as long as the out-of-sample distribution $F$ shifts from the data-generating distribution $D$, there is always a constant amount of loss $\frac{1}{2}\sqrt{p(p-c)d_{\rm OT}(F,D)}$ in terms of the performance guarantee, even as $N \to\infty~(\varepsilon_N \to 0)$.

\Cref{thm:guarantee} also suggests that the calibration for the index $\alpha_N$ of misspecification aversion is affected by both the estimation error and the distribution shift. According to~\eqref{equ:alpha-N}, when non-zero, $\alpha_N$ is increasing in the sample size $N$ (\textit{i.e.}, decreasing in $\varepsilon_N$) while decreasing in the extent of distribution shift $d_{\rm OT}(F,D)$. This implies that at the same confidence level, the newsvendor needs to be more misspecification averse in either case of a smaller amount of data or a more significant distribution shift, to guarantee the performance. Moreover, in the presence of distribution shift ($d_{\rm OT}(F,D)>0$), even when the estimation error vanishes with the sufficient data (\textit{i.e.}, $\varepsilon_N \to 0$ as $N \to \infty$), $\alpha_N \to \frac{1}{2}\sqrt{\frac{p(p-c)}{d_{\rm OT}(F,D)}}<\infty$, implying that \ref{prob:misspecification} ($\alpha_N<\infty$) shall still outperform \ref{prob:ambiguity} ($\alpha_N=\infty$). Importantly, this means that the commonly used cross-validation approach (which is based on the data-generating distribution $D$) for calibrating distributionally robust optimization models could, unfortunately, work {\it poorly} in the situation of misspecification (\textit{i.e.}, calibrating the index $\alpha$ for \ref{prob:misspecification}), as we demonstrate in the numerical study. 

\section{Extensions}\label{extension}

In this section, we extend the model \ref{prob:misspecification} along the following directions: (\textit{i}) there are multiple products, (\textit{ii}) the ambiguity set is defined via optimal transport, and (\textit{iii}) the extent of misspecification is measured by total variation distance.

\subsection{Multiple Products}\label{sec:multi product}

Our misspecification-averse model can be extended to the multi-product newsvendor problem. Consider $M$ products (each with unit price $p_i$ and cost $c_i$, $i \in [M]$) whose random demands are collectively denoted by $\tilde{\bm{u}} = (\tilde{u}_1, \dots, \tilde{u}_M)$ that follows a multi-dimensional distribution $F$. The misspecification-averse newsvendor then solves
\begin{equation}\label{prob:multi product}\tag{\sc Multiple}
\max_{\bm{q} \geq \bm{0}} \; \min_{F \in \mathcal{P}_M} \; \{\mathbb{E}_{F}[\omega(\bm{q},\tilde{\bm{u}})] + \alpha \cdot d(F,\mathcal{C})\},
\end{equation}
where $\mathcal{P}_M$ is the set of probability distributions supported on $\mathbb{R}^M_+$, the optimal-transport cost $d_{\rm OT}(\cdot,\cdot)$ is defined in \eqref{equ:optimal transport 2} with the cost function $\|\cdot\|_2^2$, $\bm{q} = (q_1,\dots,q_M)$ is the vector of order quantities, and $\omega(\bm{q}, \bm{u}) = \sum_{i = 1}^M \pi_i(q_i,u_i)= \sum_{i = 1}^M p_i \cdot \min\{q_i, u_i\} - c_i q_i$.

If the ambiguity set $\mathcal{C}$ is specified by marginal mean-variance information of each product, then \ref{prob:multi product} is separable concerning products and \Cref{thm:newsvendor transport-2} yields each product's optimal order quantity. If $\mathcal{C}$ is specified by mean and correlation information, then \ref{prob:multi product} becomes much more involved, as its ambiguity-averse counterpart is already intractable~(\citealt{hanasusanto2015distributionally, natarajan2018asymmetry}). More details on the reformulation and computational difficulty of \ref{prob:multi product} with complete covariance information are presented in \Cref{appendix:multi product}. In the following, we show that an analytical solution can be derived for a case that captures the partial correlation across products. In particular, we consider an ambiguity set with mean and sum-of-variance constraints:
\[
\mathcal{C} = \bigg\{G \in \mathcal{P}_M 
~\bigg|~ 
\mathbb{E}_G[\tilde{v}_i] = \mu_i ~\forall i \in [M],\;\;\mathbb{E}_G\bigg[\sum_{i \in [M]}\tilde{v}_i^2\bigg] \leq K
 \bigg\},
\]
where $K$ is some non-negative constant that bounds the sum of the variance of products' demands. Note that when $M = 1$ (\textit{i.e.}, there is a single product), the ambiguity set $\mathcal{C}$ reduces to the mean-variance ambiguity set $\mathcal{A}$. The analytical solution to \ref{prob:multi product} with $\mathcal{C}$ relies on characterizing the optimal dual variable, denoted by $\lambda^\star$, to its sum-of-variance constraint. We indicate that adopting a similar technique, we can also obtain the optimal solution of \ref{prob:multi product} under an additional budget constraint $\sum_{i\in[M]}q_i \leq Q$ for some $Q > 0$.

To ease our presentation, let $\bar{\lambda}_0 = 0$, $\bar{\lambda}_i = \frac{c_i}{(2\mu_i - p_i/\alpha)^+}$ for $i \in [M]$ and $\bar{\lambda}_{M+1} = +\infty$. Without loss of generality, we rearrange $\bar{\lambda}_i, i\in [M]$ ascendingly, {\it i.e.}, $\bar{\lambda}_1 \leq \bar{\lambda}_2 \leq \ldots \leq \bar{\lambda}_M$. Besides, we define
\begin{equation}\label{equ:index}
i^\star = \min\{j \in[M+1] \mid \Theta_j(\bar{\lambda}_j) < K\},
\end{equation}
where for $j \in [M+1]$ and $\lambda \geq 0$,
\begin{equation}\label{equ:Theta}
\Theta_j(\lambda) = \sum_{i \in [M] \backslash [j-1]} \frac{p_i\mu_i^2(\alpha^2 c_i + 2\alpha c_i \lambda + p_i \lambda^2)}{(p_i \lambda + \alpha c_i)^2} + \sum_{i \in [j-1]} \frac{(p_i - c_i)c_i}{4\lambda^2}.
\end{equation}
We then present the following strong-duality result that characterizes the optimal dual variable to the sum-of-variance constraint and the \textit{decomposability} of the dual problem.

\begin{theorem}\label{thm:multi product dual}
Given $K \geq 0$ and $\alpha \geq 0$, \ref{prob:multi product} is equivalent to
\begin{equation*}
\max_{\lambda \geq 0} \; \bigg\{- \lambda K + \sum_{i \in [M]} \varPi_{i, \alpha}^\star(\lambda) \bigg\}.
\end{equation*}
Here, for each $i \in [M]$, $\varPi_{i, \alpha}^\star(\lambda)$ is the optimal value of the following optimization problem:
\begin{equation}\label{equ:single product misspecification}
\varPi_{i, \alpha}^\star(\lambda)=\max_{q_i \geq 0} \; \min_{F_i \in \mathcal{P}} \; \bigg\{\mathbb{E}_{F_i}[\pi_i(q_i,\tilde{u}_i)] + \min_{G_i\in\mathcal{C}_i} \; \big\{\lambda\cdot\mathbb{E}_{G_i}[\tilde{v}_i^2] + \alpha\cdot d_{\rm OT}(F_i,G_i)\big\}\bigg\}
\end{equation}
with $\tilde{u}_i \sim F_i$, $\tilde{v}_i \sim G_i$, and $\mathcal{C}_i = \{G \in \mathcal{P} \mid \mathbb{E}_{G}[\tilde{v}] = \mu_i\}$ being a mean ambiguity set. The optimal $\lambda^\star$ is decreasing in $K$, and with $i^\star$ in \eqref{equ:index} and $\Theta_i(\cdot)$ in \eqref{equ:Theta}, it can be characterized as follows.
\begin{enumerate}
\item[(i)] If $\Theta_{i^\star}(\bar{\lambda}_{i^\star - 1}) \leq K$, then $\lambda^\star = \bar{\lambda}_{i^\star - 1}$.
\vspace{1mm}
\item[(ii)] If $\Theta_{i^\star}(\bar{\lambda}_{i^\star - 1}) > K$, then $\lambda^\star\in(\bar{\lambda}_{i^\star - 1}, \bar{\lambda}_{i^\star})$ is the solution to the equation $\Theta_{i^\star}(\lambda) = K$.
\end{enumerate}
\end{theorem}

By \Cref{thm:multi product dual}, the value of $\lambda^\star$ can be efficiently binary searched due to monotonicity. Given the value of $\lambda^\star$, \ref{prob:multi product} can be decomposed into multiple misspecification-averse {\it single-product} newsvendor problems in \eqref{equ:single product misspecification}.
Moreover, it is critical to note that the optimal dual variable $\lambda^\star$ captures not only the demand-correlation information encoded in the sum-of-variance constraint, but also the cross-product cost structure (\textit{i.e.}, $p_i$ and $c_i$, $i \in [M]$) that collectively affects those single-product problems. The following result adapts the reasoning for \Cref{thm:newsvendor transport-2} to derive the analytical solution for \ref{prob:multi product}.

\begin{theorem}[\sc Optimal Solution: Multiple Products]\label{thm:multi product}
Given $\alpha \geq 0$ and the optimal $\lambda^\star$ characterized in \Cref{thm:multi product dual}, the optimal order quantity $\bm{q}^\star_\alpha$ of \ref{prob:multi product} is, for each $i \in [M]$, $q_{i,\alpha}^\star = \mu_i + \frac{p_i - 2c_i}{4\lambda^\star} - \frac{p_i}{4\alpha}$ if $\alpha \geq \frac{p_i}{2(\mu_i - c_i/(2\lambda^\star))^+}$ and $q_{i,\alpha}^\star = \frac{\lambda^\star (\lambda^\star + \alpha) p_i \mu_i^2}{\alpha(p_i\lambda^\star/\alpha + c_i)^2}$ otherwise.
\end{theorem}

The analytical solution $\bm{q}_\alpha^\star$ captures the information of distribution characteristics, cost structure, and correlation across products. In particular, the optimal order quantity $q_{i,\alpha}^\star$ for each product $i$ is not only determined by its distribution characteristic (\textit{i.e.}, $\mu_i$) and cost structure (\textit{i.e}, $p_i$ and $c_i$), but also by those of other products via $\lambda^\star$. Furthermore, in the multi-product problem, the index of misspecification aversion $\alpha$ has a \textit{double effect} on the optimal order quantity in two ways: on the one hand, it directly affects the optimal order quantity $q_{i,\alpha}^\star$ of each product in the decomposed single-product problem~\eqref{equ:single product misspecification}; on the other hand, it also affects the optimal dual variable $\lambda^\star$ that in turn influences the optimal order quantities $q_{i,\alpha}^\star$ of all products.

Finally, we emphasize that \Cref{thm:multi product} generalizes \Cref{thm:newsvendor transport-2} from a single product to multiple products (see \Cref{EC:generalization} for a detailed derivation) as well as the forthcoming~\Cref{coro:multi product ambiguity} for the ambiguity-averse multi-product model to ambiguity and misspecification aversion. Note that as $\alpha \to \infty$, \ref{prob:multi product} immediately reduces to the ambiguity-averse only counterpart
\begin{equation}\label{prob:multi product ambiguity}
\max_{\bm{q} \geq \bm{0}} \; \min_{G \in \mathcal{C}} \; \mathbb{E}_{G}[\omega(\bm{q},\tilde{\bm{v}})].
\end{equation}
Its optimal solution, as expected, generalizes the Scarf model from a single product to multiple products. The proof is straightforward and is thus omitted.
\begin{corollary}\label{coro:multi product ambiguity}
Let the optimal dual variable $\lambda^\star$ be characterized in \Cref{thm:multi product dual}, the optimal quantity $\bm{q}^\star_\infty$ of  the ambiguity-averse multi-product newsvendor problem~\eqref{prob:multi product ambiguity} is, for each $i \in [M]$, $q_{i,\infty}^\star = \mu_i + \frac{p_i - c_i}{4\lambda^\star}$ if $\lambda^\star \geq \frac{c_i}{2\mu_i}$
and $q_{i,\infty}^\star = \frac{p_i\mu_i^2}{c_i^2}\cdot \lambda^\star$ otherwise.
\end{corollary}

\subsection{Distance-Based Ambiguity Set}\label{sec:ambiguity set}

Apart from the mean-variance ambiguity set $\mathcal{A}$, the distance-based ambiguity set $\mathcal{B}(\theta) = \{G \in \mathcal{P} \mid d_{\rm OT}(G,H) \leq \theta\}$
with a reference distribution $H$, $\theta \geq 0$, and optimal-transport cost $d_{\rm OT}(\cdot,\cdot)$ between probability distributions is a popular alternative for specifying partial distributional information, which has also been widely used in the aforementioned applications of decision theory (\citealt{petracou2022decision}), newsvendor (\citealt{chen2021regret}, \citealt{zhang2023optimal}), and risk management (\citealt{wozabal2014robustifying}). Conceptually, $\mathcal{B}(\theta)$ consists of all probability distributions in a $\theta$-neighbourhood around $H$, where the closeness is measured by $d_{\rm OT}(\cdot,\cdot)$ given in~\eqref{equ:optimal transport 2}. Since $\mathcal{B}(\theta_1) \subseteq \mathcal{B}(\theta_2)$ for any $\theta_2 \geq \theta_1 \geq 0$, a larger value of $\theta$ indicates a lower confidence in $H$. When $\theta = 0$, $\mathcal{B}(\theta)$ shrinks to a singleton containing only the reference distribution $H$, that is, $\mathcal{B}(0) = \{H\}$. It is natural to consider the following variant of~\ref{prob:misspecification} where we replace $\mathcal{A}$ with $\mathcal{B}(\theta)$:
\begin{equation}\label{prob:double misspecification}
\max_{q \geq 0} \; \min_{F \in \mathcal{P}} \; \Big\{\mathbb{E}_{F}[\pi(q,\tilde{u})] + \alpha \cdot d(F,\mathcal{B}(\theta))\Big\},
\end{equation}
which hedges against the possible misspecification over the ambiguity set $\mathcal{B}(\theta)$. Quite notably, problem~\eqref{prob:double misspecification} is essentially equivalent to hedging against misspecification over the singleton $\{H\}$ but with a {\it stronger} aversion to misspecification. To avoid a degenerate case, we assume $H^{-1}(\kappa) > 0$.

\begin{theorem}[\sc Optimal Solution: Distance-Based Ambiguity Set]\label{thm:double misspecification}
Given $\theta \geq 0$, $\alpha \geq 0$ and a reference distribution $H$, there exists some $\gamma^\star \in [0, \alpha]$ such that problem~\eqref{prob:double misspecification} is equivalent to 
\begin{equation}\label{prob:double misspecification equivalence}
\max_{\psi \geq 0} \; \min_{F \in \mathcal{P}} \; \Big\{\mathbb{E}_{F}[\pi(\psi,\tilde{u})] + \gamma^\star \cdot d_{\rm OT}(F,H)\Big\}.
\end{equation}
When $\theta = 0$, $\gamma^\star = \alpha$; otherwise, with the optimal order quantity $q^\star_H = H^{-1}(\kappa)$ of \ref{prob:nominal} under $H$ and $\beta = \int_{0}^{q^\star_H}u^2{\rm d}H(u) > 0$, $\gamma^\star$ can be characterized as follows.
\begin{enumerate}
\item[(\textit{i})] If $\theta \geq \beta$, then $\gamma^\star = 0$.
		
\item[(\textit{ii})] If $\theta < \beta$ and $\alpha(1-\sqrt{\theta/\beta}) < \frac{p}{2 q^\star_H}$, then $\gamma^\star = \alpha(1-\sqrt{\theta/\beta})$.
		
\item[(\textit{iii})] If $\theta < \beta$ and $\alpha(1-\sqrt{\theta/\beta}) \geq \frac{p}{2 q^\star_H}$, then $\gamma^\star$ is the solution to $\int_0^{\frac{p}{2x}}u^2{\rm d}H(u) + \frac{p^2}{4 x^2}\big(\kappa-H\big(\frac{p}{2 x}\big)\big) -\frac{\alpha^2\theta}{(\alpha - x)^2} =0$.
\end{enumerate}
With $\gamma^\star$, the optimal order quantity $\psi^\star_{\gamma^\star}$ of problem~\eqref{prob:double misspecification equivalence} can be characterized as $\psi^\star_{\gamma^\star} =  q^\star_H \cdot (\frac{\gamma^\star  q^\star_H}{p})$ if $\gamma^\star < \frac{p}{2 q^\star_H}$ and $\psi^\star_{\gamma^\star} = q^\star_H - \frac{p}{4 \gamma^\star}$ otherwise.
\end{theorem}

\Cref{thm:double misspecification} establishes the equivalence between problem~\eqref{prob:double misspecification}---which hedges against misspecification over a distance-based ambiguity set $\mathcal{B}(\theta)$ around the reference distribution---and problem~\eqref{prob:double misspecification equivalence} that, with a stronger aversion, hedges against misspecification to the reference distribution $H$. Note that \Cref{thm:double misspecification} states that $\gamma^\star = \alpha$ whenever $\theta = 0$---which, indeed, corresponds to the ambiguity neutrality. If we further have $\alpha \to \infty$, then problem~\eqref{prob:double misspecification} reduces to \ref{prob:nominal} under $H$, and as expected, \Cref{thm:double misspecification} concludes that $\psi^\star_{\gamma^\star} = q^\star_H$. 

\subsection{Misspecification Measured by Total Variation Distance}\label{sec:phi divergence}

Apart from the optimal-transport cost, total variation distance is also popular for measuring the closeness between probability distributions. In this section, we replace $d_{\rm OT}(\cdot,\cdot)$ in~\eqref{equ:optimal transport 2} with total variation distance when defining $d(F,\mathcal{A})$ and investigate the corresponding \ref{prob:misspecification} problem. Formally, the total variation distance between distributions $F$ and $G$ is defined as 
$d_{\rm TV}(F, G) = \sup_{A \in \mathcal{B}(\mathbb{R}_+)}|F(A) - G(A)|$,
where $\mathcal{B}(\mathbb{R}_+)$ is the Borel $\sigma$-algebra on $\mathbb{R}_+$. Equipped with the total variation distance,\footnotemark\footnotetext{Other types of $\phi$-divergence can also be applied to problem \eqref{prob:newsvendor TV} of misspecification, which, however, are less computationally appealing. For instance, if $d_{\rm OT}(\cdot, \cdot)$ is replaced with the Kullback–Leibler divergence, then problem \eqref{prob:newsvendor TV} becomes optimizing the worst-case constant absolute risk aversion, which is generally intractable (\citealp{chen2025robust}). Also, if we replace $d_{\rm OT}(\cdot, \cdot)$ with the Gini concentration index, then problem \eqref{prob:newsvendor TV} becomes optimizing worst-case mean-variance, which does not admit analytical solutions either.} we investigate the following variant of \ref{prob:misspecification}:
\begin{equation}\label{prob:newsvendor TV}
\max_{q\geq 0} \; \min_{G\in\mathcal{A}} \; \min_{F\in\mathcal{P}} \; \{\mathbb{E}_F[\pi(q,\tilde{v})] + \alpha \cdot d_{\rm TV}(F, G)\}.
\end{equation}

\begin{theorem}[\sc Optimal Solution: TV-Based Misspecification]\label{thm:newsvendor TV}
Given $\alpha\geq 0$ and the mean-variance ambiguity set $\mathcal{A}$, problem~\eqref{prob:newsvendor TV} can be equivalently reformulated as $\max_{0\leq q \leq \frac{2\alpha}{p}} \min_{G\in\mathcal{A}} \mathbb{E}_G[\pi(q, \tilde{v})]$,
and its optimal order quantity is $q_\alpha^\star = \min\{\frac{2\alpha}{p}, q_\infty^\star\}$,
where $q_\infty^\star$ is the optimal order quantity of \ref{prob:ambiguity} characterized in \eqref{equ:optimal order ambiguity}.
\end{theorem}

\Cref{thm:newsvendor TV} reveals that the optimal order quantity $q^\star_\alpha$ of problem~\eqref{prob:newsvendor TV} grows increasingly and linearly in $\alpha$, capped by the optimal order quantity $q^\star_\infty$ of~\ref{prob:ambiguity}.  When $\alpha$ is sufficiently large (that is, $\alpha \geq \frac{p}{2}q^\star_\infty$), $q_\alpha^\star$ coincides with $q^\star_\infty$. In this case, the newsvendor fully trusts the information specified in $\mathcal{A}$ and essentially becomes misspecification neutral to the ambiguity set $\mathcal{A}$. Hence, ambiguity aversion takes full charge of determining the formula of $q^\star_\alpha$. When $\alpha$ is relatively small (that is, $\alpha < \frac{p}{2}q^\star_\infty$), $q_\alpha^\star = \frac{2\alpha}{p}$ is then purely determined by misspecification aversion without being affected by the mean-variance information in $\mathcal{A}$.  

\section{Numerical Experiments with Retailing Data}\label{sec:experiments}

In this section, we demonstrate the effectiveness of incorporating misspecification, using the real-world daily demand data over one year for different stock keeping units (SKUs) of our industrial partner (a supermarket). Our goal is to compare the out-of-sample expected profit of the optimal order quantities obtained from \ref{prob:misspecification} and \ref{prob:ambiguity}, as well as~\ref{prob:nominal}.

\subsection{Data and Models}\label{sec:data and setting}

\begin{figure}[tb]
\begin{subfigure}{.33\textwidth}
\centering
\includegraphics[width=1.11\linewidth]{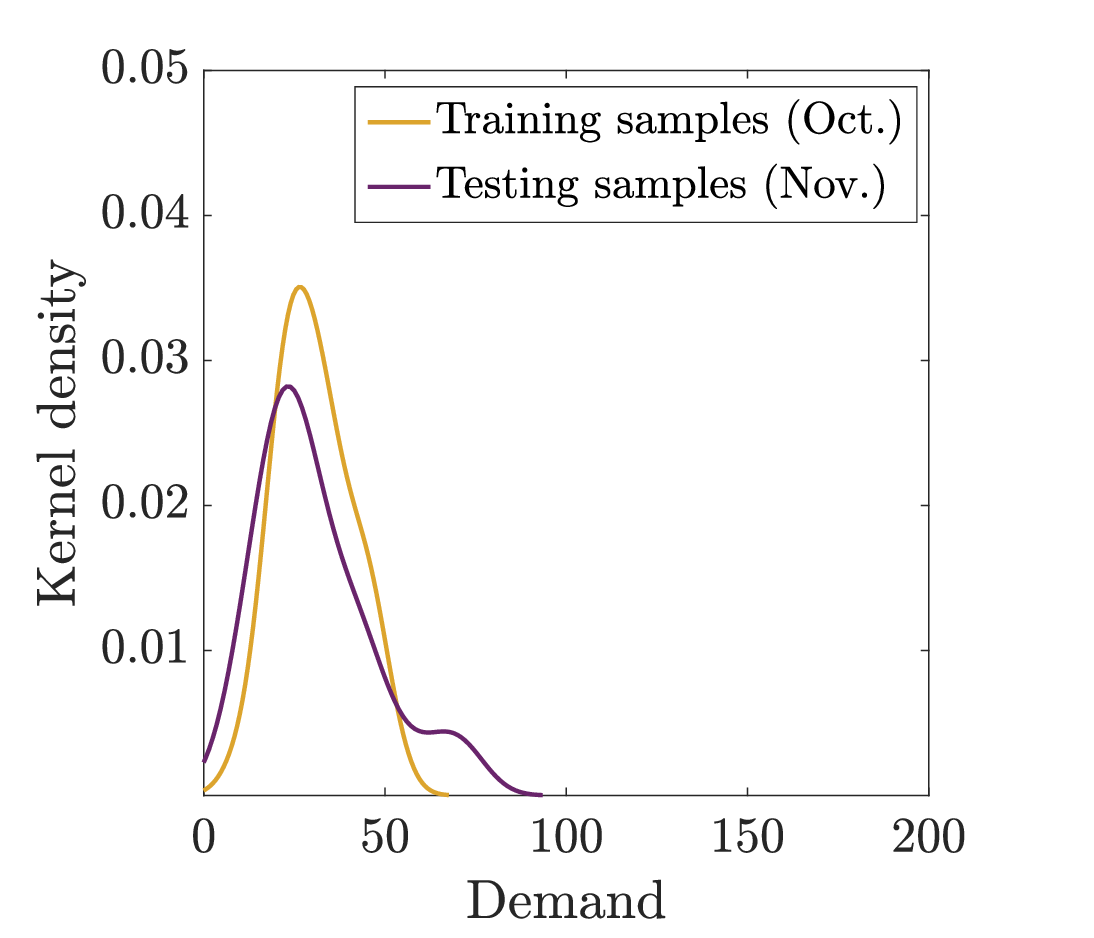}
\end{subfigure}
\begin{subfigure}{.33\textwidth}
\centering
\includegraphics[width=1.11\linewidth]{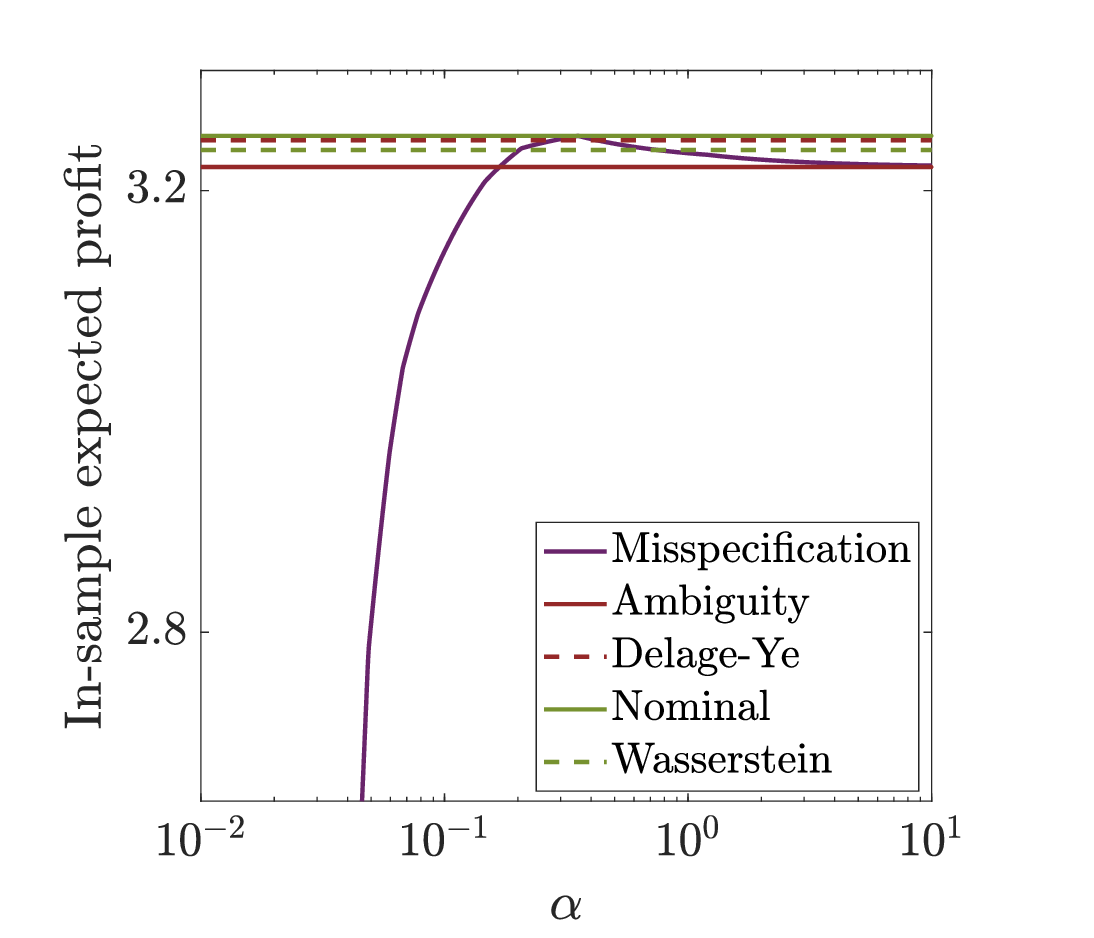}
\end{subfigure}
\begin{subfigure}{.33\textwidth}
\centering
\includegraphics[width=1.11\linewidth]{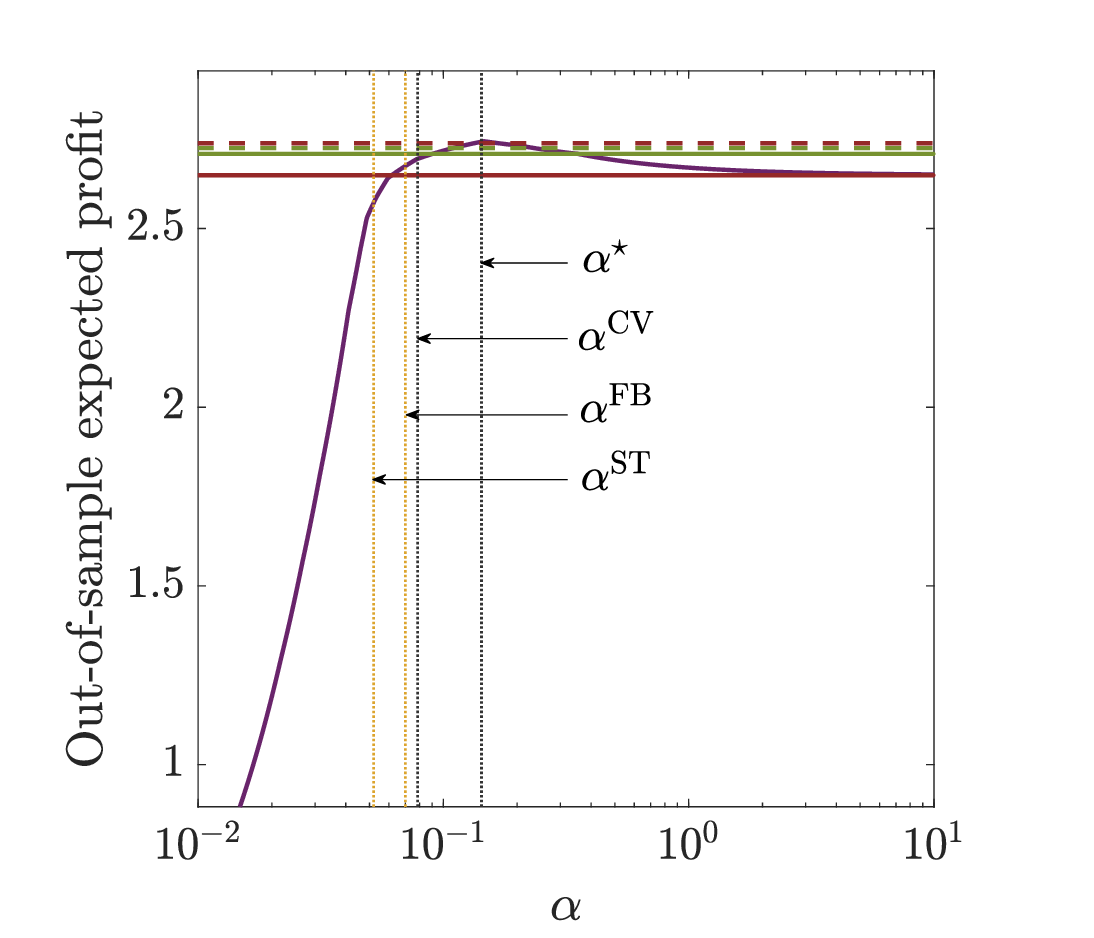}
\end{subfigure}
\begin{subfigure}{.33\textwidth}
\centering
\includegraphics[width=1.11\linewidth]{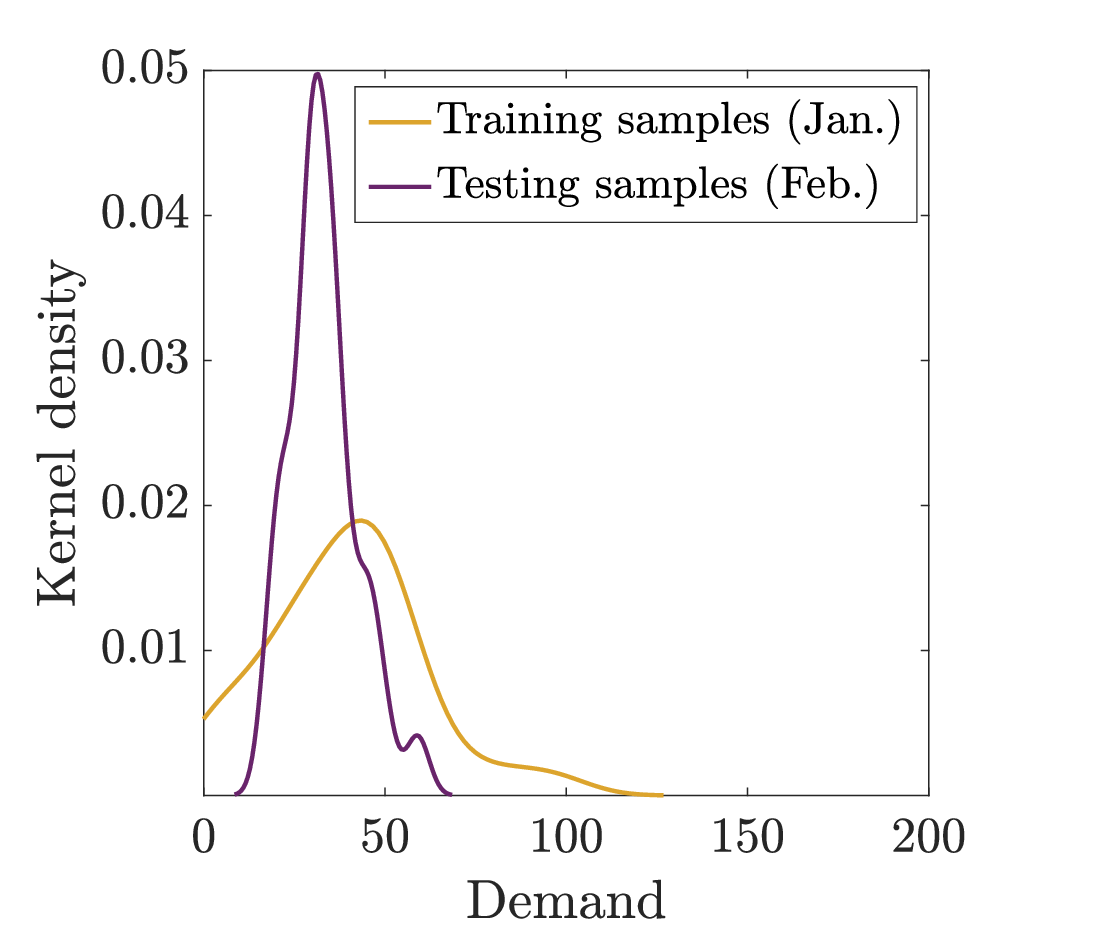}
\end{subfigure}
\begin{subfigure}{.33\textwidth}
\centering
\includegraphics[width=1.11\linewidth]{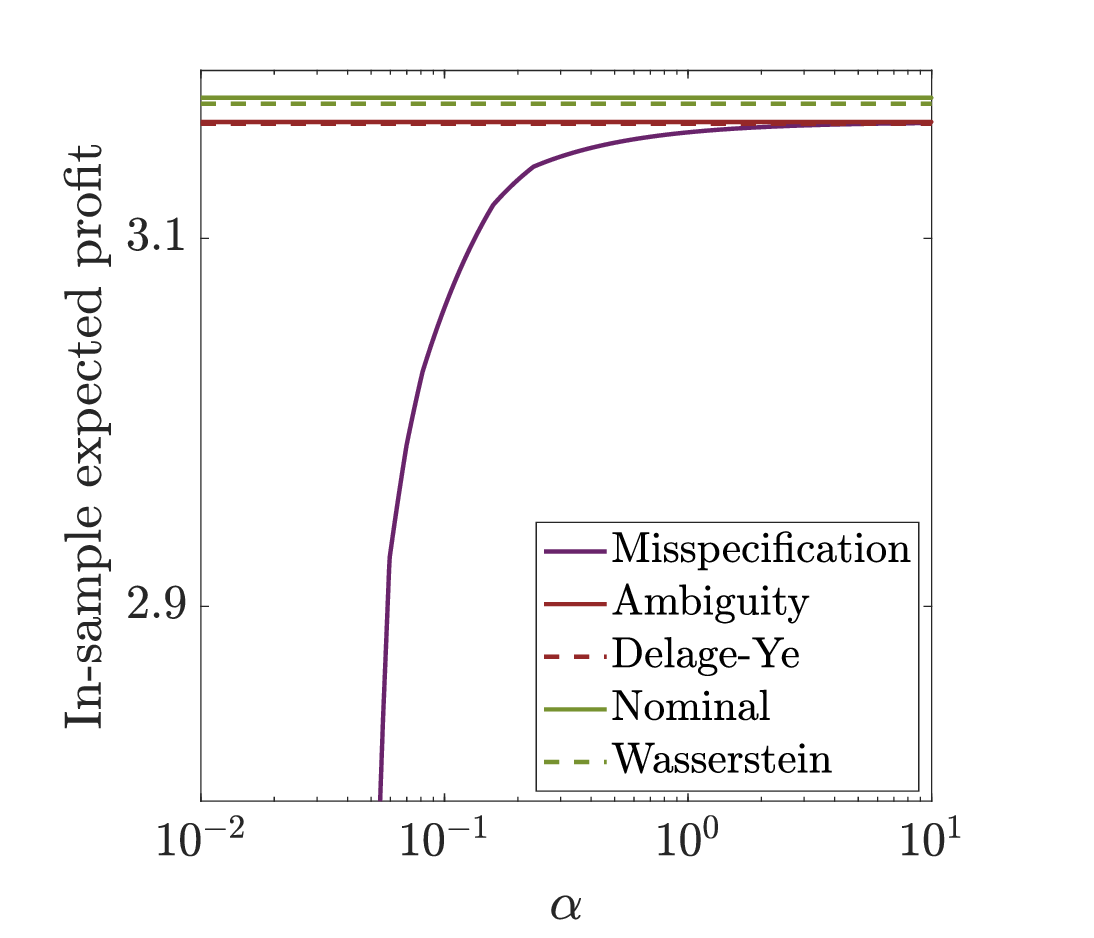}
\end{subfigure}
\begin{subfigure}{.33\textwidth}
\centering
\includegraphics[width=1.11\linewidth]{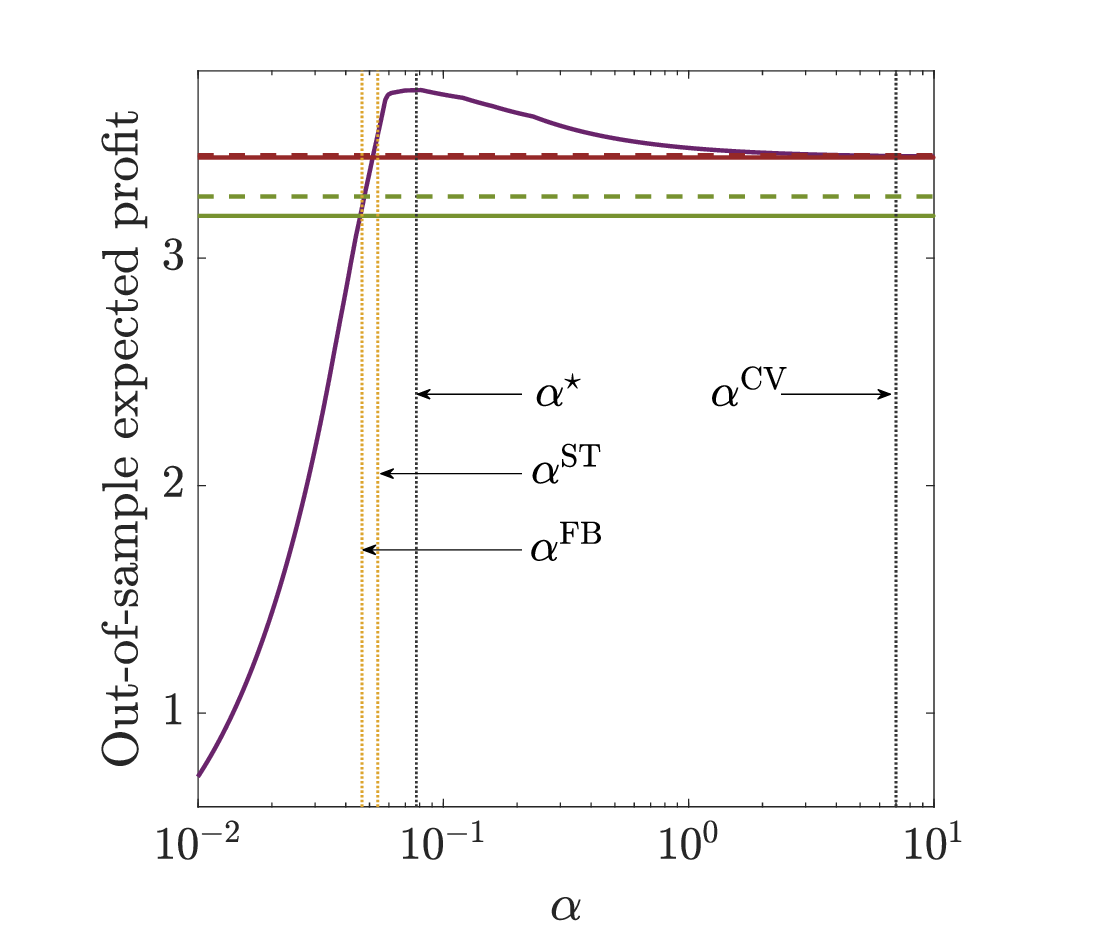}
\end{subfigure}
\begin{subfigure}{.33\textwidth}
\centering
\includegraphics[width=1.11\linewidth]{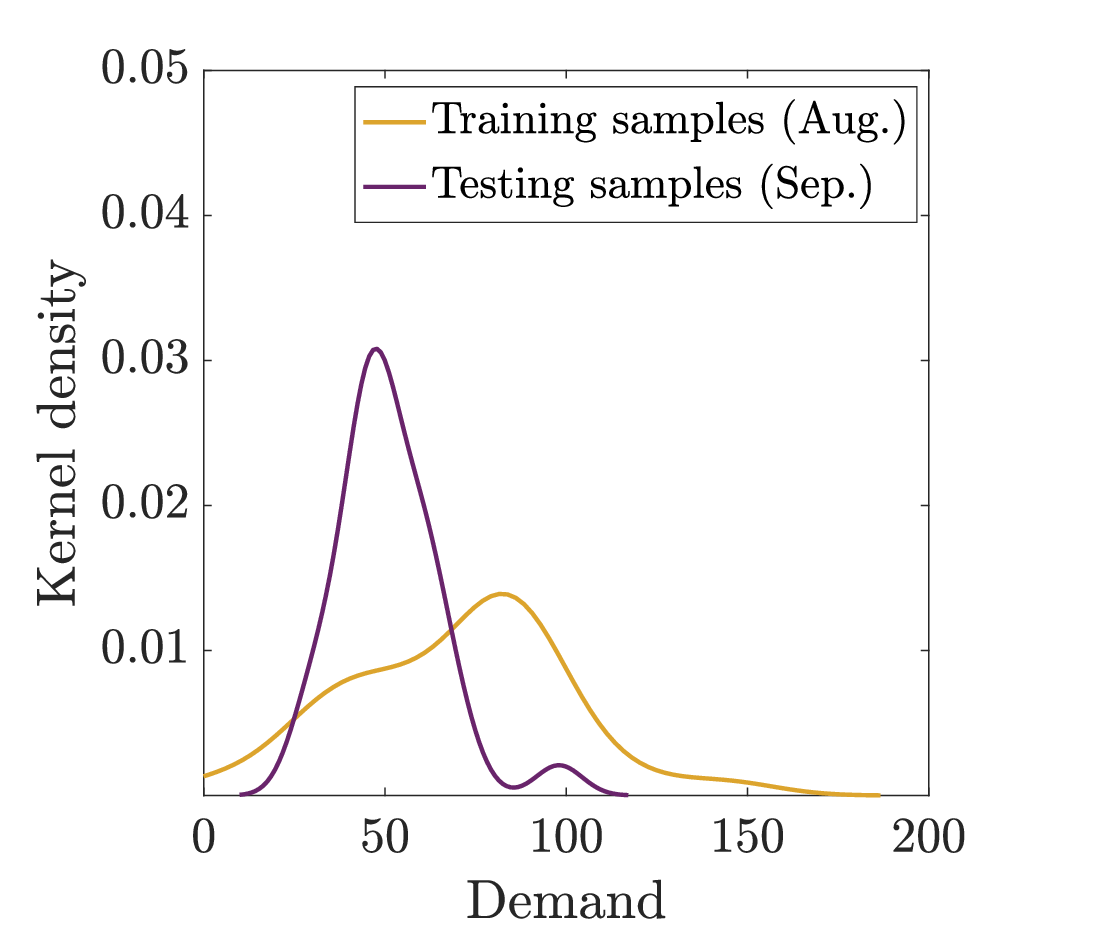}
\end{subfigure}
\begin{subfigure}{.33\textwidth}
\centering
\includegraphics[width=1.11\linewidth]{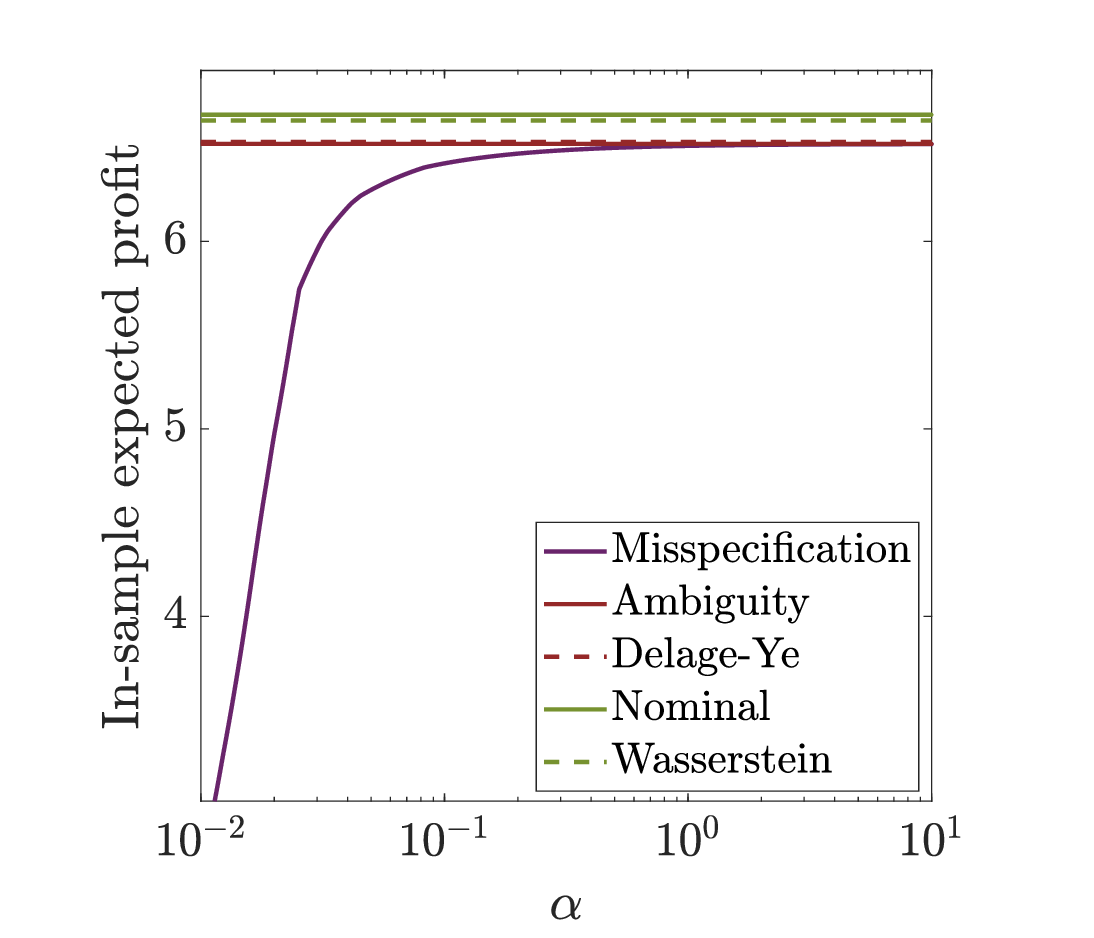}
\end{subfigure}
\begin{subfigure}{.33\textwidth}
\centering
\includegraphics[width=1.11\linewidth]{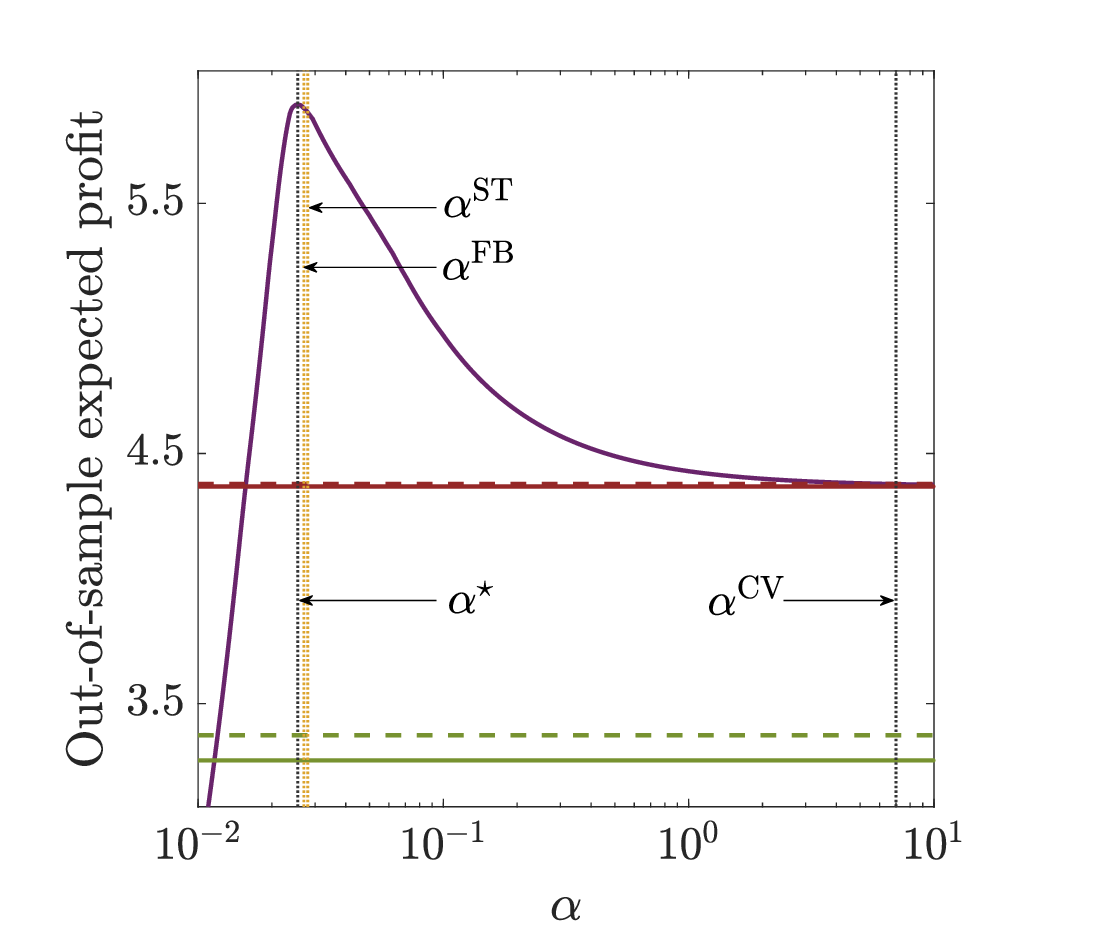}
\end{subfigure}
\caption{\textnormal{Demand density and performance of \ref{prob:misspecification}, \ref{prob:ambiguity}, \ref{prob:Delage}, \ref{prob:nominal}, and \ref{prob:wasserstein} when variability is low (first row), moderate (second row), and high (third row).}}
\vspace{-4mm}
\label{fig:Out-of-sample alpha real data}
\end{figure}

Our data consists of a set of SKUs, and we first consider the popular drinking water (as mentioned in \Cref{fig:demand data} and denoted by ${\rm SKU}_0$) that shows non-stationary characteristics, especially in the sense of monthly mean and variance. For the experimental purpose, we consider the demand data in two consecutive months as training and testing samples. In particular, according to the variability in sales series and the associated mean-variance change over the consecutive months as demonstrated in~\Cref{fig:demand data}, we identify October vs. November as a \textit{low-variability} scenario, January vs. February as a \textit{moderate-variability} scenario, and August vs. September as a \textit{high-variability} scenario. See the first column of~\Cref{fig:Out-of-sample alpha real data}, highlighting the non-stationarity of the demand.

At the end of the training month, the newsvendor obtains the demand observations (training samples) $\hat{v}_1,\dots,\hat{v}_N$ and needs to decide an order quantity to satisfy the random demand that will materialize in the testing month. The newsvendor solves \ref{prob:nominal} with the empirical distribution based on $\hat{v}_1,\dots,\hat{v}_N$. For \ref{prob:ambiguity}, the newsvendor estimates $\mu$ and $\sigma^2$ of the mean-variance ambiguity set $\mathcal{A}$ via sample mean and sample variance, that is, $\mu = \frac{1}{N}\sum_{i = 1}^N \hat{v}_i$ and $\sigma^2 = \frac{1}{N}\sum_{i = 1}^N \hat{v}_i^2 - \big(\frac{1}{N}\sum_{i =1}^N \hat{v}_i\big)^2$. However, taking these estimates as the mean and variance of the demand in the testing month may lead to misspecification (recall from \Cref{fig:demand data}). It is thus meaningful to consider \ref{prob:misspecification} with different values of $\alpha$. For a comprehensive comparison, we also consider two other benchmarks for \ref{prob:ambiguity}. The first is the Wasserstein newsvendor model, {\it i.e.},
\begin{equation}\label{prob:wasserstein}\tag{\sc Wasserstein}
\max_{q \geq 0} \; \min_{G \in \mathcal{B}(\theta)} \; \mathbb{E}_{G}[\pi(q,\tilde{v})],
\end{equation}
where $\mathcal{B}(\theta)$ is defined in Section~\ref{sec:ambiguity set} and its solution can be derived via \Cref{thm:double misspecification} with $\alpha \to \infty$. The other is the newsvendor model based on the ambiguity set proposed in \cite{delage2010distributionally}:
\begin{equation}\label{prob:Delage}\tag{\sc Delage-Ye}
\max_{q \geq 0} \; \min_{G \in \mathcal{D}} \; \mathbb{E}_{G}[\pi(q,\tilde{v})],
\end{equation}
where $\mathcal{D} = \{G \in \mathcal{P} \mid \mu - \sqrt{\gamma_1}\sigma \leq \mathbb{E}_G[\tilde{v}] \leq \mu + \sqrt{\gamma_1}\sigma,\;\mathbb{E}_G[\tilde{v}^2] \leq -\mu^2 + \gamma_2\sigma^2 + 2\mu\mathbb{E}_G[\tilde{v}]\}$, and $\gamma_1 \geq 0$ and $\gamma_2 \geq 1$ are two constants. \ref{prob:Delage} does not admit any closed-form solution and is solved using the reformulation in \cite{delage2010distributionally}. Hyperparameters $\theta$ and $(\gamma_1,\gamma_2)$ in these two benchmarks are selected via the cross-validation approach. In fact, we also test the performance of the newsvendor model with uncertain mean and variance as depicted in \Cref{remark:natarajan}, which turns out to be very similar to that of \ref{prob:Delage}. 

\subsection{Results and Discussion}\label{sec:numerical single}

\Cref{fig:Out-of-sample alpha real data} summarizes the results of different approaches in the three scenarios where the variability between training and testing samples is, respectively, low, moderate, and high. Regarding the in-sample expected profit, in all scenarios, \ref{prob:nominal} achieves the highest value, while \ref{prob:misspecification} converges to \ref{prob:ambiguity} as $\alpha$ approaches infinity. Focusing on the out-of-sample performance, \textit{first}, in the low-variability scenario, \ref{prob:nominal} outperforms \ref{prob:ambiguity}, and under most values of $\alpha$, \ref{prob:misspecification} does not yield a higher out-of-sample expected profit than \ref{prob:nominal}---an intuitive result since the demand process is quite stationary. For \ref{prob:Delage} and \ref{prob:wasserstein}, they outperform \ref{prob:ambiguity} and \ref{prob:nominal}, and almost achieve the best performance of \ref{prob:misspecification}. \textit{Furthermore}, in the moderate-variability scenario,  \ref{prob:nominal} yields \$$3.19$ in the out-of-sample expected profit while \ref{prob:misspecification} (resp., \ref{prob:ambiguity}) results in an improvement at most \$$0.54$ or 16.93\% in percentage (resp., an improvement of \$$0.25$ or 7.84\%). We also note that \ref{prob:Delage} achieves the same expected profit as \ref{prob:ambiguity} while \ref{prob:wasserstein} yields an improvement of \$$0.09$ or 2.68\% compared to \ref{prob:nominal}. \textit{Finally}, in the high-variability scenario, \ref{prob:ambiguity} and \ref{prob:misspecification} have an even larger improvement in out-of-sample expected profit: \ref{prob:nominal} yields \$$3.24$  while \ref{prob:misspecification} (resp., \ref{prob:ambiguity}) results in an improvement at most \$$2.65$ or 81.79\% (resp., an improvement of \$$1.18$ or 36.42\%). We observe that \ref{prob:Delage} also achieves the same expected profit as \ref{prob:ambiguity} while \ref{prob:wasserstein} yields an improvement of \$$0.10$ or 3.06\% compared to \ref{prob:nominal}.

In the third column of Figure~\ref{fig:Out-of-sample alpha real data}, we first emphasize two values of the index of misspecification aversion, $\alpha^{\rm CV}$ and $\alpha^{\star}$. The former is selected via cross-validation using the training data, while the latter is the one that achieves the largest out-of-sample expected profit. We would like to highlight two important observations and insights as follows.

\begin{enumerate}
\item In the low-variability scenario with stationary demand, \ref{prob:misspecification} with $\alpha^{\rm CV}$ performs quite close to \ref{prob:nominal} and \ref{prob:ambiguity}, and $\alpha^{\rm CV}$ is also close to $\alpha^{\star}$. This, not only justifies the predominance of mean-variance statistics in capturing the underlying distribution for the newsvendor's decision, but also implies the usefulness of the cross-validation method in calibrating the parameter $\alpha$ of \ref{prob:misspecification}, in the situation of stationary demand. 

\item In the moderate-variability and high-variability scenarios with non-stationary demand, the calibrated \ref{prob:misspecification} model with $\alpha=\alpha^{\rm CV}$ yields an out-of-sample performance close to that of \ref{prob:ambiguity}, which, however, is  \textit{far away} from the best performance that \ref{prob:misspecification} could achieve with $\alpha=\alpha^{\star}$. This confirms the implication of the finite-sample performance guarantee derived in \Cref{thm:guarantee}: in the situation of distribution shift where the testing samples vary highly from the training samples, cross-validation---purely relying on training samples---could be depreciative in its effectiveness for calibrating a model's parameter.
\end{enumerate}

It is generally impossible to calibrate the index $\alpha$ of misspecification aversion under arbitrary and unknown distribution shifts (see the opening discussion by \citealt{sutter2021robust}, p.1). However, in certain practical situations (for example, companies may conduct test markets before launching a new product to gather limited testing data), partial information about the distribution shift may be available through a limited number of testing samples drawn from the out-of-sample distribution (\citealp{sugiyama2012machine}). To address such cases, we propose two calibration strategies: a formula-based approach and a stress-testing approach. We provide a brief description in the following and leave the implementation details to \Cref{appendix:calibration}.

The formula-based approach directly explores formula~\eqref{equ:alpha-N} derived in \Cref{thm:guarantee} to determine $\alpha$. Specifically, we use testing samples to estimate $d_{\rm OT}(F,D)$, while treating $\epsilon_N$ as a tuning parameter and using cross-validation to select $\epsilon_N$. The stress-testing approach can be viewed as a controlled stress test that proceeds with two steps: constructing a stress-testing distribution and validating based on the constructed stress-testing distribution. In the construction step, to mimic the testing environment, we build a synthetic stress-testing distribution that matches the estimated distribution shift based on the testing samples. In the validation step, we cross-validate $\alpha$ based on the constructed testing distribution.

The third column of \Cref{fig:Out-of-sample alpha real data} shows that in all three testing scenarios, both $\alpha^{\rm FB}$ (by the formula-based approach) and $\alpha^{\rm ST}$ (by the stress-testing approach) select relatively small values that, compared to $\alpha^{\rm CV}$, are more closer to the optimal $\alpha^\star$. 

\begin{enumerate}
\item In the low-variability scenario---where the distribution shift is minimal---the formula-based $\alpha^{\rm FB}$ closely aligns with $\alpha^{\rm CV}$ (which is already good enough), while the stress-testing $\alpha^{\rm ST}$ tends to yield smaller values. This behavior is expected: although the actual shift is small, the empirical distance between the training distribution and the few testing samples often overestimates the shift. As a result, both $\alpha^{\rm ST}$ and $\alpha^{\rm FB}$ are calibrated to be relatively conservative.

\item In the moderate-variability scenario---where the distribution shift is more pronounced---the empirical distance provides a more accurate reflection of the shift compared to the low-variability scenario. We observe that the stress-testing approach yields an $\alpha^{\rm ST}$ with slightly better out-of-sample performance than $\alpha^{\rm CV}$, while the formula-based $\alpha^{\rm FB}$ performs marginally worse. Notably, both $\alpha^{\rm FB}$ and $\alpha^{\rm ST}$ are closer to the optimal $\alpha^\star$ than in the low-variability scenario. When the distribution shift is substantial, the empirical distance reliably captures the extent of the shift. Consequently, both $\alpha^{\rm FB}$ and $\alpha^{\rm ST}$ yield better-calibrated estimates that closely approximate the optimal $\alpha$, and they outperform the plain vanilla $\alpha^{\rm CV}$, which is validated solely over training samples.
\end{enumerate}

These findings demonstrate that with limited knowledge of the distribution shift, our proposed calibration strategies perform effectively in mitigating its impact---particularly in scenarios with substantial variability. We would like to acknowledge that our proposal is merely an initial attempt, and that developing a rigorous statistical framework for inferring distributional shifts remains an active and emerging area of research. As such, this aspect lies beyond the scope of the present study and is deferred to future work, as discussed in the conclusion.

\begin{table}[t!]
\centering \small
\begin{tabular}{c|c|c|cccccc}
\hline\hline
\multicolumn{1}{c|}{\multirow{2}{*}{$\alpha$}}  & \multicolumn{1}{c|}{\multirow{2}{*}{Performance Comparison}} &  \multicolumn{1}{c|}{\multirow{2}{*}{Percentage}} &  \multicolumn{3}{c}{\multirow{1}{*}{$\Pi_{\rm M}$}} &  \multicolumn{2}{c}{\multirow{1}{*}{$\Pi_{\rm A}$}} & \multicolumn{1}{c}{\multirow{1}{*}{$\Pi_{\rm S}$}} \\ \cline{5-5} \cline{7-7} \cline{9-9}
& &   &&  \multicolumn{1}{c}{\multirow{1}{*}{\hspace{-1.05mm}\!\! Mean $\mid$ STD}} & & \multicolumn{1}{c}{\multirow{1}{*}{\hspace{-1.1mm}\!\! Mean $\mid$ STD}} & & \multicolumn{1}{c}{\multirow{1}{*}{\hspace{-1.1mm}\!\! Mean $\mid$ STD}} \\ \cline{1-9} 
\multirow{2}{*}{$\alpha_{\rm low}$} &  $\Pi_{\rm M} > \Pi_{\rm S}~\mbox{and}~\Pi_{\rm M} > \Pi_{\rm A}$ & 28\% & & 24.95 $\mid$ 19.93 & & 14.49 $\mid$ 13.17 & & 15.90 $\mid$ 14.33 \\
& $\Pi_{\rm M} \leq \Pi_{\rm S}~~\mbox{or}~~\Pi_{\rm M} \leq \Pi_{\rm A}$ & 72\% & & 5.09 $\mid$ 7.22 & & 8.60 $\mid$ 8.94 & & 8.80 $\mid$ 8.45 \\ \cline{1-9}
\multirow{2}{*}{$\alpha_{\rm mid}$}
& $\Pi_{\rm M} > \Pi_{\rm S}~\mbox{and}~\Pi_{\rm M} > \Pi_{\rm A}$ & 81\% & & 12.07 $\mid$ 10.95 & & 9.74 $\mid$ 9.72 & & 9.56 $\mid$ 9.15 \\
& $\Pi_{\rm M} \leq \Pi_{\rm S}~~\mbox{or}~~\Pi_{\rm M} \leq \Pi_{\rm A}$ & 19\% & & 13.28 $\mid$ 14.55 & & 12.52 $\mid$ 13.59 & & 15.91 $\mid$ 15.51 \\ \cline{1-9}
\multirow{2}{*}{$\alpha_{\rm high}$} & $\Pi_{\rm M} > \Pi_{\rm S}~\mbox{and}~\Pi_{\rm M} > \Pi_{\rm A}$ & 69\% & & 11.49 $\mid$ 11.04 & & 10.03 $\mid$ 10.10 & & 9.58 $\mid$ 9.50 \\
& $\Pi_{\rm M} \leq \Pi_{\rm S}~~\mbox{or}~~\Pi_{\rm M} \leq \Pi_{\rm A}$ & 31\% & & 11.47 $\mid$ 11.55 & & 10.69 $\mid$ 11.39 & & 13.06 $\mid$ 12.84 \\
\hline\hline
\end{tabular}
\vspace{0.2cm}
\caption{\textnormal{Out-of-sample expected profits over $100$ SKUs. For each SKU, $\alpha_{\rm low} = p/100$, $\alpha_{\rm mid} = p/20$, and $\alpha_{\rm high} = p/10$, where $p$ is the unit price of the SKU. Here, we denote by $\Pi_{\rm M}$, $\Pi_{\rm A}$, and $\Pi_{\rm S}$ the out-of-sample expected profits of \ref{prob:misspecification}, \ref{prob:ambiguity}, and \ref{prob:nominal}, respectively.}}
\label{tab:comparison}
\vspace{-5mm}
\end{table}

We next repeat the above experiment over a pool of $100$ SKUs, for each of which we randomly select two consecutive months as training and testing samples. For each SKU, we consider the same setting as in~\Cref{sec:data and setting} to evaluate the out-of-sample expected profits of \ref{prob:misspecification}, \ref{prob:ambiguity}, and \ref{prob:nominal}. \Cref{tab:comparison} summarizes the number of SKUs that one model outperforms another, and the corresponding mean and standard deviation of the out-of-sample profits of a model over these SKUs. Under a small value of $\alpha$ (\textit{i.e.}, $\alpha_{\rm low}$), for $28\%$ of $100$ SKUs, \ref{prob:misspecification} outperforms both \ref{prob:ambiguity} and \ref{prob:nominal} with a large profit improvement but a large standard deviation; for the majority $72\%$, \ref{prob:misspecification} underperforms either \ref{prob:ambiguity} or \ref{prob:nominal} with a large profit loss and a small standard deviation. Under a medium value of $\alpha$ (\textit{i.e.}, $\alpha_{\rm mid}$),  \ref{prob:misspecification} yields superior performance than both \ref{prob:ambiguity} and \ref{prob:nominal} by noting that \ref{prob:misspecification} outperforms both \ref{prob:ambiguity} and \ref{prob:nominal} for a majority $81\%$ of $100$ SKUs. Even under a high value of $\alpha$ (\textit{i.e.}, $\alpha_{\rm high}$),  \ref{prob:misspecification} also has a fairly good out-of-sample performance, where \ref{prob:misspecification} outperforms both \ref{prob:ambiguity} and \ref{prob:nominal} for a majority $69\%$ of $100$ SKUs. In other words, for each $\alpha\in\{\alpha_{\rm low},\alpha_{\rm mid}, \alpha_{\rm high}\}$, there always quite a proportion of SKUs such that over these products \ref{prob:misspecification} has a better out-of-sample performance than both \ref{prob:ambiguity} and \ref{prob:nominal}, justifying the need of incorporating misspecification to the newsvendor problem.
\begin{figure}[tb]
\begin{subfigure}{.33\textwidth}
\centering
\includegraphics[width=1.11\linewidth]{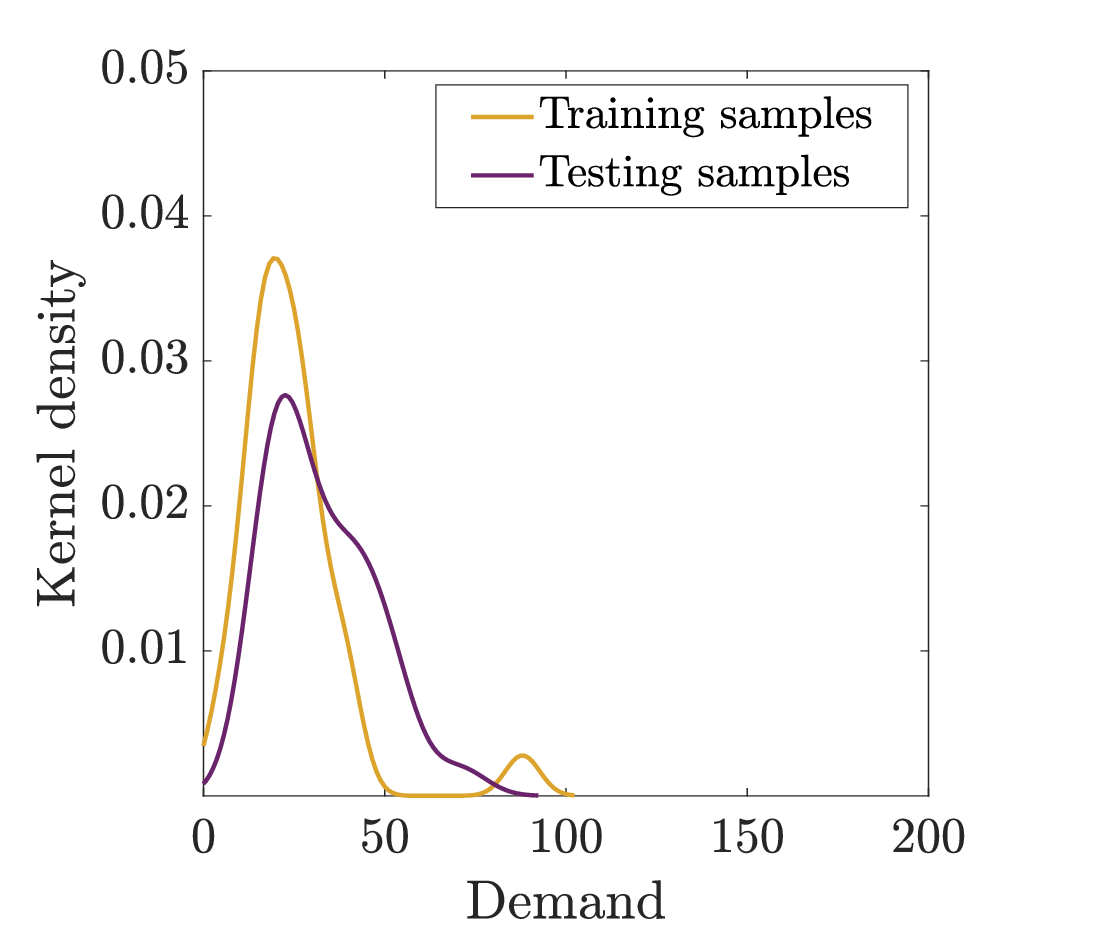}
\end{subfigure}
\begin{subfigure}{.33\textwidth}
\centering
\includegraphics[width=1.11\linewidth]{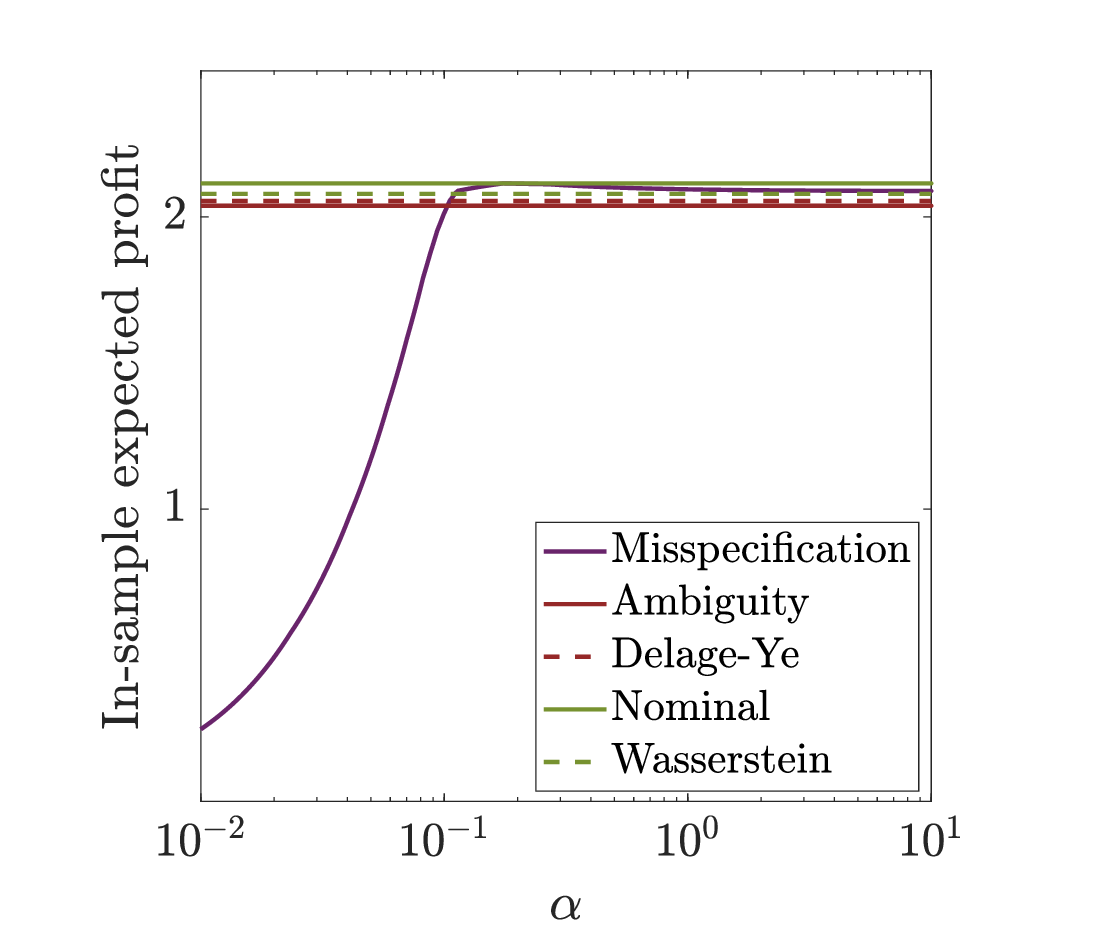}
\end{subfigure}
\begin{subfigure}{.33\textwidth}
\centering
\includegraphics[width=1.11\linewidth]{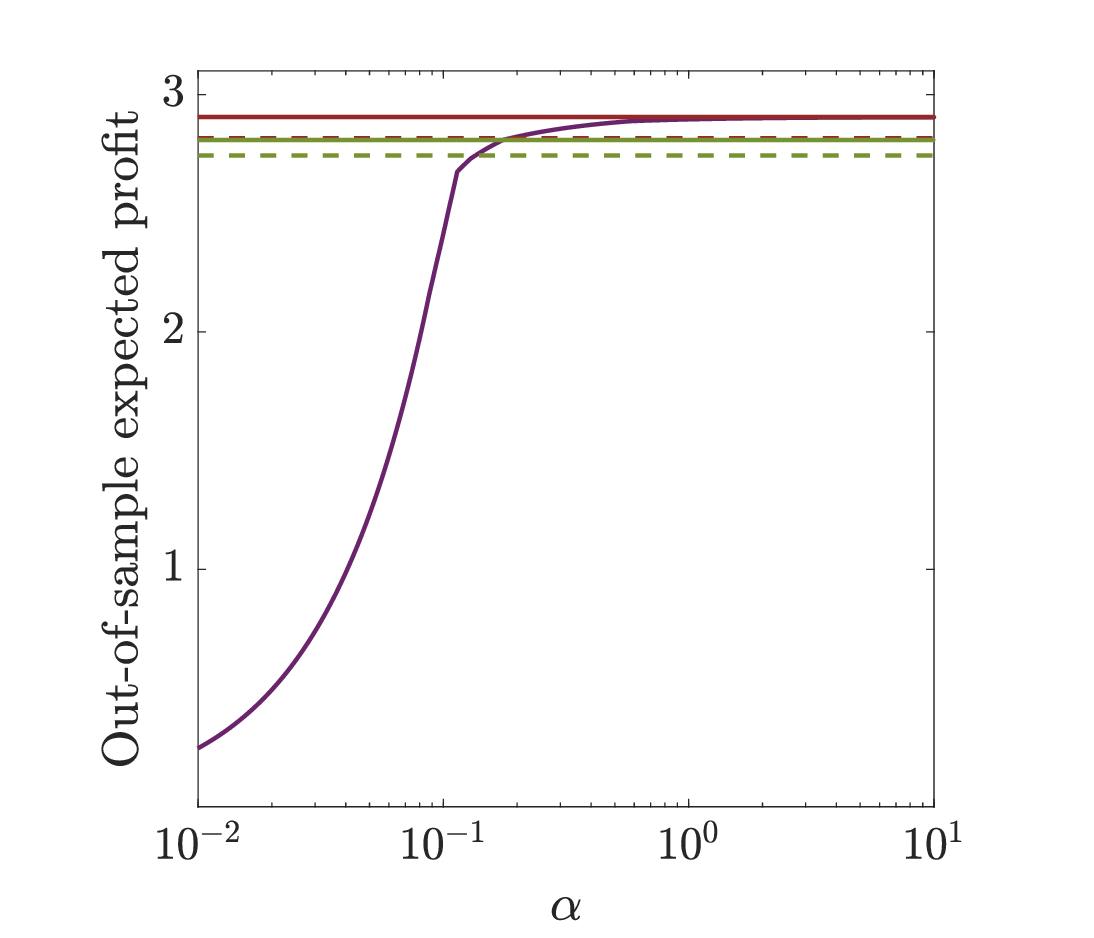}
\end{subfigure}
\caption{\textnormal{Demand density and performance of \ref{prob:misspecification}, \ref{prob:ambiguity}, \ref{prob:Delage}, \ref{prob:nominal}, and \ref{prob:wasserstein} when variability is low (first row), moderate (second row), and high (third row).}}
\vspace{-4mm}
\label{fig:inferior}
\end{figure}

We caution that misspecification aversion does not necessarily lead to performance improvement, as evidenced by the first row of \Cref{fig:Out-of-sample alpha real data} that \ref{prob:misspecification} underperforms \ref{prob:ambiguity} across a wide range of $\alpha$ and by \Cref{tab:comparison} that for certain SKUs, misspecification may result in a performance decline. To elucidate this phenomenon, we present in \Cref{fig:inferior} a representative scenario where the testing demand is modestly higher than that in the training month. In this context, \ref{prob:misspecification} exhibits inferior out-of-sample performance compared to \ref{prob:ambiguity} over all values of $\alpha$. While misspecification aversion can enhance performance when the out-of-sample environment is substantially worse than the training environment, it may adopt a conservative stance that leads to profit loss when the out-of-sample environment is relatively stable or even more favorable than the in-sample environment.

\section{Conclusion}

Since the seminal work of \cite{scarf1958min}, the mean-variance ambiguity set has been widely used for decision-making in mitigating distributional uncertainty. However, in many practical situations, the mean and variance can be misspecified, consequently resulting in inexpert newsvendor decisions. To address this issue, we introduce misspecification upon ambiguity and propose a misspecification-averse (and ambiguity-averse) newsvendor model. We investigate the impact and rationale of misspecification aversion from decision-criterion, operational, and statistical perspectives. We also extend our model to establish a comprehensive framework (multi-products, ambiguity captured by optimal transport, and misspecification measured by total variation distance) for the newsvendor under ambiguity and misspecification.

Our present study focuses on and investigates many aspects of the misspecification-averse newsvendor problem. The framework has several interesting directions remaining unexplored and can be extended to other operational problems, opening up promising avenues for future studies.

\vspace{2mm}
\noindent
\textbf{Statistical inference of distribution shift.} $\;$ As articulated in \Cref{thm:guarantee}, the performance guarantee and the index of misspecification aversion are statistically described with a term of distribution shift. How to estimate the distribution shift is statistically critical and is also practically relevant for calibrating the model of misspecification aversion. 

\vspace{2mm}
\noindent
\textbf{Misspecification in prescriptive analytics}. $\;$ Prescriptive analytics, as an emerging paradigm for data-driven decision-making, seeks a decision rule that maps the observed data to an action, which usually leverages some parametric (structural) assumptions on the uncertainty (see, {\it e.g.}, \citealp{chu2025solving}). These assumed parametric models could misspecify the ground truth. Therefore, we believe that our approach, in marriage with the prescriptive analytics framework, has the potential to mitigate the downside consequences of model misspecification.

\vspace{2mm}
\noindent
\textbf{Structuring distributional uncertainty in other operations management problems}. $\;$ 
Although this study focuses on the newsvendor problem, our analysis can also be applied in other operations management settings, for instance, inventory control, logistics, dynamic pricing, and project management, wherein the distributional uncertainty has been widely acknowledged. Following the spirit of ``all models are wrong, but some are useful'' (\citealp{box1976science}), we believe that treating ambiguity and misspecification differently in a properly structured fashion can differentiate the 
\textit{useful} wrong models and \textit{harmful} wrong models in coping with uncertainty, and therefore enhance the decisions for operations management. 

\ACKNOWLEDGMENT{
The authors thank the Department Editor, the Associate Editor, and two anonymous reviewers for their extraordinary efforts in reviewing the paper and for providing insightful comments that helped improve the paper. Zhi Chen is supported in part by the National Natural Science Foundation of China [72422002, 72394395], the Hong Kong Research Grants Council General Research Fund [CUHK-11502422], and the Asian Institute of Supply Chains and Logistics. Ruodu Wang is supported by the Natural Sciences and Engineering Research Council of Canada [CRC-2022-00141, RGPIN-2024-03728]. Shuming Wang is supported by the National Natural Science Foundation of China [Grants 72471224, 72171221, 71922020, and 71988101], the Fundamental Research Funds for the Central Universities [Grant UCAS-E2ET0808X2], and a grant from the MOE Social Science Laboratory of Digital Economic Forecasts and Policy Simulation at UCAS.
}


{\newcommand{\custombibsize}{\fontsize{11pt}{13pt}\selectfont} 
\bibliographystyle{ormsv080}
\renewcommand*{\bibfont}{\custombibsize}
\bibliography{reference}}

\vspace{2mm}

\begin{thebibliography}{}

\bibitem[Besbes et al.\hspace{-1mm} (2023)]{besbes2023contextual}
Besbes, Omar, Will Ma, Omar Mouchtaki. 2023. From contextual data to newsvendor decisions: On the actual performance of data-driven algorithms. {\it Available at} \url{https://arxiv.org/abs/2302.08424}.

	
\bibitem[{Dowson and Landau \hspace{-1mm}(1982)}]{dowson1982frechet}
Dowson, D.C., B.V. Landau. 1982. The Fr{\'e}chet distance between multivariate normal distributions. {\it Journal of Multivariate Analysis} \textbf{12}(3) 450--455.
	
\bibitem[Gao and Kleywegt \hspace{-1mm}(2023)]{gao2023distributionally} 
Gao, Rui, Anton Kleywegt. 2023. Distributionally robust stochastic optimization with Wasserstein distance. {\it Mathematics of Operations Research} \textbf{48}(2) 603–655.

\bibitem[{Gelbrich\hspace{-1mm}(1990)}]{gelbrich1990formula2}
Gelbrich, Matthias. 1990. On a formula for the $L^2$ Wasserstein metric between measures on Euclidean and Hilbert spaces. {\it Mathematische Nachrichten} \textbf{147}(1) 185--203.

\bibitem[Hanasusanto et al. \hspace{-3mm} (2015)]{hanasusanto2015distributionally1}
Hanasusanto, Grani A., Daniel Kuhn, Stein W. Wallace, Steve Zymler. 2015. Distributionally robust multi-item newsvendor problems with multimodal demand distributions. {\it Mathematical Programming} \textbf{152}(1) 1--32.

\bibitem[Kupiec \hspace{-2mm} (2002)]{kupiec2002stress}
Kupiec, Paul. 2002. Stress testing in a value at risk framework. {\it Risk management: value at risk and beyond}, 76--99. Cambridge University Press.

\bibitem[M{\"u}ller\hspace{-1mm} (1997)]{muller1997integral}
M{\"u}ller, Alfred. 1997. Integral probability metrics and their generating classes of functions. {\it Advances in Applied Probability} \textbf{29}(2) 429--443.

\bibitem[Natarajan and Teo \hspace{-2mm} (2017)]{natarajan2017reduced}
Natarajan, Karthik, Chung-Piaw Teo. 2017. On reduced semidefinite programs for second order moment bounds with applications. {\it Mathematical Programming} \textbf{161} 487--518.

\bibitem[Natarajan et al.\hspace{-1mm} (2018)]{natarajan2018asymmetry1}
Natarajan, Karthik, Melvyn Sim, Joline Uichanco. 2018. Asymmetry and ambiguity in newsvendor models. {\it Management Science} \textbf{64}(7) 3146--3167.

\bibitem[{Rigollet and Jan-Christian \hspace{-1mm}(2023)}]{rigollet2023high}
Rigollet, Philippe, H{\"u}tter Jan-Christian. 2023. High-dimensional statistics. {\it Available at} \url{https://arxiv.org/abs/2310.19244}.

\bibitem[{Scarf \hspace{-2mm} (1958)}]{scarf1958min1}
Scarf, Herbert. 1958. A min-max solution of an inventory problem. {\it Studies in the mathematical theory of inventory and production.}

\bibitem[{Villani \hspace{-1mm}(2009)}]{villani2009optimal1}
Villani, C{\'e}dric. 2009. {\it Optimal transport: Old and new}. Springer.

\bibitem[Wainwright \hspace{-1mm}(2019)]{Wainwright2019high2}
Wainwright, Martin. 2019. High-dimensional statistics: A non-asymptotic viewpoint. {\it Cambridge University Press}.

\bibitem[Zhang et al.\hspace{-1mm} (2025)]{zhang2025short}
Zhang, Luhao, Jincheng Yang, Rui Gao. 2025. A short and general duality proof for Wasserstein distributionally robust optimization. {\it Operations Research} \textbf{73}(4) 2146--2155.

\end{thebibliography}

\newpage
\ECSwitch 
\vspace{-6mm}
\ECHead{\centering { E-Companion to \\``Newsvendor under Ambiguity and Misspecification"}}
{\SingleSpacedXI \small

\AUTHOR{Feng Liu, Zhi Chen, Ruodu Wang, and Shuming Wang}

\section{Technical Lemmas}\label{sec:techniques}

\begin{lemma}\label{lemma:convexity}
For any $F, G\in\mathcal{P}$, $d_{\rm OT}(F,G)$ in~\eqref{equ:optimal transport 2} is jointly convex in $F$ and $G$.
\end{lemma}

\noindent{\bf Proof.} $\;$ 
Given $F_1, F_2, G_1, G_2\in\mathcal{P}$, let $\Gamma_1\in\mathcal{W}(F_1, G_1)$ and $\Gamma_2\in\mathcal{W}(F_2, G_2)$ be the joint distributions that solve the corresponding optimal transport. That is, $d_{\rm OT}(F_1, G_1)=\int_{\mathbb{R}_+\times\mathbb{R}_+} |u - v|^2 \Gamma_1({\rm d} u, {\rm d}v)$ and $d_{\rm OT}(F_2, G_2)=\int_{\mathbb{R}_+\times\mathbb{R}_+} |u - v|^2\Gamma_2({\rm d} u, {\rm d}v)$.
We first show that the joint distribution $\Gamma_\lambda=(1-\lambda) \Gamma_1+\lambda \Gamma_2$ has marginals $F_\lambda = (1-\lambda) F_1 + \lambda F_2$ and $G_\lambda = (1-\lambda) G_1 + \lambda G_2$. In fact, for any Borel set $\mathfrak{B} \subseteq \mathbb{R}_+$, we have
\[
\begin{array}{ll}
\Gamma_\lambda(\mathfrak{B} \times \mathbb{R}_+) = (1-\lambda) \Gamma_1(\mathfrak{B} \times \mathbb{R}_+) + \lambda \Gamma_2(\mathfrak{B} \times \mathbb{R}_+) = (1-\lambda) F_1(\mathfrak{B})+ \lambda F_2(\mathfrak{B}) \\
\Gamma_\lambda(\mathbb{R}_+ \times \mathfrak{B}) = (1-\lambda) \Gamma_1(\mathbb{R}_+ \times \mathfrak{B}) + \lambda \Gamma_2(\mathbb{R}_+ \times \mathfrak{B}) = (1-\lambda) G_1(\mathfrak{B})+ \lambda G_2(\mathfrak{B}),
\end{array}
\]
which implies that $\Gamma_\lambda\in\mathcal{W}(F_\lambda, G_\lambda)$. As a result,
\[
\begin{array}{cl}
d_{\rm OT}(F_\lambda, G_\lambda)
\leq \int_{\mathbb{R}_+\times\mathbb{R}_+} |u - v|^2 \;\Gamma_\lambda({\rm d} u, {\rm d} v)=(1-\lambda)\cdot d_{\rm OT}(F_1, G_1) + \lambda\cdot d_{\rm OT}(F_2, G_2),
\end{array}
\]
which concludes the proof.  
\hfill $\square$
\vspace{2mm}

\begin{lemma}\label{lemma:misspecification constrained}
Suppose that $0 < \varepsilon < \infty$, and let $\Upsilon^\star_\varepsilon$ denote the optimal value of problem~\eqref{prob:misspecification constrained} corresponding to $\varepsilon$, and $\Pi^\star_\alpha$ the optimal value of \ref{prob:misspecification} corresponding to $\alpha$. Then,
\begin{equation}\label{prob:misspecification dual}
\Upsilon^\star_\varepsilon = \max_{\alpha \geq 0} \; \{\Pi^\star_\alpha - \varepsilon\alpha\}.
\end{equation}
The optimal value $\alpha^\star$ of problem~\eqref{prob:misspecification dual} is achieved, and the optimal solution $(q_{\alpha^\star}^\star, F_{\alpha^\star}^\star)$ of \ref{prob:misspecification} associated with $\alpha=\alpha^\star$ is also the optimal solution to problem~\eqref{prob:misspecification constrained}.
\end{lemma}

\noindent{\bf Proof.} $\;$ 
We start by fixing the order quantity. For any fixed $q \geq 0$, we define a function $\varphi(q, \varepsilon) = \min_{d(F,\mathcal{A}) \leq \varepsilon}\mathbb{E}_F[\pi(q,\tilde{v})]$. First, it is clear that $\varphi(q, \varepsilon)$ is decreasing in $\varepsilon$. Second, $\varphi(q, \varepsilon)$ is bounded on $\mathbb{R}_+$ because for any $\varepsilon \in (0,+\infty)$, $-cq \leq \varphi(q, \varepsilon) \leq (p-c)q$. Third, $\varphi(q,\varepsilon)$ is convex in $\varepsilon$ on $\mathbb{R}_{++}$. To see this, we fix $F_1, F_2\in\mathcal{P}$ and $\varepsilon_1, \varepsilon_2\geq 0$ such that $d(F_1,\mathcal{A}) \leq \varepsilon_1$ and $d(F_2,\mathcal{A})\leq\varepsilon_2$. Note that $d(F,\mathcal{A}) = \min_{G\in\mathcal{A}}d_{\rm OT}(F,G)$ is convex in $F$ because $d_{\rm OT}(F,G)$ is jointly convex in $F$ and $G$ (\Cref{lemma:convexity}) and maximization over a convex set preserves convexity. For any $\lambda\in[0,1]$, it holds that $d(\lambda F_1 + (1-\lambda) F_2, \mathcal{A})\leq \lambda d(F_1,\mathcal{A}) + (1-\lambda) d(F_2,\mathcal{A}) \leq \lambda\varepsilon_1 + (1-\lambda)\varepsilon_2$,
where the first inequality follows from the convexity of $d(\cdot, \mathcal{A})$. Thus, we have $\varphi(q, \lambda\varepsilon_1 + (1-\lambda)\varepsilon_2) \geq \mathbb{E}_{\lambda F_1 + (1-\lambda)F_2}[\pi(q,\tilde{u})] = \lambda\cdot\mathbb{E}_{F_1}[\pi(q,\tilde{u})] + (1 - \lambda)\cdot \mathbb{E}_{F_2}[\pi(q,\tilde{u})]$.
Taking the minimum over $F_1$ and $F_2$ yields $\varphi(q,\lambda\varepsilon_1 + (1-\lambda)\varepsilon_2) \geq \lambda\varphi(q,\varepsilon_1) + (1-\lambda)\varphi(q,\varepsilon_2)$,
which establishes the convexity (and continuity) of $\varphi(q,\varepsilon)$.
Finally, given $\alpha>0$, the Legendre transform of the convex function $\varphi(q, \cdot)$ is
\[
\begin{array}{ll}
\displaystyle \varphi^*(q, \alpha) = \min_{\varepsilon\geq 0} \; \{\alpha\varepsilon + \varphi(q, \varepsilon)\} & = \displaystyle \min_{\varepsilon \geq 0}\;\min_{F\in\mathcal{P}}\;\{ \alpha\varepsilon + \mathbb{E}_F[\pi(q,\tilde{u})] : d(F,\mathcal{A}) \leq \varepsilon\} = \displaystyle \min_{F\in\mathcal{P}}\;\{\alpha\cdot d(F, \mathcal{A}) + \mathbb{E}_F[\pi(q,\tilde{u})]\},
\end{array}
\]
which is concave in $\alpha$. Note that the above relation also holds for $\alpha = 0$. For any $\varepsilon > 0$, applying the Legendre transform on the concave function $\varphi^*(q, \cdot)$ yields
\[
(\varphi^*(q, \varepsilon))^* = \max_{\alpha \geq 0} \; \{\varphi^*(q, \alpha) - \alpha \varepsilon\} = \max_{\alpha \geq 0} \; \min_{F\in\mathcal{P}} \; \{\mathbb{E}_F[\pi(q,\tilde{u})] + \alpha\cdot d(F, \mathcal{A}) - \varepsilon \alpha \}.
\]
Since $\varphi(q, \varepsilon)$ is bounded, convex, and continuous, for any $\varepsilon> 0$, it holds that $\varphi(q, \varepsilon) = (\varphi^*(q, \varepsilon))^*$. We now optimize the order quantity. Maximizing over $q\geq 0$ yields
\begin{align}
\displaystyle \Upsilon^\star_\varepsilon = \max_{q \geq 0} \; \varphi(q, \varepsilon)
& \displaystyle = \max_{q \geq 0} \; \min_{d(F,\mathcal{A}) \leq \varepsilon} \; \mathbb{E}_F[\pi(q,\tilde{u})] \label{equ:constrained} \\ \label{equ:dual 1}
& \displaystyle = \max_{\alpha \geq 0} \; \max_{q \geq 0} \; \big\{\min_{F\in\mathcal{P}} \; \{\mathbb{E}_F[\pi(q,\tilde{u})] + \alpha\cdot d(F, \mathcal{A})\} - \varepsilon \alpha \big\} \\ \label{equ:dual 3}
& \displaystyle = \max_{\alpha \geq 0} \; \{\Pi^\star_\alpha - \varepsilon \alpha\}, 
\end{align}
where \eqref{equ:dual 3} follows from the definition of $\Pi^\star_\alpha$. Indeed, the optimal value of \eqref{equ:dual 3} is achieved at some finite $\alpha^\star \geq 0$ because the function $\Pi^\star_\alpha - \varepsilon \alpha$  is continuous in $\alpha$ with $\Pi^\star_0 - \varepsilon\cdot 0 = 0$ and $\Pi^\star_\alpha - \varepsilon\alpha\to -\infty$ as $\alpha \to \infty$. Finally, note that the optimal solution $(\alpha^\star, q^\star, F^\star)$ of \eqref{equ:dual 1} has its part $(q^\star, F^\star)$ being the optimal solution of \ref{prob:misspecification} with $\alpha = \alpha^\star$, which we have denoted by $(q^\star_\alpha, F^\star_\alpha)$.  Hence, $(q^\star_{\alpha^\star}, F^\star_{\alpha^\star})$ is optimal to \eqref{equ:constrained}, \textit{i.e.}, problem~\eqref{prob:misspecification constrained}. This concludes the proof.
\hfill $\square$
\vspace{2mm}

\begin{lemma}[{\sc Interchangeability}, \textnormal{\citealp{zhang2025short}}]\label{lemma:interchangeability}
Given $G \in \mathcal{P}_M$ and a $G$-measurable function $h(\cdot):\mathbb{R}_+^M \mapsto \mathbb{R}$ with $\mathbb{E}_{G}[h(\bmt{u})] < +\infty$, we have
\begin{equation*}
\min_{F\in\mathcal{P}_M,\;\Gamma\in\mathcal{W}(F,G)}\mathbb{E}_{\Gamma}\big[h(\bmt{u})+\alpha\cdot \|\bmt{u}-\bmt{v}\|_2^2\big] \;=\; \mathbb{E}_{G}\Big[\min_{\bm{u} \geq 0} \; \big\{h(\bm{u}) + \alpha \cdot \|\bm{u} - \bmt{v}\|_2^2 \big\}\Big].
\end{equation*}
\end{lemma}

\noindent{\bf Proof.} $\;$ 
Note that the $(\mathbb{R}_+^M, \|\cdot\|_2)$ is a metric space and every Borel probability measure $G \in \mathcal{P}(\mathbb{R}_+^M)$ is tight. The result is then an immediate consequence of proposition~2 in \cite{zhang2025short}. 
\hfill $\square$
\vspace{2mm}

\begin{lemma}\label{lemma:newsvendor type-2}
Let $\ell(\alpha,q,v) = \min_{u \geq 0} \; \{\pi(q, u) + \alpha(u-v)^2\}$. For any fixed $v\geq 0$,
\begin{enumerate}
\item[(\textit{i})] if $0\leq q\leq \frac{p}{4\alpha}$, then $\ell(\alpha,q,v) = \min\{\alpha v^2, pq\} - cq$;  

\item[(\textit{ii})] if $q> \frac{p}{4\alpha}$, then
\[
\ell(\alpha,q,v)=\left\{
\begin{array}{ll}
\alpha v^2-cq & 0 \leq v\leq \frac{p}{2\alpha}\\ [1mm]
p\cdot\min\left\{v - \frac{p}{4\alpha}, q\right\} - cq & v > \frac{p}{2\alpha}.
\end{array}
\right.
\]
\end{enumerate}
\end{lemma}

\noindent {\bf Proof.} $\;$ 
The function $\pi(q, v) + \alpha(v-u)^2$ can be written as a piecewise quadratic function as follows:
\[
g(u,q,v) = p \min \{q, u\} - cq + \alpha(u-v)^2 = 
\left\{
\begin{array}{ll}
\ubar{g}(u) = pu - cq + \alpha(u-v)^2 & u\leq q\\
\bar{g}(u) = pq - cq + \alpha(u-v)^2 & u > q.
\end{array}
\right.
\]
Before proceeding, we denote the left and right derivatives of $g(u,q,v)$ at $u = u_0$ by $g'_-(u_0)$ and $g'_+(u_0)$, respectively. Based on the value of $v$, there are three cases to consider. 
\begin{enumerate}
\item[(\textit{i})] If $v \geq q + \frac{p}{2\alpha}$, \textit{i.e.}, $\ubar{g}_{-}'(q)=2\alpha(q-v) + p \leq 0$, then $\ubar{g}(u)$ is decreasing over $(0,q]$. Note that $\bar{g}'(u) = 2\alpha(u-v)=0$ admits a unique solution $u=v$, implying that $\bar{g}(u)$ is decreasing in $(q,v)$ and increasing in $[v,+\infty)$. Therefore, it holds that $\min_{u \geq 0} g(u,q,v) = \bar{g}(v) = pq - cq$.

\item[(\textit{ii})] If $v \leq q$, \textit{i.e.}, $\bar{g}_{+}'(q) = 2\alpha(q-v) \geq 0$, and we have $\ubar{g}'(u) = 2\alpha(u-v) + p = 0$ admits a unique solution $u = v - \frac{p}{2\alpha}$, implying that $\ubar{g}(u)$ is first decreasing in $[0,(v-p/2\alpha)^+]$ and then increasing in $((v-p/2\alpha)^+,q]$, and $\bar{g}(u)$ is increasing in $(q,+\infty)$. Therefore, it holds that $\min\limits_{u \geq 0} g(u,q,v)  = \ubar{g}((v-\frac{p}{2\alpha})^+)$.
	
\item[(\textit{iii})] If $q \leq v \leq q + \frac{p}{2\alpha}$, {\it i.e.}, 
$\ubar{g}_{-}'(q) \geq 0$ and $\bar{g}_{+}'(q) \leq 0$,  then both $\ubar{g}(u)$ and $\bar{g}(u)$ admit a corresponding minimizer in the domain. It follows that $\min\limits_{u \geq 0} g(u,q,v) = \min\{\ubar{g}((v-\frac{p}{2\alpha})^+), \bar{g}(v)\}.$
\end{enumerate}

Note that when $v<\frac{p}{2\alpha}$, we have $\ubar{g}\big((v-\frac{p}{2\alpha})^+\big) = 
\alpha v^2 - cq$ and when $\frac{p}{2\alpha}\leq v \leq q+\frac{p}{2\alpha}$, $\ubar{g}\big((v-\frac{p}{2\alpha})^+\big) = 
p( v - \frac{p}{4\alpha}) - cq$.
We next discuss based on the value of $q$. When $q \geq \frac{p}{2\alpha}$, we have
\[
\begin{array}{ll}
\displaystyle \ell(\alpha,q,v) = \min_{u \geq 0} \; g(u,q,v) & = 
\left\{
\begin{array}{ll}
\alpha v^2-cq &  v\leq \frac{p}{2\alpha}\\ 
p(v-\frac{p}{4\alpha})-cq & \frac{p}{2\alpha} \leq v\leq q\\ [1mm]
p\cdot\min\{v-\frac{p}{4\alpha},q\}-cq & q\leq v \leq q+\frac{p}{2\alpha}\\
pq-cq & v\geq q+\frac{p}{2\alpha}
\end{array}
\right.
=\displaystyle
\left\{
\begin{array}{ll}
\alpha v^2-cq & v\leq \frac{p}{2\alpha}\\ 
p(v-\frac{p}{4\alpha})-cq & \frac{p}{2\alpha}\leq v\leq q+\frac{p}{4\alpha}\\ 
pq-cq & v\geq q+\frac{p}{4\alpha};
\end{array}
\right.
\end{array}
\]
When $q \leq \frac{p}{2\alpha}$, we have
\[ 
\ell(\alpha,q,v) = \min_{u \geq 0}\; g(u,q,v) = \displaystyle
\left\{
\begin{array}{ll}
\alpha v^2-cq & v\leq q\\ 
\min\{\alpha v^2, pq\}-cq & q\leq v\leq \frac{p}{2\alpha}\\ 
p\cdot \min\{v-\frac{p}{4\alpha},q\}-cq & \frac{p}{2\alpha}\leq v \leq q+\frac{p}{2\alpha}\\ 
pq-cq & v\geq q+\frac{p}{2\alpha}.
\end{array}
\right.
\]
If $\frac{p}{4\alpha}\leq q\leq \frac{p}{2\alpha}$,  then $\min\{\alpha v^2, pq\} = \alpha v^2$ for $u\in\big[q,\frac{p}{2\alpha}\big]$, resulting in
\[
\ell(\alpha,q,v)=\left\{
\begin{array}{ll}
\alpha v^2-cq & v\leq \frac{p}{2\alpha}\\
p(v-\frac{p}{4\alpha})-cq & \frac{p}{2\alpha}\leq v\leq q+\frac{p}{4\alpha}\\
pq-cq & v\geq q+\frac{p}{4\alpha}.
\end{array}
\right.
\]
If $q \leq \frac{p}{4\alpha}$, then $p\cdot\min\big\{v-\frac{p}{4\alpha},q\big\}-cq=pq-cq$ for $v\in\big[\frac{p}{2\alpha},q+\frac{p}{2\alpha}\big]$. Correspondingly,
\[ 
\ell(\alpha,q,v)=
\left\{
\begin{array}{ll}
\alpha v^2-cq & v \leq \sqrt{\frac{pq}{\alpha}}\\
 pq-cq & v\geq \sqrt{\frac{pq}{\alpha}}
\end{array}
\right.= \min\{\alpha v^2, pq\} - cq.
\]
Consolidating these results based on the three ranges of $q$ then completes the proof.
\hfill $\square$
\vspace{2mm}



\section{Proofs.}\label{appendix:proof}

\noindent{\bf Proof of \Cref{thm:distributional transform}.} $\;$ In the~\ref{prob:misspecification} problem, given $q\geq 0$ and $G\in\mathcal{A}$, it holds that
\[
\min_{F \in \mathcal{P}} \; \{\mathbb{E}_F[\pi(q,\tilde{u})] + \alpha \cdot d_{\rm OT}(F,G)\} = \mathbb{E}_G\big[\min_{u \geq 0} \; \{\pi(q,u) + \alpha(u-\tilde{v})^2\}\big]  = \mathbb{E}_G[\ell(\alpha,q,\tilde{v})].
\]
Here, the first equality follows from the interchangeability principle (\Cref{lemma:interchangeability}) and the second equality follows from $\ell(\alpha,q,v) = \min_{u \geq 0} \{\pi(q,u) + \alpha(u-v)^2\}$ (see \Cref{lemma:newsvendor type-2} for its closed-form expression). It suffices to verify for any $q \geq 0$ that $\int_{\mathbb{R}_+} \pi(q,v) {\rm d} T_{\varphi_\alpha}[G](v) = \int_{\mathbb{R}_+} \ell(\alpha, q ,v) {\rm d} G(v)~~\forall G \in \mathcal{A}$.
In view of (\textit{i}), that is, $q < \frac{p}{4\alpha}$, we have $\int_{\mathbb{R}_+} \pi(q,v) {\rm d} T_{\varphi_\alpha}[G](v) = \int_{\mathbb{R}_+} \pi(q,v) {\rm d} G\big(\sqrt{\frac{p}{\alpha}v}\big) = \int_{\mathbb{R}_+} \pi\big(q,\frac{\alpha}{p} v^2\big) {\rm d} G(v) = \int_{\mathbb{R}_+} \ell(\alpha, q, v) {\rm d} G(v)$,
where the second equality follows from the variable substitution $v\leftarrow\sqrt{\frac{p}{\alpha}v}$ and the third equality follows from the fact that $\ell(\alpha,q,v) = \min\{\alpha v^2, pq\} - cq$ when $q<\frac{p}{4\alpha}$. 

As for (\textit{ii}), that is, $q\geq \frac{p}{4\alpha}$, we have
\[
\begin{array}{ll}
\int_{\mathbb{R}_+} \pi(q,v) \; {\rm d} T_{\varphi_\alpha}[G](v) & = \int_0^{\frac{p}{4\alpha}}\pi(q, v){\rm d}G\big(\sqrt{\frac{p}{\alpha}v}\big) + \int_{\frac{p}{4\alpha}}^\infty \pi(q, v){\rm d}G\big(v+\frac{p}{4\alpha}\big)\\ 
& = \int_0^{\frac{p}{2\alpha}}\pi\big(q, \frac{\alpha}{p}v^2\big){\rm d}G(v) + \int_{\frac{p}{2\alpha}}^\infty \pi\big(q,  v - \frac{p}{4\alpha}\big){\rm d}G(v) \\
& = \int_0^\infty \ell(\alpha, q, v){\rm d}G(v),
\end{array}
\]
where the second equality follows from the variable substitution: $v \leftarrow \sqrt{\frac{p}{\alpha}v}$ for $v < \frac{p}{4\alpha}$ and $v \leftarrow v + \frac{p}{4\alpha}$ for $v\geq\frac{p}{4\alpha}$ and the third equality follows from the fact that $\ell(\alpha,q,v) = \min\{\alpha v^2, pq\} - cq = \alpha v^2 - cq$ when $v<\frac{p}{2\alpha}$ and $\ell(\alpha,q,v) = p\min\big\{v - \frac{p}{4\alpha}, q\big\} - cq$ when $v\geq\frac{p}{2\alpha}$.
\hfill $\square$
\vspace{2mm}

\noindent{\bf Proof of \Cref{prop:transformed expectation}.} $\;$ Note that for any $q\geq 0$, it holds that $L(q) = \min_{F\in\mathcal{P},\;G\in\mathcal{A}}\;\big\{\mathbb{E}_{F}[\pi(q,\tilde{u})]+\alpha \cdot d_{\rm OT}(F,G)\big\} = \min_{G\in\mathcal{A}}\;\mathbb{E}_{G}[\Psi(\alpha,q,\tilde{v})]$,
where $\Psi(\alpha,q,v) = \pi(q, \varphi_\alpha(v))$. Given $q\geq 0$, $L(q)$ is a moment problem:
\begin{equation}\label{prob:primal}\tag{\sc Primal}
\begin{array}{cll}
\min\limits_{G} & \int_{\mathbb{R}_+} \; \Psi(\alpha, q, v) \; {\rm d}G(v) \\
{\rm s.t.} & \int_{\mathbb{R}_+} \; v \; {\rm d}G(v) = \mu &~~\cdots~s_\alpha \\
& \int_{\mathbb{R}_+}\; v^2 \; {\rm d}G(v) = \mu^2 + \sigma^2 &~~\cdots~r_\alpha \\
& \int_{\mathbb{R}_+}\;{\rm d}G(v) = 1 &~~\cdots~t_\alpha  \\
& G \in \mathcal{M}_+,
\end{array}
\end{equation}
whose dual is 
\begin{equation}\label{prob:dual}\tag{\sc Dual}
\begin{array}{cll}
\displaystyle \max_{s_\alpha,\;r_\alpha,\;t_\alpha} & \mu s_\alpha - (\mu^2+\sigma^2) r_\alpha - t_\alpha \\[1mm]
{\rm s.t.} & v s_\alpha - v^2 r_\alpha - t_\alpha \leq \Psi(\alpha, q, v) &~\forall v \geq 0.
\end{array}
\end{equation}
We next derive the expression of $L(q)$ by constructing a pair of primal and dual feasible solutions that attain the same objective value (\textit{i.e.}, strong duality holds between \ref{prob:primal} and \ref{prob:dual}). We first find a dual feasible solution $(s_\alpha, r_\alpha, t_\alpha)$ such that the quadratic function $g(v) = v s_\alpha - v^2 r_\alpha - t_\alpha$ touches $\Psi(\alpha, q, v)$ at exactly two points, where there are three scenarios in total (as shown in Figure~\ref{fig:dual tangent}). We then find a primal feasible solution $G_\alpha \in \mathcal{A}$ whose objective value matches the dual objective under $(s_\alpha, r_\alpha, t_\alpha)$. This is achieved by setting the support of $G_\alpha$ as the tangent points of $g(v)$ and $\Psi(\alpha, q, v)$ and then solving for the corresponding probabilities based on the moment conditions in $\mathcal{A}$. In the following, we consider these three scenarios based on the value of $q$.

\begin{figure}[tb]
\begin{subfigure}{0.33\textwidth}
\centering
\includegraphics[width=1.11\linewidth]{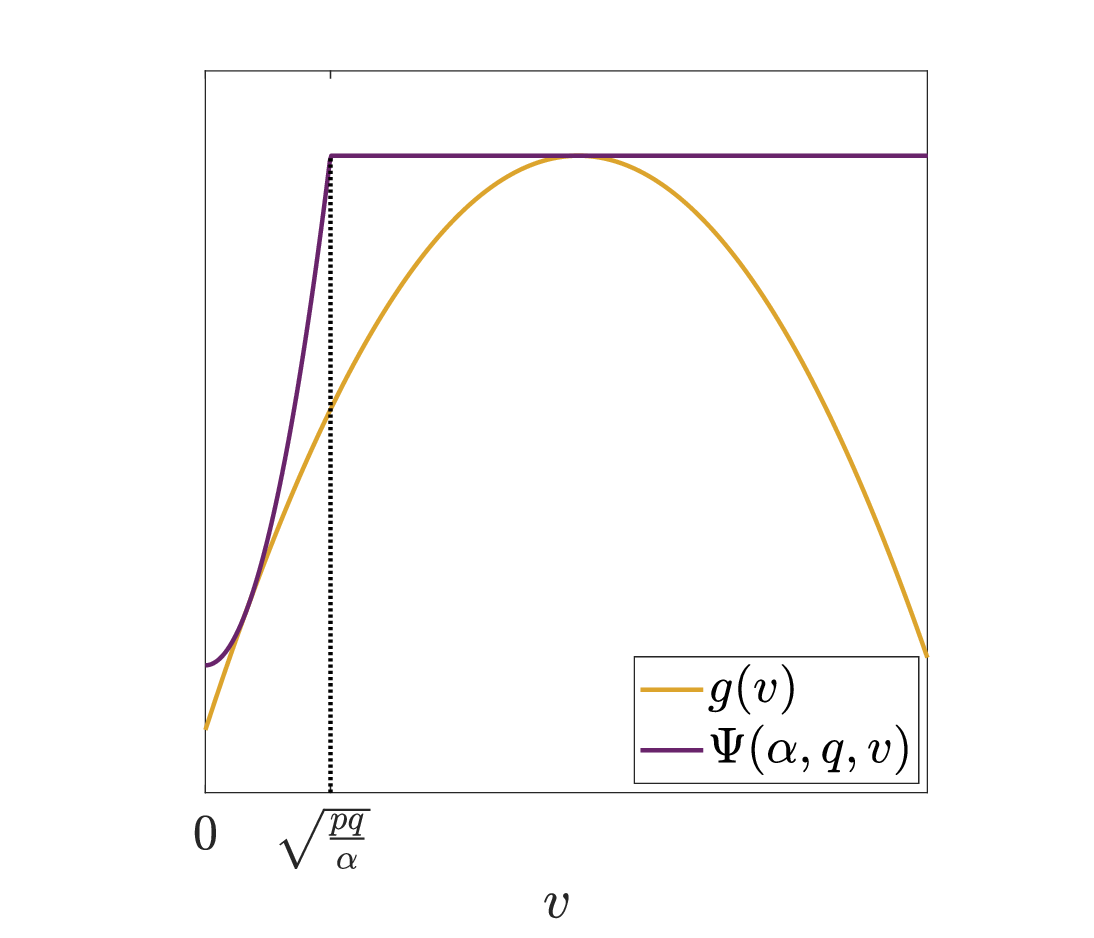}
\end{subfigure}
\begin{subfigure}{0.33\textwidth}
\centering
\includegraphics[width=1.11\linewidth]{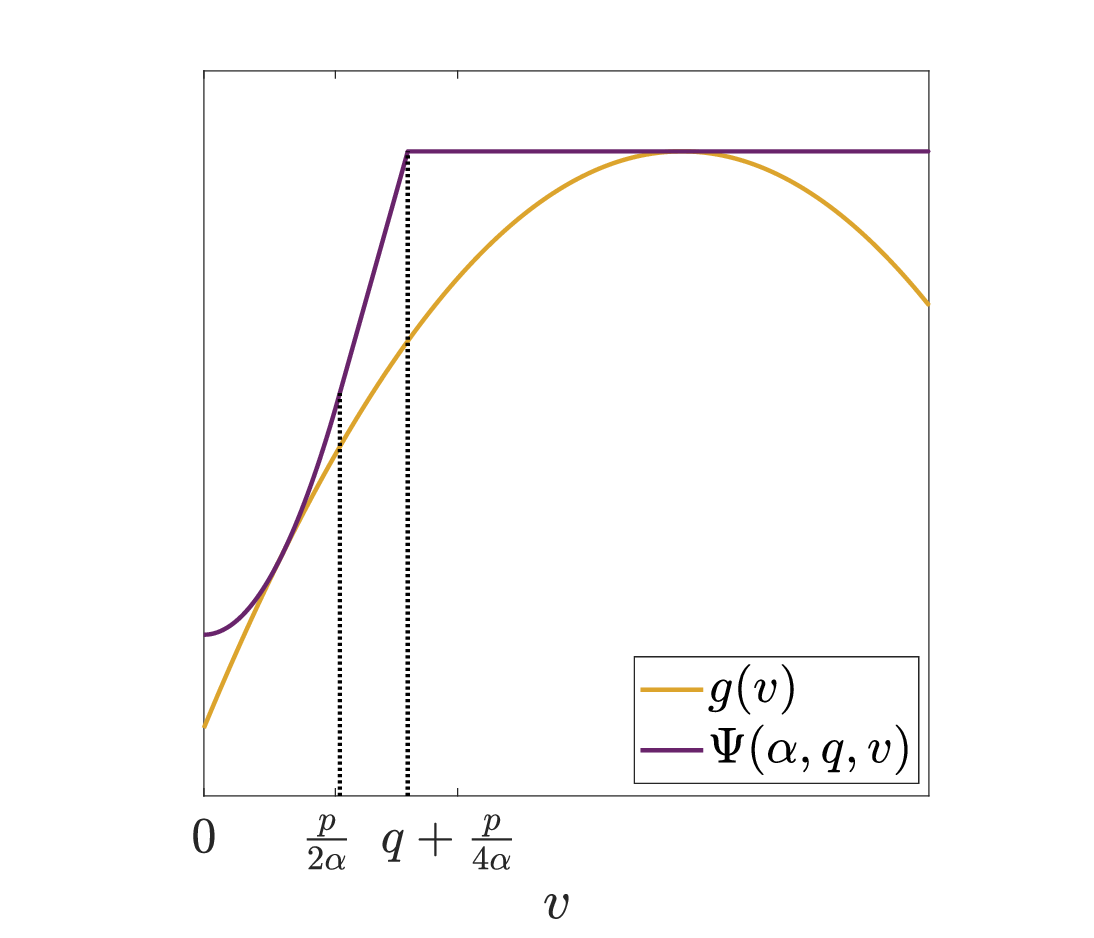}
\end{subfigure}
\begin{subfigure}{0.33\textwidth}
\centering
\includegraphics[width=1.11\linewidth]{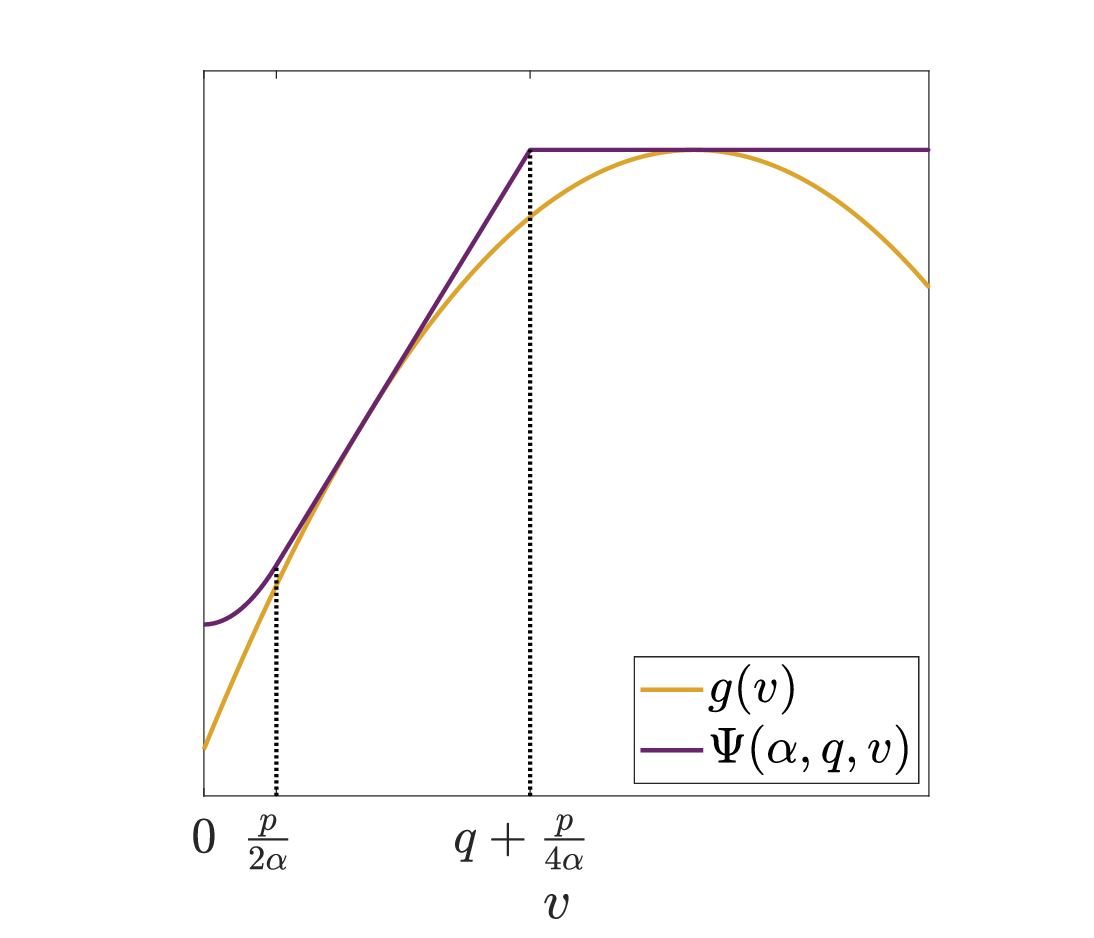}
\end{subfigure}
\caption{\textnormal{Illustration of the dual problem. Here, $g(v) = vs_\alpha - v^2r_\alpha - t_\alpha$.}}
\label{fig:dual tangent}
\vspace{-6mm}
\end{figure}

\textit{Scenario~1}. $\;$ 
When $0 \leq q \leq \frac{p}{4\alpha}$, we first construct a feasible distribution to \ref{prob:primal} as follows:
\begin{equation}\label{equ:primal solution 1}
\textstyle G_\alpha = \bigg(\frac{1}{2} - \frac{\mu^2 - \sigma^2 - \frac{pq}{\alpha}}{2\sqrt{(\frac{pq}{\alpha} + \mu^2 + \sigma^2)^2 - 4\mu^2\frac{pq}{\alpha}}}\bigg) \cdot \delta_{v_1} + \bigg(\frac{1}{2} + \frac{\mu^2 - \sigma^2 - \frac{pq}{\alpha}}{2\sqrt{(\frac{pq}{\alpha} + \mu^2 + \sigma^2)^2 - 4\mu^2\frac{pq}{\alpha}}}\bigg) \cdot \delta_{v_2},
\end{equation}
where the support points are 
\[
\begin{array}{ll}
\textstyle v_1 = \frac{1}{2\mu} \big(\frac{pq}{\alpha} + \mu^2 + \sigma^2 - \sqrt{(\frac{pq}{\alpha} + \mu^2 + \sigma^2)^2 - 4\mu^2\frac{pq}{\alpha}}\big) 
~{\rm and}~
v_2 = \frac{1}{2\mu} \big(\frac{pq}{\alpha} + \mu^2 + \sigma^2 + \sqrt{(\frac{pq}{\alpha} + \mu^2 + \sigma^2)^2 - 4\mu^2\frac{pq}{\alpha}}\big).
\end{array}
\]
One can verify that $G_\alpha \in \mathcal{A}$ and the corresponding primal objective value under $G_\alpha$ is equal to $\mathbb{E}_{G_\alpha}[\Psi(\alpha,q,\tilde{u})] = \frac{\alpha}{2}\big(\frac{pq}{\alpha}+\mu^2+\sigma^2-\sqrt{(\frac{pq}{\alpha}+\mu^2+\sigma^2)^2-4\mu^2\frac{pq}{\alpha}}\big)-cq$. 
We next construct a dual feasible solution that attains the same dual objective value. Note that $\Psi(\alpha, q, v) = \min\{\alpha v^2, pq\} - cq ~\forall v \geq 0$ when $0 \leq q \leq \frac{p}{4\alpha}$ (\Cref{lemma:newsvendor type-2}). Hence, \ref{prob:dual} becomes
\[
\begin{array}{cll}
\displaystyle \max_{s_\alpha,\;r_\alpha,\;t_\alpha} & \displaystyle  s_\alpha \mu - r_\alpha (\mu^2+\sigma^2) - t_\alpha \\ 
{\rm s.t.} & \displaystyle s_\alpha v - r_\alpha v^2 -t_\alpha \leq \alpha v^2 - cq &~~~ \displaystyle \forall v \geq 0\\ 
& \displaystyle s_\alpha v - r_\alpha v^2 - t_\alpha \leq pq - cq &~~~ \displaystyle \forall v \geq 0.
\end{array}
\]
Consider the following solution:
\begin{equation}\label{equ:dual solution 1}\mspace{-20mu}
\begin{array}{ll}
s_\alpha = \frac{2\mu pq}{\sqrt{(\frac{pq}{\alpha} + \mu^2 + \sigma^2)^2 - 4\mu^2\frac{pq}{\alpha}}},
r_\alpha = \frac{\alpha}{2} \bigg(\frac{\frac{pq}{\alpha} + \mu^2 + \sigma^2}{\sqrt{(\frac{pq}{\alpha} + \mu^2 + \sigma^2)^2 - 4\mu^2\frac{pq}{\alpha}}} - 1\bigg),
t_\alpha = \frac{pq}{2} \bigg(\frac{\frac{pq}{\alpha} + \mu^2 + \sigma^2}{\sqrt{(\frac{pq}{\alpha} + \mu^2 + \sigma^2)^2 - 4\mu^2\frac{pq}{\alpha}}} - 1\bigg) + cq,
\end{array}
\end{equation}
which satisfies $s_\alpha \mu - r_\alpha (\mu^2 + \sigma^2) - t_\alpha = \frac{\alpha}{2}\big(\frac{pq}{\alpha}+\mu^2+\sigma^2-\sqrt{(\frac{pq}{\alpha}+\mu^2+\sigma^2)^2-4\mu^2\frac{pq}{\alpha}}\big)-cq$.
It remains to argue that this solution is feasible to \ref{prob:dual}. If $q = 0$, then $s_\alpha = r_\alpha = t_\alpha = 0$, naturally feasible to \ref{prob:dual}. We next investigate $q > 0$. In this case, $r_\alpha > 0$. The first semi-infinite constraint of \ref{prob:dual} is equivalent to $\max_{v \geq 0} \{ s_\alpha v - (r_\alpha + \alpha) v^2 - t_\alpha + cq\} \leq 0$. For the left-hand side maximization, the optimal solution is $v^\star = \frac{s_\alpha}{2(r_\alpha + \alpha)} \geq 0$, which attains an optimal value of $\frac{s_\alpha^2}{4(r_\alpha + \alpha)} - t_\alpha + cq = 0$. Hence, the first semi-infinite constraint is satisfied. Similarly, the second semi-infinite constraint of \ref{prob:dual} is equivalent to $\max_{v \geq 0} \{ s_\alpha v - r_\alpha v^2 - t_\alpha - pq + cq\} \leq 0$.
For the left-hand side, the optimal solution is $v^\star = \frac{s_\alpha}{2r_\alpha}$ and the corresponding optimal value is $\frac{s^2_\alpha}{4r_\alpha} - t_\alpha - pq + cq = \frac{pq(r_\alpha + \alpha)}{\alpha} - \frac{pq}{\alpha}r_\alpha - pq = 0$. Hence, the second semi-infinite constraint is also satisfied, concluding that solution~\eqref{equ:dual solution 1} is feasible to \ref{prob:dual} and establishing the strong duality.

\textit{Scenario~2}. $\;$ 
When $q\geq \frac{p}{4\alpha}$ and $(2\mu-\frac{p}{\alpha})q < \mu^2 + \sigma^2 - \frac{p\mu}{2\alpha}$, \ref{prob:dual} becomes
\[
\begin{array}{cll}
\displaystyle \max_{s_\alpha,\;r_\alpha,\;t_\alpha} & \displaystyle  s_\alpha \mu - r_\alpha (\mu^2+\sigma^2) -t_\alpha \\ [1mm]
{\rm s.t.} & s_\alpha v - r_\alpha v^2 - t_\alpha \leq \alpha v^2 - cq &~~~ \forall 0\leq v \leq \frac{p}{2\alpha}\\ 
& s_\alpha v - r_\alpha v^2 - t_\alpha \leq p\big(v-\frac{p}{4\alpha}\big) - cq &~~~\forall v \geq \frac{p}{2\alpha}\\ 
& s_\alpha v - r_\alpha v^2 - t_\alpha \leq pq - cq &~~~ \forall v \geq \frac{p}{2\alpha}.
\end{array}
\]
Consider the pair of primal feasible solution~\eqref{equ:primal solution 1} and dual solution~\eqref{equ:dual solution 1}. Upon the results in {\it Step~1}, we can establish the feasibility of the constructed dual solution similar to {\it Scenario 1}.

\textit{Scenario~3}. $\;$ 
When $q\geq \frac{p}{4\alpha}$ and $(2\mu-\frac{p}{\alpha})q\geq\mu^2+\sigma^2 - \frac{p\mu}{2\alpha}$, we first construct a feasible primal solution:
\begin{equation}\label{equ:primal solution 2}
\textstyle G_\alpha = \frac{1}{2}\bigg(1+\frac{q+\frac{p}{4\alpha}-\mu}{\sqrt{(q+\frac{p}{4\alpha}-\mu)^2+\sigma^2}}\bigg)\cdot\delta_{v_1} + \frac{1}{2}\bigg(1-\frac{q+\frac{p}{4\alpha}-\mu}{\sqrt{(q+\frac{p}{4\alpha}-\mu)^2+\sigma^2}}\bigg)\cdot\delta_{v_2},
\end{equation}
where $v_1 = q + \frac{p}{4\alpha} - \sqrt{(q + \frac{p}{4\alpha}-\mu)^2+\sigma^2}~~{\rm and}~~ v_2 = q + \frac{p}{4\alpha} + \sqrt{(q + \frac{p}{4\alpha}-\mu)^2+\sigma^2}$.
One can verify that $G_\alpha \in \mathcal{A}$ and $\mathbb{E}_{G_\alpha}[\Psi(\alpha,q,\tilde{u})] = \frac{p}{2}\big(\mu - q - \frac{p}{4\alpha} - \sqrt{(q + \frac{p}{4\alpha}-\mu)^2 + \sigma^2}\big) + (p-c)q$.
We next construct a dual feasible solution that attains the same dual objective value. Consider the following solution:
\begin{equation}\label{equ:dual solution 2}
\textstyle s_\alpha = \frac{p}{2} + 2r_\alpha\big(q + \frac{p}{4\alpha}\big), 
t_\alpha = \frac{p^2}{16r_\alpha}+r_\alpha\big(q+\frac{p}{4\alpha}\big)^2 + \frac{p}{2}\big(q+\frac{p}{4\alpha}\big) - (p-c)q, r_\alpha = \frac{p}{4\sqrt{(q + \frac{p}{4\alpha}-\mu)^2 + \sigma^2}},
\end{equation}
which is feasible to \ref{prob:dual} and attains the same primal objective value under $G_\alpha$.

Note that {\it Scenario 1} and {\it Scenario 2} correspond to the same objective function. Combining the results in these three scenarios, we then obtain the desired result immediately. 
\hfill $\square$
\vspace{2mm}

\noindent{\bf Proof of \Cref{thm:newsvendor transport-2}.} $\;$ Defining $L(q) = \min_{G\in\mathcal{A}} \mathbb{E}_{G}[\Psi(\alpha,q,\tilde{v})]$, \ref{prob:misspecification} is then equivalent to $\max_{q\geq 0} L(q)$.
Note that when $\kappa < \frac{\sigma^2}{\mu^2+\sigma^2}$, $\max_{q\geq 0}\min_{G\in\mathcal{A}}\mathbb{E}_{G}[\Psi(\alpha,q,\tilde{v})] \leq \max_{q\geq 0}\min_{G\in\mathcal{A}}\mathbb{E}_{G}[\pi(q,\tilde{v})] = 0$,
where the inequality holds since $\Psi(\alpha,q,v) = \pi(q,\varphi_\alpha(v)) \leq \pi(q,v) ~\forall q\geq 0$ and the equality follows from the tight lower bound derived in \cite{scarf1958min1}. The order quantity $q = 0$ satisfies $\mathbb{E}_{G}[\Psi(\alpha, 0,\tilde{v})] = 0 ~\forall G\in\mathcal{A}$ (\textit{i.e.}, $\min_{G\in\mathcal{A}}\mathbb{E}_{G}[\Psi(\alpha, 0,\tilde{v})] = 0$). Hence, $q_\alpha^\star = 0$. This corresponds to case~(\textit{ii}) in the statement. In the remainder of the proof, we focus on $\kappa \geq \frac{\sigma^2}{\mu^2+\sigma^2}$. As shown in \Cref{prop:transformed expectation}, we have
\[
L(q) = \begin{cases}
\frac{p}{2}\big(\mu - q - \frac{p}{4\alpha} - \sqrt{(q + \frac{p}{4\alpha}-\mu)^2 + \sigma^2}\big) + (p-c)q & {\rm if}~ q \in \mathcal{Q}\\[2mm]
\frac{\alpha}{2}\big(\frac{pq}{\alpha}+\mu^2+\sigma^2-\sqrt{(\frac{pq}{\alpha}+\mu^2+\sigma^2)^2-4\mu^2\frac{pq}{\alpha}}\big)-cq & {\rm otherwise}
\end{cases}
\]
where $\mathcal{Q}= \{q \in \mathbb{R}_+ \mid q \geq \frac{p}{4\alpha}, \;(2\mu-\frac{p}{\alpha})q \geq \mu^2+\sigma^2 - \frac{p\mu}{2\alpha}\}$.
Based on the value of $\alpha$, we divide the arguments into three different scenarios.

\textit{Scenario~1}. $\;$
Suppose that $\alpha < \frac{p}{2\mu}$, \textit{i.e.}, $2\mu - \frac{p}{\alpha} < 0$. For any $q \geq \frac{p}{4\alpha}$, it holds that
$(2\mu-\frac{p}{\alpha}) q \leq (2\mu-\frac{p}{\alpha})\frac{p}{4\alpha} = \frac{p}{\alpha}(\mu - \frac{p}{4\alpha}) - \frac{p\mu}{2\alpha} < \mu^2 + \sigma^2 - \frac{p\mu}{2\alpha}$, where the strict inequality follows from the fact that $\frac{p}{\alpha}(\mu - \frac{p}{4\alpha}) < \mu^2$ when $\alpha < \frac{p}{2\mu}$. That is to say, for any $q \geq 0$ we have $L(q) = \frac{\alpha}{2}\big(\frac{pq}{\alpha}+\mu^2+\sigma^2 -\sqrt{(\frac{pq}{\alpha}+\mu^2+\sigma^2)^2-4\mu^2\frac{pq}{\alpha}}\big) - cq$. Setting the derivative of $L(q)$ to $0$, we then obtain $q^\star_\alpha = \big(\mu^2 - \sigma^2 + 2f(1-\kappa)\mu\sigma\big) \cdot \frac{\alpha}{p}$.

\textit{Scenario~2}. $\;$ 
Suppose that $\frac{p}{2\mu}\leq \alpha< \frac{p}{2(\mu-\sigma\sqrt{(1-\kappa)/\kappa})}$. Note that $(2\mu - \frac{p}{\alpha})\frac{p}{4\alpha} + \frac{p\mu}{2\alpha} = \frac{p}{\alpha}(\mu - \frac{p}{2\alpha}) \leq \mu^2 + \sigma^2$, where the inequality is due to the fact that $\frac{p}{\alpha}(\mu - \frac{p}{2\alpha}) \leq \frac{\mu^2}{2}$. Hence, $\frac{p}{4\alpha} \leq \frac{\mu}{2} + \frac{\sigma^2}{2\mu - \frac{p}{\alpha}}$, yielding
\begin{equation}\label{equ:objective 2}
L(q) = 
\left\{
\begin{array}{ll}
\frac{\alpha}{2}\big(\frac{pq}{\alpha}+\mu^2+\sigma^2-\sqrt{(\frac{pq}{\alpha}+\mu^2+\sigma^2)^2-4\mu^2\frac{pq}{\alpha}}\big)-cq &~~~ q \leq \frac{\mu}{2} + \frac{\sigma^2}{2\mu - \frac{p}{\alpha}}\\ [1mm]
\frac{p}{2}\big(\mu - q - \frac{p}{4\alpha} - \sqrt{(q + \frac{p}{4\alpha}-\mu)^2 + \sigma^2}\big) + (p-c)q &~~~ q \geq \frac{\mu}{2} + \frac{\sigma^2}{2\mu - \frac{p}{\alpha}}.
\end{array}
\right.
\end{equation}
Setting $L'(q)$ to $0$ then yields $q^\star_\alpha = \big(\mu^2 - \sigma^2 + 2f(1-\kappa)\mu\sigma\big) \cdot \frac{\alpha}{p} \leq \frac{\mu}{2} + \frac{\sigma^2}{2\mu - \frac{p}{\alpha}}$,
where the inequality follows from the condition $\alpha <  \frac{p}{2(\mu-\sigma\sqrt{(1-\kappa)/\kappa})}$.

\textit{Scenario~3}. $\;$ 
Suppose that $\alpha \geq \frac{p}{2(\mu-\sigma\sqrt{(1-\kappa)/\kappa})}$. Then $L(q)$ is given in \eqref{equ:objective 2}, which is concave. Setting $L'(q^\star_\alpha) = 0$ yields $q^\star_\alpha = \mu + \sigma f(1-\kappa) - \frac{p}{4\alpha}\geq \frac{\mu}{2} + \frac{\sigma^2}{2\mu - \frac{p}{\alpha}}$.

By noting that {\it Scenario~1} and {\it Scenario~2} correspond to $\alpha <\frac{p}{2(\mu - \sigma \sqrt{(1-\kappa)/\kappa})} $ and {\it Scenario~3} corresponds to $\alpha \geq \frac{p}{2(\mu - \sigma \sqrt{(1-\kappa)/\kappa})}$, we then complete the proof.
\hfill $\square$
\vspace{2mm}

\noindent{\bf Proof of~\Cref{prop:price}.} $\;$ When $p\geq \max\big\{\frac{\mu^2+\sigma^2}{\mu^2}c, 2\alpha\mu\big\}$, it holds that $q_\alpha^\star = \big(\mu^2 - \sigma^2 + \mu\sigma\frac{1-2\kappa}{\sqrt{\kappa(1-\kappa)}}\big)\cdot \frac{\alpha}{p} = \frac{\alpha}{c} \cdot \frac{(\mu^2-\sigma^2+\mu\sigma(1-1/x))}{(x^2+1)}$,
where we denote $x = \sqrt{\kappa/(1-\kappa)}$ to have the second equality. The  derivative of $q_\alpha^\star$ with respect to $x$ is $\frac{\partial q_\alpha^\star}{\partial x} = \frac{\alpha}{c} \cdot \frac{4\mu\sigma-\mu\sigma x^2-2x(\mu^2-\sigma^2)+\mu\sigma/x^2}{(x^2+1)^2}$.
To determine the sign of $\frac{\partial q_\alpha^\star}{\partial x}$, it suffices to focus on the term $4\mu\sigma-\mu\sigma x^2-2x(\mu^2-\sigma^2)+\frac{\mu\sigma}{x^2}$. Note that when $x\geq 1$ ({\it i.e.}, $\kappa>\frac{1}{2}$), we have $4\mu\sigma-\mu\sigma x^2-2x(\mu^2-\sigma^2)+\frac{\mu\sigma}{x^2}\leq 5\mu\sigma-\mu\sigma x^2-2x(\mu^2-\sigma^2) = -\mu\sigma(x+\frac{\mu^2-\sigma^2}{\mu\sigma})^2+5\mu\sigma-\frac{(\mu^2-\sigma^2)}{\mu\sigma}$.
Since $\mu\sigma>0$, there must exist some $x_1\geq 1$ such that when $x> x_1$, it holds that $4\mu\sigma-\mu\sigma x^2-2x(\mu^2-\sigma^2)+\frac{\mu\sigma}{x^2}\leq  -\mu\sigma(x+\frac{\mu^2-\sigma^2}{\mu\sigma})^2+5\mu\sigma-\frac{(\mu^2-\sigma^2)}{\mu\sigma}< 0$.
As $x$ is strictly increasing in ($\kappa$ and) $p$, it is clear that $q_\alpha^\star$ is strictly decreasing in $p$ when $p > p_\alpha^\star = \max\{\frac{\mu^2+\sigma^2}{\mu^2}c, 2\alpha\mu, c(x_1^2+1)\}$,
completing the proof.
\hfill $\square$
\vspace{2mm}

\noindent{\bf Proof of~\Cref{prop:variance}.} $\;$
The derivative of $q_\alpha^\star$ with respect to $\sigma$ is
\[
\frac{\partial q_\alpha^\star}{\partial \sigma} = \left\{
\begin{array}{ll}
 f(1-\kappa) &~~~ \sigma < (\mu-\frac{p}{2\alpha})\sqrt{\frac{\kappa}{1-\kappa}}\\ [1mm]
\big( - 2\sigma + 2\mu f(1-\kappa)\big)\cdot \frac{\alpha}{p} &~~~ (\mu-\frac{p}{2\alpha})\sqrt{\frac{\kappa}{1-\kappa}}\leq\sigma\leq \mu\sqrt{\frac{\kappa}{1-\kappa}},
\end{array}
\right.
\]
which is decreasing in $\sigma$. When $\kappa \geq \frac{1}{2}$, $\frac{\partial q_\alpha^\star}{\partial \sigma} = f(1-\kappa)\geq 0~~\forall \sigma \in \big[0,(\mu - \frac{p}{2\alpha})\sqrt{\frac{\kappa}{1-\kappa}}\big]$. We proceed by dividing the argument into two cases. On the one hand, if $\alpha\geq\frac{p-c}{\mu}$, we have $\mu f(1-\kappa) = \frac{\mu}{2}(\sqrt{\frac{\kappa}{1-\kappa}} - \sqrt{\frac{1-\kappa}{\kappa}}) \leq (\mu - \frac{p}{2\alpha})\sqrt{\frac{\kappa}{1-\kappa}}$, which indicates that $\frac{\partial q_\alpha^\star}{\partial \sigma} = ( - 2\sigma + 2\mu f(1-\kappa))\cdot \frac{\alpha}{p} \leq 0~~\forall \sigma \in [(\mu - \frac{p}{2\alpha})\sqrt{\frac{\kappa}{1-\kappa}},\mu\sqrt{\frac{\kappa}{1-\kappa}}]$. On the other hand, if $\alpha<\frac{p-c}{\mu}$, it holds that $\mu f(1-\kappa)  > (\mu - \frac{p}{2\alpha})\sqrt{\frac{\kappa}{1-\kappa}}$, which implies that $\frac{\partial q_\alpha^\star}{\partial \sigma} \geq 0~~\forall \sigma \in [(\mu - \frac{p}{2\alpha})\sqrt{\frac{\kappa}{1-\kappa}},\mu f(1-\kappa)]$ and $\frac{\partial q_\alpha^\star}{\partial \sigma} \leq 0~~\forall \sigma \in [\mu f(1-\kappa),\mu\sqrt{\frac{\kappa}{1-\kappa}}]$.
Combining the two cases together then yields the desired result.
\hfill $\square$
\vspace{2mm}

\noindent{\bf Proof of \Cref{lemma:distance}.} $\;$
We proceed with the proof by dividing the argument into two scenarios.

{\it Scenario 1.1.} If $\frac{\hat{\mu}}{\hat{\sigma}} \geq \frac{\mu}{\sigma}$, by theorem 2.1 in \cite{gelbrich1990formula2}, for any $G \in \mathcal{A}_N$, it holds that $d_{\rm OT}(D,G) \geq (\hat{\mu} - \mu)^2 + (\hat{\sigma} - \sigma)^2$, and hence $d(D, \mathcal{A}_N) \geq (\hat{\mu} - \mu)^2 + (\hat{\sigma} - \sigma)^2$. In the following, we show that the inequality is tight. Consider the transformed random variable $\tilde{w} = \frac{\hat{\sigma}}{\sigma}\tilde{u} + \hat{\mu} - \frac{\hat{\sigma}\mu}{\sigma} \sim G^\dag$ where $\tilde{u}\sim D$. Since the $\tilde{w}$ is a linear transformation of $\tilde{u}$, we have $d_{\rm OT}(D, G^\dag) = (\hat{\mu} - \mu)^2 + (\hat{\sigma} - \sigma)^2$ (\citealt{dowson1982frechet}). Moreover, we can verify $\mathbb{E}_{G^\dag}[\tilde{w}] = \hat{\mu}$, $\mathbb{E}_{G^\dag}[\tilde{w}^2] = \hat{\mu}^2 + \hat{\sigma}^2$, and $G^\dag\{\tilde{w} \in [0,+\infty)\} = 1$ (since $\frac{\hat{\mu}}{\hat{\sigma}} \geq \frac{\mu}{\sigma}$). Hence, $G^\dag \in \mathcal{A}_N$, concluding $d(D, \mathcal{A}_N) = (\hat{\mu} - \mu)^2 + (\hat{\sigma} - \sigma)^2$.

{\it Scenario 1.2.} If $\frac{\hat{\mu}}{\hat{\sigma}} < \frac{\mu}{\sigma}$, we construct an upper bound for $d(D, \mathcal{A}_N)$. Consider the transformed random variable $\tilde{w} = k^\dag\max\{0, \tilde{u} - t^\dag\} \sim G^\dag$ with $\tilde{u}\sim D$, where $k^\dag > 0$ and $t^\dag \geq 0$ satisfy
\begin{equation}\label{EC:moment constraints}
\textstyle \int_{t^\dag}^{+\infty}k^\dag(u - t^\dag){\rm d}D(u) = \hat{\mu}
~~~{\rm and}~~~
\int_{t^\dag}^{+\infty}k^{\dag2}(u - t^\dag)^2{\rm d}D(u) = \hat{\mu}^2 + \hat{\sigma}^2.
\end{equation}
In the following, we first show that there exists $(k^\dag, t^\dag)$ satisfying~\eqref{EC:moment constraints}. Eliminating the variable $k$ in~\eqref{EC:moment constraints}, it suffices to check whether there exists $t \geq 0$ such that
\[
\textstyle \hat{\mu}^2\int_{t}^{+\infty}(u - t)^2{\rm d}D(u) - (\hat{\mu}^2 + \hat{\sigma}^2)(\int_{t}^{+\infty}(u - t){\rm d}D(u))^2 = 0.
\]
Define $h(t) = \hat{\mu}^2\int_t^{+\infty}(u - t)^2{\rm d}D(u) - (\hat{\mu}^2 + \hat{\sigma}^2)(\int_t^{+\infty}(u - t){\rm d}D(u))^2$. Setting the derivative $h'(t) = 2\int_t^{+\infty}(u - t){\rm d}D(u)((\hat{\mu}^2 + \hat{\sigma}^2)(1-D(t)) - \hat{\mu}^2)$ to $0$ then yields $t^\diamond = D^{-1}(\frac{\hat{\sigma}^2}{\hat{\mu}^2 + \hat{\sigma}^2})$. It is straightforward to see that $h(0) = \hat{\mu}^2 (\mu^2 + \sigma^2) - \mu^2 (\hat{\mu}^2 + \hat{\sigma}^2) < 0$ and $\lim_{t\to+\infty}h(t) = 0$, which implies that $h(t^\diamond) > \lim_{t\to+\infty}h(t) = 0$ since $h(t)$ is decreasing in $(t^\diamond, +\infty)$. Therefore, there must exist some $t^\dag \in [0, t^\diamond]$ such that $h(t^\dag) = 0$, verifying the feasibility of~\eqref{EC:moment constraints}. This indicates that $\mathbb{E}_{G^\dag}[\tilde{w}] = \hat{\mu}$ and $\mathbb{E}_{G^\dag}[\tilde{w}^2] = \hat{\mu}^2 + \hat{\sigma}^2$. Additionally, since $G^\dag\{\tilde{w} \in [0,+\infty)\} = 1$, it is immediate to see that $G^\dag \in \mathcal{A}_N$. Subsequently, we identify the upper bounds for $k^\dag$ and $t^\dag$. Note that for any $t \in [0, t^\diamond]$, we have $h''(t) = 2\hat{\mu}^2(1-D(t)) - 2d(t)(\hat{\mu}^2 + \hat{\sigma}^2)\int_t^{+\infty}(u - t){\rm d}D(u) - 2(1-D(t))^2(\hat{\mu}^2 + \hat{\sigma}^2) \leq 0$, where $d(t)$ is the density function of the distribution $D$. Hence, it holds that
\[
\textstyle \frac{h(t^\dag)-h(0)}{t^\dag - 0} \geq \frac{h(t^\diamond) - h(t^\dag)}{t^\diamond - t^\dag}~\Longrightarrow~t^\dag \leq \frac{-h(0)t^\diamond}{h(t^\diamond) - h(0)}.
\]
Note that as $N \to +\infty$, $h(0)\to 0$ and hence $t^\dag \to 0$. This implies that for sufficiently large $N$, $t^\dag \leq \mu$. Since $\int_{t^\dag}^{+\infty}(u - t^\dag){\rm d}D(u)  = \mu - \int_0^{t^\dag}u{\rm d}D(u) - t^\dag (1 - D(t^\dag)) \geq \mu - t^\dag (1 - D(0))$, we then have
\[
\textstyle \frac{\hat{\mu}}{k^\dag} \geq \mu - t^\dag (1 - D(0)) ~\Longrightarrow~ k^\dag \leq \frac{\hat{\mu}}{\mu - t^\dag (1 - D(0))} \leq \frac{\hat{\mu}}{\mu + \frac{h(0)t^\diamond}{h(t^\diamond) - h(0)}} \leq \frac{\hat{\mu}}{\mu + \frac{h(0)}{h'(0)}},
\]
where the last inequality follows from the fact that $\frac{h(t^\diamond) - h(0)}{t^\diamond} \leq h'(0)$ since $h''(t) \leq 0$ for any $t \in [0, t^\diamond]$. Plugging the expressions of $h(0)$ and $h'(0)$, it is then immediate to see that $k^\dag \leq \frac{\hat{\mu}}{\mu + \frac{h(0)}{h'(0)}} = \frac{2\mu^2\hat{\sigma}^2}{\mu^2\hat{\sigma}^2 + \hat{\mu}^2\sigma^2}$.
For the optimal-transport cost, note that the objective of the Kantorovich formulation as defined in~\eqref{equ:optimal transport 2} is no larger than that of the Monge formulation (\citealp{villani2009optimal1}), {\it i.e.}, $
\textstyle d_{\rm OT}(D, G^\dag) \leq \inf_{\psi: T_{\psi}[D] = G^\dag}\;\int_0^{+\infty}(u - \psi(u))^2{\rm d}D(u)$. Note that the function $\psi^\dag(u) = k^\dag\cdot \max\{u - t^\dag,0\}$ is feasible to the right-hand side problem. Hence,
\[
\begin{array}{rl}
d(D,\mathcal{A}_N) \leq d_{\rm OT}(D,G^\dag) & \leq \int_{0}^t u^2 {\rm d}D(u) + \int_{t}^{+\infty} (k^\dag(u - t^\dag) - u)^2 {\rm d}D(u)\\ [1mm]
& = \mu^2 + \sigma^2 + \hat{\mu}^2 + \hat{\sigma}^2 - 2k^\dag\int_{t^\dag}^{+\infty}(u - t^\dag)u {\rm d}D(u)\\ [1mm]
& \leq \mu^2 + \sigma^2 + \hat{\mu}^2 + \hat{\sigma}^2 - \frac{(\mu^2\hat{\sigma}^2 + \hat{\mu}^2\sigma^2)(\hat{\mu}^2 + \hat{\sigma}^2)}{\mu^2\hat{\sigma}^2}\\ [1mm]
& = (\mu - \hat{\mu})^2 + (\sigma - \hat{\sigma})^2 + 2\mu\hat{\mu} + 2\sigma\hat{\sigma} - \frac{(\mu^2\hat{\sigma}^2 + \hat{\mu}^2\sigma^2)(\hat{\mu}^2 + \hat{\sigma}^2)}{\mu^2\hat{\sigma}^2},
\end{array}
\]
where the third line follows from $u \geq u-t^\dag$ and $k^\dag \leq \frac{2\mu^2\hat{\sigma}^2}{\mu^2\hat{\sigma}^2 + \hat{\mu}^2\sigma^2}$. Besides,
\[
\begin{array}{rl}
2\mu\hat{\mu} + 2\sigma\hat{\sigma} - \frac{(\mu^2\hat{\sigma}^2 + \hat{\mu}^2\sigma^2)(\hat{\mu}^2 + \hat{\sigma}^2)}{\mu^2\hat{\sigma}^2} & = \frac{\hat{\mu}(\mu^2\hat{\sigma}^2 - \hat{\mu}^2\sigma^2)}{\mu\hat{\sigma}^2} + \frac{\sigma(\mu\hat{\sigma} - \hat{\mu}\sigma)}{\mu} + \frac{\hat{\sigma}(\hat{\mu}\sigma - \mu\hat{\sigma})}{\hat{\mu}} \leq  \frac{\mu^2\hat{\sigma}^2 - \hat{\mu}^2\sigma^2}{\sigma\hat{\sigma}},
\end{array}
\]
where the inequality is due to $\frac{\hat{\mu}}{\hat{\sigma}} < \frac{\mu}{\sigma}$. Hence, we have $d(D,\mathcal{A}_N) \leq (\mu - \hat{\mu})^2 + (\sigma - \hat{\sigma})^2 + \frac{\mu^2\hat{\sigma}^2 - \hat{\mu}^2\sigma^2}{\sigma\hat{\sigma}}$. Combining the results in these two scenarios then completes the proof. 
\hfill $\square$
\vspace{2mm}

\noindent{\bf Proof of \Cref{prop:concentration}.} $\;$ We proceed in two steps. In the first step, we derive the concentration inequalities of the sample mean and variance, respectively. In the second step, we establish the concentration inequality for the mean-variance ambiguity set.

{\it Step 1.} Note that for any $x \in \mathbb{R}$, $\mathbb{E}_{D^N}[{\rm exp}(\frac{x}{N}\sum_{i=1}^N(\hat{v}_i - \mu))] = \prod_{i=1}^N\mathbb{E}_{D}[{\rm exp}(\frac{x}{N}(\hat{v}_i - \mu))] \leq \mathbb{E}_{D}[{\rm exp}(-\frac{x^2\nu^2}{2N})]$,
where the expectation is taken with respect to the random sample $\hat{v}_i$, the equality follows from the fact that $\hat{v}_1,\ldots,\hat{v}_N$ are i.i.d, and the inequality follows from the fact that $D$ is sub-Gaussian with variance proxy $\nu^2$. Hence, $\hat{\mu} - \mu = \frac{1}{N}\sum_{i=1}^N(\hat{v}_i - \mu)$ is sub-Gaussian with variance proxy $\frac{\nu^2}{N}$. According to the concentration inequality of the sample mean for a sub-Gaussian distribution characterized in lemma 1.3 of \cite{rigollet2023high}, with probability at least $1 - \eta$, we have
\begin{equation}\label{equ:concentration mean}
\textstyle |\hat{\mu} - \mu|  \leq \nu\sqrt{\frac{2\log(2/\eta)}{N}}.
\end{equation}
By theorem 6.5 in \cite{Wainwright2019high2},  for any $\delta > 0$ there exist some constants $C_1$, $C_2$ and $C_3$ such that with probability at least $1 - \eta$, it holds that
\[
\textstyle |\hat{\mu}^2 + \hat{\sigma}^2 - \mu^2 - \sigma^2| \leq \nu^2 C_1\big(\frac{1}{N} +\sqrt{\frac{1}{N}}\big) + \nu^2\max\big\{\sqrt{\frac{\log(C_2/\eta)}{C_3N}}, \frac{\log(C_2/\eta)}{C_3N}\big\}.
\]
Since $\frac{1}{N}\leq \frac{1}{\sqrt{N}}$ and $\sqrt{x} \leq 1 + x$ for any $x > 0$, with probability at least $1 - \eta$, we further have
\begin{equation}\label{equ:concentration second-order moment}
\textstyle |\hat{\mu}^2 + \hat{\sigma}^2 - \mu^2 - \sigma^2| \leq \frac{\nu^2}{\sqrt{N}}\big(2C_1 + 1 + \frac{\log(C_2/\eta)}{C_3}\big).
\end{equation}
Note that
\begin{equation}\label{EC:variance decomposition}
|\hat{\sigma}^2 - \sigma^2| \leq |\hat{\mu}^2 + \hat{\sigma}^2 - \mu^2 - \sigma^2| + |\hat{\mu}^2 - \mu^2| \leq |\hat{\mu}^2 + \hat{\sigma}^2 - \mu^2 - \sigma^2| + (\hat{\mu} - \mu)^2 + 2\mu|\hat{\mu} - \mu|,
\end{equation}
where the first inequality follows from the triangle inequality, and the second inequality follows from the fact that $|\hat{\mu}^2 - \mu^2| = |(\hat{\mu} - \mu)^2 + 2\mu(\mu - \hat{\mu})| \leq (\hat{\mu} - \mu)^2 + 2\mu|\hat{\mu} - \mu|$.  Hence, it holds that
\[
|\hat{\sigma} - \sigma| = \frac{|\hat{\sigma}^2 - \sigma^2|}{\hat{\sigma} + \sigma} \leq \frac{|\hat{\sigma}^2 - \sigma^2|}{\sigma} \leq \frac{|\hat{\mu}^2 + \hat{\sigma}^2 - \mu^2 - \sigma^2| + (\hat{\mu} - \mu)^2 + 2\mu|\hat{\mu} - \mu|}{\sigma}.
\]
Applying the reverse union bound to inequalities~\eqref{equ:concentration mean} and~\eqref{equ:concentration second-order moment}, we then have with probability at least $1 - \eta$,
\begin{equation}\label{equ:concentration varaince}
\begin{array}{ll}
|\hat{\sigma} - \sigma| & \leq \frac{1}{\sigma}\big(\frac{\nu^2}{\sqrt{N}}\big(2C_1 + 1 + \frac{\log(2C_2/\eta)}{C_3}\big) + 2\nu^2\frac{\log(4/\eta)}{N} + 2\mu\nu\sqrt{\frac{2\log(4/\eta)}{N}}\big)  \leq \frac{\xi_1 + \xi_2 \log(1/\eta)}{\sqrt{N}},
\end{array}
\end{equation}
where $\xi_1 = \frac{\nu}{\sigma}\big(\nu\big(2C_1 + 1 + \frac{\log(2C_2)}{C_3}\big) + 2(\nu + 2\mu)\log(4) + 2\mu\big)$ and $\xi_2 = \frac{\nu}{\sigma}\big(\frac{\nu}{C_3} + 2\nu + 4\mu\big)$. Here, the second inequality follows from the fact that $\frac{1}{N} \leq \frac{1}{\sqrt{N}}$ and $\sqrt{x} \leq 1 + x$ for any $x \geq 0$.

{\it Step 2.} To derive the concentration property for the mean-variance ambiguity, we divide the argument into two cases based on the expression of $d(D, \mathcal{A}_N)$.

{\it Scenario 1.} If $\frac{\hat{\mu}}{\hat{\sigma}} \geq \frac{\mu}{\sigma}$, according to \Cref{lemma:distance}, we have $d(D,\mathcal{A}_N) = (\hat{\mu} - \mu)^2 + (\hat{\sigma} - \sigma)^2$. By the concentration inequalities for sample mean and sample variance derived in equations~\eqref{equ:concentration mean} and~\eqref{equ:concentration varaince}, it is then immediate to see that with probability at least $1 - \eta$, it holds that
\begin{equation}\label{equ:mean-variance bound}
\begin{array}{rl}
(\hat{\mu} - \mu)^2 + (\hat{\sigma} - \sigma)^2 \leq & \nu^2\frac{4\log(2/\eta)}{N} + \frac{(\xi_1 + \xi_2 \log(1/\eta))^2}{N}
\leq \frac{\xi_1^2 + 2\log(4)\nu^2 + (2\xi_1\xi_2 + 2\nu^2)\log(1/\eta) + \xi_2^2(\log(1/\eta))^2}{\sqrt{N}}.
\end{array}
\end{equation}
Let $\zeta_1 = \xi_1^2 + 2\log(4)\nu^2$ and $\zeta_2 = 2\xi_1\xi_2 + 2\nu^2$. When $\zeta_2^2 \geq 4\zeta_1\xi_2^2$, we have $\zeta_1 + \zeta_2\log(1/\eta) + \xi_2^2(\log(1/\eta))^2 \leq (\frac{\zeta_1}{\xi_2} + \xi_2 \log(1/\eta))^2$. Setting $c_1 = \frac{\xi_1\xi_2 + \nu^2}{\xi_2}$ and $c_2 = \xi_2$ then yields the result. When $\zeta_2^2 < 4\zeta_1\xi_2^2$, we have $\zeta_1 + \zeta_2\log(1/\eta) + \xi_2^2(\log(1/\eta))^2 \leq (\sqrt{\zeta_1} + \xi_2 \log(1/\eta))^2$. Setting $c_1 = \sqrt{\zeta_1}$ and $c_2 = \xi_2$ then yields the desired result.

{\it Scenario 2.} If $\frac{\hat{\mu}}{\hat{\sigma}} < \frac{\mu}{\sigma}$, it holds that
\[
\begin{array}{ll}
\frac{\mu^2\hat{\sigma}^2 - \hat{\mu}^2\sigma^2}{\sigma\hat{\sigma}} & \leq \frac{\mu|\hat{\sigma}^2 - \sigma^2| + \sigma^2|\hat{\mu} - \mu|}{\sigma\hat{\sigma}} \leq \frac{\mu(|\hat{\mu}^2 + \hat{\sigma}^2 - \mu^2 - \sigma^2| + (\hat{\mu} - \mu)^2 + 2\mu|\hat{\mu} - \mu|) + \sigma^2|\hat{\mu} - \mu|}{\sigma\hat{\sigma}},
\end{array}
\]
where the second inequality follows from inequality~\eqref{EC:variance decomposition}. According to inequalities~\eqref{equ:concentration mean} and~\eqref{equ:concentration second-order moment}, with probability at least $1 - \eta$, it holds that
\[
\begin{array}{rl}
& \mu\big(|\hat{\mu}^2 + \hat{\sigma}^2 - \mu^2 - \sigma^2| + (\hat{\mu} - \mu)^2 + 2\mu|\hat{\mu} - \mu|\big) + \sigma^2|\hat{\mu} - \mu| \\ [1mm]
\leq & \frac{\nu^2\mu}{\sqrt{N}}\big(2C_1 + 1 + \frac{\log(2C_2/\eta)}{C_3}\big) + (2\mu^2 + \sigma^2)\nu\sqrt{\frac{2\log(4/\eta)}{N}} + \mu \nu^2\frac{2\log(4/\eta)}{N}\\ [2mm]
\leq & \frac{\delta_1 + \delta_2 \log(1/\eta)}{\sqrt{N}},
\end{array}
\]
where $\delta_1 = \nu^2\mu\big(2C_1 + 1 + \frac{\log(2C_2)}{C_3}\big) + 2\log(4)(2\mu^2 + \sigma^2 + \nu\mu)\nu + (2\mu^2 + \sigma^2)\nu$ and $\delta_2 =  \frac{\nu^2\mu}{C_3} + 2(2\mu^2 + \sigma^2)\nu + 2\mu\nu^2$. When $N$ is sufficiently large, $\sigma > \frac{\xi_1 + \xi_2 \log(1/\eta)}{\sqrt{N}}$, and then \eqref{equ:concentration varaince} implies that $\hat{\sigma} \leq \sigma - \frac{\xi_1 + \xi_2 \log(1/\eta)}{\sqrt{N}}$ with probability at least $1-\eta$. Hence, we have
\[
\textstyle \frac{\mu^2\hat{\sigma}^2 - \hat{\mu}^2\sigma^2}{\sigma\hat{\sigma}} \leq \frac{\frac{\delta_1 + \delta_2 \log(1/\eta)}{\sqrt{N}}}{\sigma(\sigma - \frac{\xi_1 + \xi_2 \log(1/\eta)}{\sqrt{N}})} = \frac{\delta_1 + \delta_2 \log(1/\eta)}{\sqrt{N}\sigma^2 - \sigma(\xi_1 + \xi_2 \log(1/\eta))},
\]
holds with probability at least $1-\eta$. Note that for sufficiently large $N$, we have $\sqrt{N}\sigma^2 > \delta_1 + \delta_2 \log(1/\eta) + \sigma(\xi_1 + \xi_2 \log(1/\eta))$, which implies that $\frac{\delta_1 + \delta_2 \log(1/\eta)}{\sqrt{N}\sigma^2 - \sigma(\xi_1 + \xi_2 \log(1/\eta))} < \frac{\delta_1 + \delta_2 \log(1/\eta) + \sigma(\xi_1 + \xi_2 \log(1/\eta))}{\sqrt{N}\sigma^2}$.
This inequality, together with~\eqref{equ:mean-variance bound}, implies that with probability at least $1 - \eta$, for sufficiently large $N$, it holds that
\[
\begin{array}{ll}
d(D,\mathcal{A}_N) & \leq (\hat{\mu} - \mu)^2 + (\hat{\sigma} - \sigma)^2 + \frac{\mu^2\hat{\sigma}^2 - \hat{\mu}^2\sigma^2}{\sigma\hat{\sigma}} \leq \frac{(\xi_1^2 + 2\log(2)\nu^2)\sigma^2 + \delta_1 + \sigma\xi_1 + (2\xi_1\xi_2\sigma^2 + 2\nu^2\sigma^2 + \delta_2 + \sigma\xi_2)\log(1/\eta) + \xi_2^2\sigma^2(\log(1/\eta))^2}{\sqrt{N}\sigma^2}.
\end{array}
\]
Using a similar argument as in {\it Scenario 1}, we can also show that with probability at least $1 - \eta$, for sufficiently large $N$, $d(D,\mathcal{A}_N) \leq \frac{(c_1 + c_2 \log(1/\eta))^2}{\sqrt{N}}$ where $c_1$ and $c_2$ are some constants that only depend on $\mu$, $\sigma$ and $\nu$. 
\hfill $\square$
\vspace{2mm}

\noindent{\bf Proof of \Cref{thm:guarantee}.} We proceed in two steps. In the first step, leveraging the equivalence between
\begin{equation}\label{EC:constrained}
\Upsilon^\star_\varepsilon = \max_{q \geq 0} \; \min_{d(F,\mathcal{A}) \leq \varepsilon} \; \mathbb{E}_F[\pi(q,\tilde{u})]
\end{equation}
and \ref{prob:misspecification} (as shown in~\Cref{lemma:misspecification constrained}), we characterize the relationship between $\varepsilon$ and $\alpha$. In the second step, we translate the finite-sample performance guarantee of~\eqref{EC:constrained} characterized in \Cref{prop:concentration} as the performance guarantee of \ref{prob:misspecification}.

{\it Step 1.} $\;$ Given $N$ and $\varepsilon_N + d_{\rm OT}(F,D)$, consider the constrained problem~\eqref{EC:constrained} with $\mathcal{A} = \mathcal{A}_N$ and $\varepsilon = \varepsilon_N + d_{\rm OT}(F,D)$. Suppose that $\alpha_N$ is an index of misspecification aversion such that \ref{prob:misspecification} with $\mathcal{A} = \mathcal{A}_N$ and $\alpha = \alpha_N$ has the same optimal solution as that of the corresponding constrained problem \eqref{EC:constrained}. Denote by $\Pi_{\alpha_N}^\star$ the optimal value of \ref{prob:misspecification} with $\alpha = \alpha_N$. By \Cref{lemma:misspecification constrained}, we have
\begin{equation*}
\Upsilon^\star_{\varepsilon_N + d_{\rm OT}(F,D)} = \max_{\alpha \geq 0} \; \big\{\Pi^\star_\alpha - (\varepsilon_N + d_{\rm OT}(F,D))\alpha\big\} = \Pi^\star_{\alpha_N} - (\varepsilon_N + d_{\rm OT}(F,D)) \alpha_N.
\end{equation*}
In the following, we characterize the expression of $\alpha_N$. Plugging the expression of the optimal order quantity $q_\alpha^\star$ into the worst-case transformed expectation characterized in \Cref{prop:transformed expectation}, we have
\[
\Pi_\alpha^\star - (\varepsilon_N + d_{\rm OT}(F,D))\alpha = \left\{
\begin{array}{ll}
\kappa\alpha \hat{v}^{\star2} - (\varepsilon_N + d_{\rm OT}(F,D))\alpha &~~ \alpha < \frac{p}{2\hat{v}^\star}\\
(p-c)\hat{v}^\star - \frac{p(p-c)}{4\alpha} - (\varepsilon_N + d_{\rm OT}(F,D))\alpha &~~ \alpha \geq \frac{p}{2\hat{v}^\star},
\end{array}
\right.
\]
where $\hat{v}^\star = \hat{\mu}-\hat{\sigma}\sqrt{\frac{1-\kappa}{\kappa}}$.
Note that $\Pi_\alpha^\star - (\varepsilon_N + d_{\rm OT}(F,D))\alpha$ is concave in $\alpha$. If $\varepsilon_N + d_{\rm OT}(F,D) \geq \kappa \hat{v}^{\star2}$, then $\Pi_\alpha^\star - (\varepsilon_N + d_{\rm OT}(F,D))\alpha$ is decreasing in $\alpha$, and hence $\max_{\alpha \geq 0} \{\Pi^\star_\alpha - (\varepsilon_N + d_{\rm OT}(F,D)))\alpha\} = 0$ with $\alpha_N = 0$. If $\varepsilon_N + d_{\rm OT}(F,D) < \kappa \hat{v}^{\star2}$, then by the first-order optimality condition, the maximum is attained at $\alpha_N = \frac{\sqrt{p(p-c)}}{2\sqrt{\varepsilon_N + d_{\rm OT}(F,D)}}$
and $\max_{\alpha \geq 0} \{\Pi^\star_\alpha - (\varepsilon_N + d_{\rm OT}(F,D))\alpha\} = \Pi_{\alpha_N}^\star - \frac{1}{2}\sqrt{p(p-c)(\varepsilon_N + d_{\rm OT}(F,D))}$.
To summarize the second step, we have obtained
\[
\alpha_N = \left\{
\begin{array}{l@{\quad}l}
\frac{\sqrt{p(p-c)}}{2\sqrt{\varepsilon_N + d_{\rm OT}(F,D)}} &~\varepsilon_N + d_{\rm OT}(F,D) < \kappa \hat{v}^{\star2}  \\
0 & ~\varepsilon_N + d_{\rm OT}(F,D) \geq \kappa \hat{v}^{\star2}
\end{array}
\right.
~~{\rm and}~~
\textstyle \Upsilon_{\varepsilon_N + d_{\rm OT}(F,D)}^\star = \Big(\Pi_{\alpha_N}^\star - \frac{1}{2}\sqrt{p(p-c)(\varepsilon_N + d_{\rm OT}(F,D))}\Big)^+.
\]

{\it Step 2.} $\;$ Denote by $q_{\alpha_N}^\star$ the optimal solution to \ref{prob:misspecification} with $\alpha = \alpha_N$, which, by \Cref{lemma:misspecification constrained}, is also the optimal solution to the constrained problem~\eqref{EC:constrained}. For any $F$ such that $d(F,\mathcal{A}_N) \leq \varepsilon_N + d_{\rm OT}(F,D)$, we have $\mathbb{E}_{F}[\pi(q^\star_{\alpha_N},\tilde{v})] \geq \Upsilon_{\varepsilon_N + d_{\rm OT}(F,D)}^\star$. Therefore,
\[
\begin{array}{ll}
\mathbb{P}_{D^N}\big\{\mathbb{E}_{F}[\pi(q^\star_{\alpha_N},\tilde{v})] \geq \Upsilon_{\varepsilon_N + d_{\rm OT}(F,D)}^\star\big\} & \displaystyle \geq \mathbb{P}_{D^N}\{d(F,\mathcal{A}_N) \leq \varepsilon_N + d_{\rm OT}(F,D)\} \geq \mathbb{P}_{D^N}\{d(D,\mathcal{A}_N) \leq \varepsilon_N\},
\end{array}
\]
where the second inequality follows since $d(D,\mathcal{A}_N) \leq \varepsilon_N$ implies $d(F,\mathcal{A}_N) \leq d(D,\mathcal{A}_N) + d_{\rm OT}(F,D) \leq \varepsilon_N + d_{\rm OT}(F,D)$.
This, together with the inequality in \Cref{prop:concentration} implies that $\mathbb{P}_{D^N}\{\mathbb{E}_{F}[\pi(q^\star_{\alpha_N},\tilde{u})] \geq \Upsilon_{\varepsilon_N + d_{\rm OT}(F,D)}^\star\} \geq 1- \eta$.
Plugging the expressions of $\alpha_N$ and $\Upsilon^\star_{\varepsilon_N + d_{\rm OT}(F,D)}$ given in the first step, we obtain the desired result.   
\hfill $\square$
\vspace{2mm}

\noindent{\bf Proof of \Cref{thm:multi product dual}.} $\;$
For ease of notation, define $\mathcal{G} = \{G\in\mathcal{P}_M\mid\mathbb{E}_G[\tilde{v}_i] = \mu_i ~~\forall i \in [M]\}$ and $\mathcal{C}_i = \{G\in\mathcal{P}\mid\mathbb{E}_G[\tilde{v}_i] = \mu_i\}$ for $i \in [M]$. Introducing the dual variable $\lambda \geq 0$ to the sum-of-variance constraint in the ambiguity set, then the \ref{prob:multi product} model can be equivalently reformulated as
\[
\textstyle \max\limits_{\lambda \geq 0} \; \bigg\{-\lambda K + \max\limits_{\bm{q} \geq 0} \; \min\limits_{F\in\mathcal{P}_M,\; G\in\mathcal{G}} \; \mathbb{E}_F\Big[\sum\limits_{i\in[M]}\pi(q_i,\tilde{u}_i)\Big] + \lambda\cdot \mathbb{E}_{G}\Big[\sum\limits_{i\in[M]}\tilde{v}_i^2\Big] + \alpha \cdot d_{\rm OT}(F,G)\Big\},
\]
which, by noting that $\mathcal{G}$ is decomposable with respect to multiple products, further reduces to
\begin{equation}\label{EC:dual multi product}
\textstyle \max\limits_{\lambda \geq 0} \; \Big\{-\lambda K + \sum\limits_{i\in[M]} \max\limits_{q_i \geq 0} \; \min\limits_{F_i \in \mathcal{P},\; G_i \in\mathcal{C}_i} \; \{ \mathbb{E}_{F_i}[\pi(q_i,\tilde{u}_i)] + \lambda \cdot \mathbb{E}_{G_i}[\tilde{v}_i^2] + \alpha \cdot d_{\rm OT}(F_i, G_i)\} \Big\}.
\end{equation}
In the remainder of the proof, we solve for the optimal $\lambda^\star$ and $\bm{q}_\alpha^\star$ of problem \eqref{EC:dual multi product}. Invoking the interchangeability principle characterized in \Cref{lemma:interchangeability}, then problem \eqref{EC:dual multi product} becomes
\[
\textstyle \max\limits_{\lambda \geq 0} \; \Big\{-\lambda K + \sum\limits_{i\in[M]}\max\limits_{q_i \geq 0} \; \min\limits_{G_i\in\mathcal{C}_i} \; \mathbb{E}_{G_i}\Big[\min\limits_{u_i \geq 0}\{\pi(q_i,u_i)+ \lambda \cdot \tilde{v}_i^2 + \alpha \cdot (u_i-\tilde{v}_i)^2\}\Big]\Big\},
\]
which, by defining $\digamma_i(\lambda, q_i, v_i) = \lambda \cdot v_i^2 + \min_{u_i \geq 0}\{\pi(q_i,u_i)+ \alpha \cdot (u_i-v_i)^2\}$ (see \Cref{lemma:newsvendor type-2} for its closed-form expression) for each $i\in[M]$ and $u_i \geq 0$, can be equivalently written as
\begin{equation}\label{EC:dual multi product 2}
\textstyle \max\limits_{\lambda \geq 0} \; \Big\{-\lambda K + \sum\limits_{i\in[M]}\max\limits_{q_i \geq 0} \; \min\limits_{G_i\in\mathcal{C}_i} \; \mathbb{E}_{G_i} [\digamma_i(\lambda, q_i, \tilde{v}_i)]\Big\}.
\end{equation}
In the following, our remaining proof proceeds in three steps: deriving the expression for $L_i(q_i) = \min_{G_i\in\mathcal{C}_i} \; \mathbb{E}_{G_i} [\digamma_i(\lambda, q_i, \tilde{v}_i)]$ (\textit{Step} 1), optimizing over $q_i$ to solve $\max_{q_i \geq 0} L_i(q_i)$ for each $i \in [M]$ (\textit{Step} 2), and finally, optimizing over $\lambda \geq 0$ (\textit{Step} 3). Note that given $\lambda \geq 0$, we solve the inner maximization of problem \eqref{EC:dual multi product 2} over $q_i \geq 0$ for each $i \in [M]$ in \textit{Step} 1 and \textit{Step} 2. 

{\it Step 1.} $\;$ We drop the subscript `$i$' to avoid clutter. Given $q \geq 0$, $L(q)$ is a classical moment problem:
\begin{equation}\label{prob:primal multi product}\tag{\sc Primal}
\begin{array}{cll}
\displaystyle \min_G & \int_{\mathbb{R}_+} \; \digamma(\lambda, q, v) \; {\rm d}G(v) \\[1mm]
{\rm s.t.} & \int_{\mathbb{R}_+} \; v \; {\rm d}G(v) = \mu &~~\cdots~s_\alpha \\[1mm]
& \int_{\mathbb{R}_+}\;{\rm d}G(v) = 1 &~~\cdots~t_\alpha  \\[1mm]
& G \in \mathcal{M}_+,
\end{array}
\end{equation}
whose dual is given by
\begin{equation}\label{prob:dual multi product}\tag{\sc Dual}
\begin{array}{cll}
\displaystyle \max_{s_\alpha,\;t_\alpha} & \mu s_\alpha - t_\alpha \\[1mm]
{\rm s.t.} & v s_\alpha - t_\alpha \leq \digamma(\lambda, q, v) &~\forall v \geq 0.
\end{array}
\end{equation}
We next derive the expression of $L(q)$ by constructing a pair of primal and dual feasible solutions that attain the same objective value (that is, strong duality holds between \ref{prob:primal multi product} and \ref{prob:dual multi product}). The argument breaks into nine scenarios based on the value of $q$.

{\it Scenario 1.1}. $\;$ When $q \leq \frac{p}{4\alpha}$ and $\frac{\alpha \lambda \mu^2}{p(\lambda + \alpha)} \leq q \leq \frac{\alpha (\lambda + \alpha) \mu^2}{p\lambda}$, we first construct a feasible distribution to \ref{prob:primal multi product} as follows:
\begin{equation}\label{equ:primal solution multi product}
\textstyle G_\alpha = \textstyle \bigg(\frac{\sqrt{\frac{\lambda+\alpha}{\lambda}} - \mu\sqrt{\frac{\alpha}{pq}}}{\sqrt{\frac{\lambda+\alpha}{\lambda}} - \sqrt{\frac{\lambda}{\lambda+\alpha}}}\bigg)\cdot \delta_{\sqrt{\frac{\lambda pq}{\alpha(\lambda+\alpha)}}} + \bigg(\frac{\mu\sqrt{\frac{\alpha}{pq}} - \sqrt{\frac{\lambda}{\lambda+\alpha}}}{\sqrt{\frac{\lambda+\alpha}{\lambda}} - \sqrt{\frac{\lambda}{\lambda+\alpha}}}\bigg)\cdot \delta_{\sqrt{\frac{(\lambda+\alpha) pq}{\alpha \lambda}}}.
\end{equation}
One can verify that $G_\alpha \in \mathcal{A}$ and the corresponding primal objective value is equal to $\mathbb{E}_{G_\alpha}[\digamma(\lambda,q,\tilde{v})] = 2\mu\sqrt{\lambda(\lambda+\alpha)pq/\alpha} -\frac{\lambda pq}{\alpha} - cq$.
We next construct a dual feasible that attains the same dual objective value. 
Consider the following solution
\begin{equation}\label{equ:dual solution multi product}
\textstyle s_\alpha = 2\sqrt{\frac{\lambda(\lambda+\alpha)pq}{\alpha}},~~t_\alpha =  \frac{\lambda pq}{\alpha} + cq,
\end{equation}
which satisfies $s_\alpha \mu - t_\alpha = 2\mu\sqrt{\lambda(\lambda+\alpha)pq/\alpha} -\frac{\lambda pq}{\alpha} - cq$.
Similarly to the proof of Proposition~\ref{prop:transformed expectation}, one can verify that solution~\eqref{equ:dual solution multi product} is feasible to \ref{prob:dual multi product} and this establishes the strong duality.

{\it Scenario 1.2}. $\;$ When $q \leq \frac{p}{4\alpha}$ and $q \leq \frac{\alpha\lambda\mu^2}{p(\lambda+\alpha)}$, we construct a primal feasible solution $G_\alpha = \delta_\mu$ with a primal objective value $\mathbb{E}_{G_\alpha}[\digamma(\lambda, q, \tilde{v})] = \lambda \mu^2 + pq - cq$. Consider the solution $s_\alpha = 2\lambda \mu$ and $t_\alpha = \lambda \mu^2 - (p-c)q$, which is feasible and satisfies $s_\alpha \mu - t_\alpha = \lambda \mu^2 + pq - cq$.

{\it Scenario 1.3}. $\;$ When $q \leq \frac{p}{4\alpha}$ and $q \geq \frac{\alpha(\lambda+\alpha)\mu^2}{p\lambda}$, we construct a primal feasible solution $G_\alpha = \delta_\mu$ with a primal objective value $\mathbb{E}_{G_\alpha}[\digamma(\lambda, q, \tilde{v})] = (\lambda+\alpha)\mu^2 - cq$. Consider the solution $s_\alpha = 2(\lambda+\alpha) \mu$ and $t_\alpha = (\lambda+\alpha)\mu^2 + cq$, which is feasible and satisfies $s_\alpha \mu - t_\alpha = (\lambda+\alpha)\mu^2 - cq$.

{\it Scenario 1.4}. $\;$ When $\frac{p}{4\alpha} \leq q \leq \frac{p}{4\alpha} + \frac{p}{4\lambda}$ and $\frac{\alpha \lambda \mu^2}{p(\lambda + \alpha)} \leq q \leq \frac{\alpha (\lambda + \alpha) \mu^2}{p\lambda}$, we consider the pair of primal feasible solution~\eqref{equ:primal solution multi product} and dual solution~\eqref{equ:dual solution multi product}. Upon the results established in \textit{Scenario}~1.1, we can verify that solution~\eqref{equ:dual solution 1} is feasible to \ref{prob:dual}. 

{\it Scenario 1.5}. $\;$ When $\frac{p}{4\alpha} \leq q \leq \frac{p}{4\alpha} + \frac{p}{4\lambda}$ and $q \leq \frac{\alpha\lambda\mu^2}{p(\lambda+\alpha)}$, we construct a primal feasible solution $G_\alpha = \delta_\mu$ with a primal objective value $\mathbb{E}_{G_\alpha}[\digamma(\lambda, q, \tilde{v})] = \lambda \mu^2 + pq - cq$. Consider the solution $s_\alpha = 2\lambda \mu$ and $t_\alpha = \lambda \mu^2 - (p-c)q$, which is feasible and satisfies $s_\alpha \mu - t_\alpha = \lambda \mu^2 + pq - cq$.

{\it Scenario 1.6}. $\;$ When $\frac{p}{4\alpha} \leq q \leq \frac{p}{4\alpha} + \frac{p}{4\lambda}$ and $q \geq \frac{\alpha(\lambda+\alpha)\mu^2}{p\lambda}$, we construct a primal feasible solution $G_\alpha = \delta_\mu$ with a primal objective value $\mathbb{E}_{G_\alpha}[\digamma(\lambda, q, \tilde{v})] = (\lambda+\alpha)\mu^2 - cq$. Consider the solution $s_\alpha = 2(\lambda+\alpha) \mu$ and $t_\alpha = (\lambda+\alpha)\mu^2 + cq$, which is feasible and satisfies $s_\alpha \mu - t_\alpha = (\lambda+\alpha)\mu^2 - cq$.

{\it Scenario 1.7}. $\;$ When $q \geq \frac{p}{4\alpha} + \frac{p}{4\lambda}$ and $\mu - \frac{p}{4\alpha} - \frac{p}{4\lambda} \leq q \leq \mu - \frac{p}{4\alpha} + \frac{p}{4\lambda}$, we first construct a primal feasible solution: $G_\alpha = \big(\frac{1}{2} - \frac{2\lambda}{p}(\mu - q - \frac{p}{4\alpha})\big)\cdot\delta_{q+\frac{p}{4\alpha}-\frac{p}{4\lambda}} + \big(\frac{1}{2} + \frac{2\lambda}{p}(\mu - q - \frac{p}{4\alpha})\big) \cdot \delta_{q+\frac{p}{4\alpha}+\frac{p}{4\lambda}}$.
One can verify that $G_\alpha \in \mathcal{A}$ and the corresponding objective value under $G_\alpha$ is equal to $\mathbb{E}_{G_\alpha}[\digamma(\lambda,q,\tilde{v})] = \big(2\lambda q + \frac{(\lambda + \alpha)p}{2\alpha}\big)\mu - \lambda q^2 - \frac{(\lambda + \alpha)pq}{2\alpha} - \frac{(\lambda + \alpha)^2 p^2}{16\alpha^2\lambda} + (p-c)q$.
Moreover, we consider the solution: $s_\alpha = 2\lambda q + \frac{(\lambda + \alpha)p}{2\alpha},\; t_\alpha = \lambda q^2 + \frac{(\lambda + \alpha)pq}{2\alpha} + \frac{(\lambda + \alpha)^2 p^2}{16\alpha^2\lambda} - (p-c)q$, which is dual feasible.

{\it Scenario 1.8}. $\;$ When $q \geq \frac{p}{4\alpha}$ and $q \leq \mu - \frac{p}{4\alpha} - \frac{p}{4\lambda}$, we construct a primal feasible solution $G_\alpha = \delta_\mu$ with a primal objective value $\mathbb{E}_{G_\alpha}[\digamma(\lambda, q, \tilde{v})] = \lambda \mu^2 + pq - cq$. Consider the solution $s_\alpha = 2\lambda \mu$ and $t_\alpha = \lambda \mu^2 - (p-c)q$, which satisfies $s_\alpha \mu - t_\alpha = \lambda \mu^2 + pq - cq$. The feasibility of $(s_\alpha, t_\alpha)$ can be verified easily.

{\it Scenario 1.9}. $\;$ When $q \geq \frac{p}{4\alpha}$ and $q \geq \mu - \frac{p}{4\alpha} + \frac{p}{4\lambda}$, we first construct a primal feasible solution $G_\alpha = \delta_\mu$ with a primal objective value $\mathbb{E}_{G_\alpha}[\digamma(\lambda, q, \tilde{v})] = \lambda\mu^2 + p(\mu - \frac{p}{4\alpha}) - cq$. Consider the solution $s_\alpha = 2\lambda\mu + p$ and $t_\alpha = \lambda\mu^2 + \frac{p^2}{4\alpha} + cq$, which is feasible and satisfies $s_\alpha \mu - t_\alpha = \lambda\mu^2 + p(\mu - \frac{p}{4\alpha}) - cq$.

To summarize \textit{Step~1}, we note that when the constructed primal feasible distribution is $G_\alpha = \delta_\mu$ (\textit{i.e.}, {\it Scenarios 1.2, 1.3, 1.5, 1.6, 1.8, 1.9}), the objective function $L(q)$ is either increasing or decreasing in $q$, implying that the maximum can not be attained in these scenarios. Therefore, we only need to focus on the remaining scenarios (\textit{i.e.}, {\it Scenarios 1.1, 1.4, 1.7}) where
\[
L(q) = \left\{
\begin{array}{ll} 2\mu\sqrt{\frac{\lambda(\lambda+\alpha)pq}{\alpha}} -\frac{\lambda pq}{\alpha} - cq &~~~ q \in \mathcal{Q}_1 \\ [2mm]
\big(2\lambda q + \frac{(\lambda + \alpha)p}{2\alpha}\big)\mu - \lambda q^2 - \frac{(\lambda + \alpha)pq}{2\alpha} - \frac{(\lambda + \alpha)^2 p^2}{16\alpha^2\lambda} + (p-c)q &~~~q \in \mathcal{Q}_2 
\end{array}
\right.
\]
with $\mathcal{Q}_1=\{q\mid q \leq \frac{p}{4\alpha} + \frac{p}{4\lambda},\;\frac{\alpha \lambda \mu^2}{p(\lambda + \alpha)} \leq q \leq \frac{\alpha (\lambda + \alpha) \mu^2}{p\lambda}\}$ and $\mathcal{Q}_2 = \{q\mid q \geq \frac{p}{4\alpha} + \frac{p}{4\lambda},\; \mu - \frac{p}{4\alpha} - \frac{p}{4\lambda} \leq q \leq \mu - \frac{p}{4\alpha} + \frac{p}{4\lambda}\}$.

{\it Step 2.} $\;$ We consider three scenarios based on the values of $\alpha$ and $\lambda$ to solve for the optimal $q$.

{\it Scenario 2.1.} $\;$ Suppose that $\lambda \leq \frac{c}{2\mu}$. For any $q \in \mathcal{Q}_1$, setting the derivative of $L(q)$ to $0$ yields
\begin{equation}\label{EC:quantity 1 multi product}
\textstyle q_\alpha^\star(\lambda) = \frac{\lambda (\lambda + \alpha) p \mu^2}{\alpha(\lambda p/\alpha + c)^2}.
\end{equation}
One can verify that $\frac{\alpha \lambda \mu^2}{p(\lambda + \alpha)} \leq q_\alpha^\star \leq \frac{\alpha (\lambda + \alpha) \mu^2}{p\lambda}$. We next show that $q_\alpha^\star(\lambda) \leq \frac{p}{4\lambda} + \frac{p}{4\alpha}$. Since given $\alpha \geq 0$, $\frac{\partial}{\partial \alpha}\big(q_\alpha^\star(\lambda) - \frac{p}{4\alpha}\big) = \frac{p\lambda^2\mu^2(\alpha (p-c) + \alpha p + \lambda p)}{(\alpha c + \lambda p)^3} + \frac{p}{4\alpha^2} \geq 0$, and hence
we have $q_\alpha^\star(\lambda) - \frac{p}{4\alpha} \leq \lim_{\alpha \to \infty} \big(q_\alpha^\star(\lambda) - \frac{p}{4\alpha}\big) = \frac{\lambda p \mu^2}{c^2} \leq \frac{p}{4\lambda}$,
where the last inequality follows from $\lambda \leq \frac{c}{2\mu}$. Since $L(q)$ is concave, $q_\alpha^\star(\lambda)$ is indeed optimal.

{\it Scenario 2.2.} $\;$ Suppose that $\lambda \geq \frac{c}{2\mu}$ and $\alpha < \frac{p}{2(\mu - c/(2\lambda))}$. For $q \in \mathcal{Q}_1$, setting $L'(q)$ to $0$ yields $\textstyle q_\alpha^\star(\lambda) = \frac{\lambda (\lambda + \alpha) p \mu^2}{\alpha(\lambda p/\alpha + c)^2} \in [\frac{\alpha \lambda \mu^2}{p(\lambda + \alpha)}, \frac{\alpha (\lambda + \alpha) \mu^2}{p\lambda}]$. Since $q_\alpha^\star(\lambda) - \frac{p}{4\alpha}$ is increasing in $\alpha$, we have $q_\alpha^\star(\lambda) - \frac{p}{4\alpha} < \lim_{\alpha \to \frac{p}{2(\mu - c/(2\lambda))}} \big(q_\alpha^\star(\lambda) - \frac{p}{4\alpha}\big) = \frac{p}{4\lambda}$.
This implies that $q_\alpha^\star(\lambda) \in \mathcal{Q}_1$. Hence, $q_\alpha^\star(\lambda)$ is optimal.

{\it Scenario 2.3.} $\;$ Suppose that $\lambda \geq \frac{c}{2\mu}$ and $\alpha \geq \frac{p}{2(\mu - c/(2\lambda))}$. For $q \in \mathcal{Q}_2$, setting $L'(q)$ to $0$ yields
\begin{equation}\label{EC:quantity 2 multi product}
\textstyle q_\alpha^\star(\lambda) = \mu + \frac{p-2c}{4\lambda} - \frac{p}{4\alpha}.
\end{equation}
One can verify that $\mu - \frac{p}{4\alpha} - \frac{p}{4\lambda} \leq q_\alpha^\star(\lambda) \leq \mu - \frac{p}{4\alpha} + \frac{p}{4\lambda}$. Since $\lambda \geq \frac{c}{2\mu}$ and $\alpha \geq \frac{p}{2(\mu - c/(2\lambda))}$, we have $q_\alpha^\star(\lambda) = \mu + \frac{p-2c}{4\lambda} - \frac{p}{2\alpha} + \frac{p}{4\alpha} \geq \frac{p}{4\alpha} + \frac{p}{4\lambda}$,
which indicates that $q_\alpha^\star(\lambda) \in \mathcal{Q}_2$. Therefore, $q_\alpha^\star(\lambda)$ is optimal.

To summarize {\it Step 2}, we note that {\it Scenario~2.1} and {\it Scenario~2.2} correspond to $\alpha < \frac{p}{2(\mu - c/(2\lambda))^+}$ and {\it Scenario~2.3} corresponds to $\alpha \geq \frac{p}{2(\mu - c/(2\lambda))^+}$. Hence, we have
\[
q_\alpha^\star = 
\left\{
\begin{array}{l@{\quad}l}
\mu + \frac{p -2c}{4\lambda} - \frac{p}{4\alpha} & \alpha \geq \frac{p}{2(\mu - c/(2\lambda))^+}~~\big(\lambda \geq \frac{c}{(2\mu - p/\alpha)^+}\big)\\ [2mm]
\frac{\lambda (\lambda + \alpha) p \mu^2}{\alpha(\lambda p/\alpha + c)^2} & \alpha < \frac{p}{2(\mu - c/(2\lambda))^+}~~\big(\lambda < \frac{c}{(2\mu - p/\alpha)^+}\big).
\end{array}
\right.
\]

{\it Step 3.} $\;$ We next solve the optimal $\lambda^\star$ of the outer maximization of problem~\eqref{EC:dual multi product 2}, whose objective function we denote by
\[
\textstyle Q(\lambda) = -\lambda K + \sum\limits_{i\in[M]}\max\limits_{q_i \geq 0} \; \min\limits_{G_i\in\mathcal{C}_i} \; \mathbb{E}_{G_i} [\digamma_i(\lambda, q_i, \tilde{v}_i)] = -\lambda K + \sum\limits_{i\in[M]}L_i(q_{i,\alpha}^\star(\lambda)).
\]
Note that given $v_i \geq 0$ and $i \in [M]$, $\digamma_i(\lambda, q_i, v_i)$ is jointly concave in $(\lambda, q_i)$. Since concavity is preserved under expectation and maximization, 
$Q(\lambda)$ is concave in $\lambda$. When $\lambda < \bar{\lambda}_i = \frac{c_i}{(2\mu_i - p_i/\alpha)^+}$ ({\it Scenario~2.1} and {\it Scenario~2.2}), plugging the expression of $q_{i,\alpha}^\star(\lambda)$ in~\eqref{EC:quantity 1 multi product}, we have $L_i(q_{i,\alpha}^\star(\lambda)) = \frac{\lambda p_i \mu_i^2\alpha(\lambda+\alpha)}{\lambda p_i + \alpha c_i}$ and $\frac{\partial L_i(q_{i,\alpha}^\star(\lambda))}{\partial \lambda} = \frac{\mu_i^2 p_i (\alpha^2 c_i + 2 \alpha c_i \lambda + \lambda^2 p_i)}{(\lambda p_i + \alpha c_i)^2}$.
When $\lambda \geq \bar{\lambda}_i = \frac{c_i}{(2\mu_i - p_i/\alpha)^+}$ ({\it Scenario 2.3}), plugging the expression of $q_{i,\alpha}^\star(\lambda)$ in~\eqref{EC:quantity 2 multi product}, we have $L_i(q_{i,\alpha}^\star(\lambda)) = \frac{\lambda (c_i - p_i) p_i + \alpha (c^2 + 4 \lambda \mu_i (\mu_i + p_i) - c (4 \lambda \mu_i + p_i))}{4 \alpha \lambda}$ and $\frac{\partial L_i(q_{i,\alpha}^\star(\lambda))}{\partial \lambda} = \frac{c_i (p_i - c_i)}{4 \lambda^2}$.
For any $\lambda \in (\bar{\lambda}_{j-1},\bar{\lambda}_{j})$ with $j\in[M+1]$, we then have
\[
\begin{array}{ll}
Q(\lambda) & = -\lambda K + \sum\limits_{i \in [M]\backslash [j - 1]}\frac{\lambda p_i \mu_i^2\alpha(\lambda+\alpha)}{\lambda p_i + \alpha c_i} + \sum\limits_{i \in [j - 1]}\frac{4\alpha\lambda \mu_i (\mu_i + p_i - c_i) - (\lambda p_i + \alpha c_i)(p_i - c_i)}{4 \alpha \lambda} \\ [3mm]
Q'(\lambda) & = -K + \sum\limits_{i \in [M]\backslash [j - 1]}\frac{\mu_i^2 p_i (\alpha^2 c_i + 2 \alpha c_i \lambda + \lambda^2 p_i)}{(\lambda p_i + \alpha c_i)^2} + \sum\limits_{i \in [j - 1]}\frac{c_i (p_i - c_i)}{4 \lambda^2} = -K + \Theta_j(\lambda).
\end{array}
\]
It can be noted that $Q(\lambda)$ is concave and $Q'(\lambda)$ is always decreasing in $(\bar{\lambda}_{j-1},\bar{\lambda}_{j})$, where at the end points it holds that $Q'_-(\bar{\lambda}_{j}) \geq Q'_+(\bar{\lambda}_{j})$ for any $j\in[M+1]$. On the one hand, if $Q'_+(\bar{\lambda}_{i^\star - 1}) = -K + \Theta_{i^\star}(\bar{\lambda}_{i^\star - 1}) \leq 0$,
then when $i^\star = 1$, we have $Q'(\lambda) \leq Q'_+(\bar{\lambda}_0) \leq 0$ for $\lambda \geq 0$, and hence $\lambda^\star = \bar{\lambda}_0 = 0$. When $i^\star > 1$, by the definition of $i^\star$, we have $Q'_-(\bar{\lambda}_{i^\star - 1}) \geq 0$. By the optimality condition of concave functions, it is then clear that $\lambda^\star = \bar{\lambda}_{i^\star - 1}$. On the other hand, if $Q'_+(\bar{\lambda}_{i^\star - 1}) = -K + \Theta_{i^\star}(\bar{\lambda}_{i^\star - 1}) > 0$,
then since (\textit{i}) $Q'_-(\bar{\lambda}_{i^\star}) < 0$ by the definition of $i^\star$ and (\textit{ii}) $Q'(\lambda)$ is continuous in $(\bar{\lambda}_{i^\star - 1}, \bar{\lambda}_{i^\star})$, there must exist some $\lambda^\star \in (\bar{\lambda}_{i^\star-1},\bar{\lambda}_{i^\star})$ such that $Q'(\lambda^\star) = 0$ ({\it i.e.}, $\Theta_{i^\star}(\lambda^\star) = K$), concluding \textit{Step~3}. 
\hfill $\square$
\vspace{2mm}

\noindent{\bf Proof of \Cref{thm:multi product}.} $\;$
With the optimal dual variable $\lambda^\star$, we can determine the optimal order quantity leveraging \eqref{EC:quantity 1 multi product} and \eqref{EC:quantity 2 multi product} established in {\it Step~2} of the proof for \Cref{thm:multi product dual}. \hfill $\square$
\vspace{2mm}

\noindent{\bf Proof of \Cref{thm:double misspecification}.} $\;$
By the interchangeability principle (see \Cref{lemma:interchangeability}), we have
\[
\max_{\psi \geq 0} \; \min_{G\in\mathcal{B}(\theta)} \; \min_{F\in\mathcal{P}} \; \{\mathbb{E}_{F}[\pi(\psi,\tilde{u})] + \alpha\cdot d_{\rm OT}(F,G)\} = \max_{\psi \geq 0} \; \min_{G\in\mathcal{B}(\theta)} \;  \mathbb{E}_{G}\bigg[\min_{u\geq0}\{\pi(\psi,u) + \alpha(u-\tilde{v})^2\}\bigg],
\]
Using a standard duality argument (see, \textit{e.g.}, \citealp{gao2023distributionally}), the right-hand side problem becomes
\[
\begin{array}{rl}
&\max\limits_{\psi \geq 0}\;\sup\limits_{t \geq 0} \; \big\{-t\theta + \mathbb{E}_{H}\big[\min\limits_{u\geq 0,\; v\geq 0}\;\{\pi(\psi,u) + \alpha(u-v)^2 + t(v-\tilde{w})^2\}\big]\big\}\\
= & \max\limits_{\psi \geq 0} \; \sup\limits_{t\geq 0}\;\big\{-t\theta + \mathbb{E}_{H}\big[\min\limits_{u\geq 0}\;\{\pi(\psi,u) + \frac{\alpha t}{\alpha + t}(u-\tilde{w})^2\}\big]\big\},
\end{array}
\]
where the equality follows since $\min_{v\geq 0}\{\alpha(u-v)^2 + t(v-w)^2\} = \frac{\alpha t}{\alpha + t}(u-w)^2$ for any $u \geq 0$ and $w \geq 0$. Interchanging the ``max'' over $\psi$ and $t$ and applying the variable substitution $\gamma \leftarrow \frac{\alpha t}{\alpha + t}$, we arrive at
\begin{equation}\label{equ:dual OT}
\textstyle \max\limits_{0 \leq  \gamma \leq \alpha}\;\big\{-\frac{\alpha \gamma}{\alpha - \gamma}\theta + \max\limits_{\psi \geq 0} \; \mathbb{E}_{H}\big[\min\limits_{u \geq 0}\;\big\{\pi(\psi,u) + \gamma(u-\tilde{w})^2\big\}\big]\big\}.
\end{equation}
Here, $\gamma = \alpha$ corresponds to $t = \infty$, and $\lim_{\gamma \to \alpha-} -\frac{\alpha \gamma}{\alpha - \gamma}\theta = 0$ if $\theta = 0$ and $\lim_{\gamma \to \alpha-} -\frac{\alpha \gamma}{\alpha -\gamma}\theta = \infty$ if $\theta > 0$. Hence, there must exists $\gamma^\star \in [0, \alpha]$ such that problem~\eqref{equ:dual OT} can be equivalently solved by
\[
\begin{array}{ll}
\displaystyle \max_{\psi \geq 0} \; \mathbb{E}_{H}\bigg[\min_{u\geq 0}\;\big\{\pi(\psi,u) + \gamma^\star \cdot (u-\tilde{w})^2\big\}\bigg] & = \displaystyle \max_{\psi \geq 0} \; \min_{F\in\mathcal{P},\;\Gamma\in\mathcal{W}(F,H)}\mathbb{E}_{\Gamma}[\pi(\psi,\tilde{u})+\gamma^\star \cdot (\tilde{u}-\tilde{w})^2]\\ 
& = \displaystyle \max_{\psi \geq 0} \; \min_{F\in\mathcal{P}} \; \{\mathbb{E}_F[\pi(\psi, \tilde{u})] + \gamma^\star \cdot d_{\rm OT}(F,H)\},
\end{array}
\]
where the first and second equalities follow from \Cref{lemma:interchangeability} and the definition of $d_{\rm OT}(F,H)$, respectively. Hence, problem~\eqref{prob:double misspecification} is equivalent to $\max_{\psi \geq 0} \min_{F \in \mathcal{P}} \{\mathbb{E}_{F}[\pi(\psi,\tilde{u})] + \gamma^\star \cdot d_{\rm OT}(F,H)\}$ for some $\gamma^\star \in [0, \alpha]$.

In the remainder of the proof, we solve the optimal $\gamma^\star$ and $\psi^\star_{\gamma^\star}$ of problem \eqref{equ:dual OT}. For any fixed $w \geq 0$, define $\Psi(\gamma,q,w) = \min_{u \geq 0} \{\pi(q,u) + \gamma(u-w)^2\}$. Problem \eqref{equ:dual OT} is then equivalent to
\begin{equation}\label{equ:dual OT reformulation}
\textstyle \max\limits_{0 \leq \gamma \leq \alpha} \; \big\{-\frac{\alpha \gamma}{\alpha - \gamma}\theta + \max\limits_{\psi \geq 0}\;\mathbb{E}_{H}\big[\Psi(\gamma, \psi, \tilde{w})\big]\big\}.
\end{equation}

We first solve, by the first-order optimality condition, the inner maximization of problem~\eqref{equ:dual OT reformulation} given $\gamma \in [0,\alpha]$. Let $Z(\psi) = \mathbb{E}_{H}[\Psi(\gamma,\psi,\tilde{w})]$. Given $w\geq 0$, by the concavity of $\Psi(\gamma,\psi,w)$, $Z(\psi)$ is also concave with a decreasing derivative. In particular, if $\psi<\frac{p}{4\gamma}$, we have $Z'(\psi) = p \cdot \mathbb{P}_{H}\{\gamma  \tilde{w}^2\geq p\psi\}-c$; if $\psi \geq \frac{p}{4\gamma}$, we have $Z'(\psi) = p \cdot \mathbb{P}_{H}\{\tilde{w}-\frac{p}{4\gamma}\geq \psi\} - c$. Note that $Z'(0) = p-c > 0$ and $\lim_{\psi\to\infty}Z'(\psi)=-c<0$. Let $\gamma_0 = \frac{p}{2 q^\star_H}$. If $\gamma<\gamma_0$, then $Z'(\frac{p}{4\gamma}) = p(\kappa - H(\frac{p}{2\gamma}))<0$. Hence, the maximum of $Z(\psi)$ is attained in $[0,\frac{p}{4\gamma}]$. Setting the derivative to $0$ yields
\begin{equation}\label{equ:q(z) 1}
\textstyle \psi^\star_\gamma = (q^\star_H)^2\cdot\frac{\gamma}{p} < \frac{p}{4\gamma}.
\end{equation}
If $\gamma \geq \gamma_0$, then $Z'(\frac{p}{4\gamma}) = p(\kappa - H(\frac{p}{2\gamma}))\geq 0$, implying that the maximum of $Z(\psi)$ must be attained in $[\frac{p}{4\gamma},\infty)$. Setting the derivative to $0$ then yields
\begin{equation}\label{equ:q(z) 2}
\textstyle \psi^\star_\gamma = q^\star_H - \frac{p}{4\gamma} \geq \frac{p}{4\gamma}.
\end{equation}

We next solve the optimal $\gamma^\star$ of the outer maximization of problem~\eqref{equ:dual OT reformulation}, whose objective is
\[
\textstyle Q(\gamma) = -\frac{\alpha \gamma}{\alpha - \gamma}\theta + \max\limits_{\psi\geq 0}\;\mathbb{E}_{H}\big[\Psi(\gamma, \psi, \tilde{w})\big] = -\frac{\alpha \gamma}{\alpha -\gamma}\theta + \mathbb{E}_{H}\big[\Psi(\gamma, \psi^\star_\gamma, \tilde{w})\big].
\]
Note that given $w \geq 0$, $\Psi(\gamma, \psi, w)$ is jointly concave in $(\gamma, \psi)$. Since concavity is preserved under expectation and maximization, $Q(\gamma)$ is also concave in $\gamma$. When $\gamma < \gamma_0$, plugging the expressions of $\psi^\star_\gamma$ in~\eqref{equ:q(z) 1} and $\Psi(\alpha,\psi,w)$ in~\Cref{lemma:interchangeability}, we have
\[
\textstyle Q(\gamma) = -\frac{\alpha \gamma}{\alpha - \gamma}\theta + \mathbb{E}_{H}\big[\Psi(\gamma, \psi^\star_\gamma, \tilde{w})\big] = -\frac{\alpha \gamma}{\alpha - \gamma}\theta + \int_{0}^{q^\star_H}zw^2{\rm d}H(w), ~~Q'(\gamma) = -\frac{\alpha^2}{(\alpha - \gamma)^2}\theta + \int_{0}^{q^\star_H}w^2{\rm d}H(w).
\]
When $\gamma \geq \gamma_0$, plugging the expression of $\psi^\star_\gamma$ in~\eqref{equ:q(z) 2}, we have
\[
\textstyle Q(\gamma) = -\frac{\alpha \gamma}{\alpha - \gamma}\theta + \int_{0}^{\frac{p}{2\gamma}}\gamma w^2{\rm d}H(w) + \int_{\frac{p}{2z}}^{q^\star_H} p(w-\frac{p}{4\gamma}){\rm d}H(w), ~~Q'(\gamma) = -\frac{\alpha^2\theta}{(\alpha-\gamma)^2} + \frac{p^2}{4 \gamma^2}(\kappa-H(\frac{p}{2 \gamma})) + \int_0^{\frac{p}{2\gamma}}w^2{\rm d}H(w).
\]
Note that $Q'(\gamma)$ is always decreasing in $\gamma$. When $\theta = 0$, it is straightforward to see that $Q'(\gamma) \geq 0$. Hence, $\gamma^\star = \alpha$ (which corresponds to $t^\star = \infty$ in the dual reformulation). In the following, we focus on $\theta > 0$. Based on the sign of $Q'(0)$, we divide the problem into two scenarios: $\theta \geq \beta$ and $\theta < \beta$, with $\beta = \int_{0}^{q^\star_H}u^2{\rm d}H(u)$. 

For the former scenario of $\theta \geq \beta$, $Q'(\gamma)\leq Q'(0) = -\theta + \beta \leq 0$ for any $\gamma \in [0,\alpha]$. Hence, the maximum of $Q(\gamma)$ is attained at $\gamma^\star=0$. Equation~\eqref{equ:q(z) 1} then yields $\psi^\star_{\gamma^\star} = 0$. This corresponds to the case~(\textit{i}) in the statement. For the latter scenario of $\theta < \beta$, $Q'(0) = -\theta + \beta > 0$. To proceed, we further consider two situations based on the sign of $Q'(\gamma_0)$. 
If $1-\sqrt{\theta/\beta}<\frac{\gamma_0}{\alpha}$, \textit{i.e.}, $\sqrt{\theta/\beta}<\frac{\alpha}{\alpha - \gamma_0}$, then $Q'(\lambda_0) = \beta - \frac{\alpha^2}{(\alpha - \gamma_0)^2}\theta<0$.
Since $Q'(0) = -\theta + \beta>0$, the maximum of $Q(\gamma)$ must be attained in $[0,\gamma_0]$. Setting the derivative of $Q(\gamma)$ to $0$ yields $\gamma^\star=\alpha(1-\sqrt{\theta/\beta}) < \gamma_0$. By equation~\eqref{equ:q(z) 1} we have $\psi^\star_{\gamma^\star} = (q^\star_H)^2(1-\sqrt{\theta/\beta})\cdot\frac{\alpha}{p}$.
This corresponds to case~(\textit{ii}) in the statement. If $1-\sqrt{\theta/\beta}\geq\frac{\gamma_0}{\alpha}$, then $Q'(\gamma_0)\geq 0$. Moreover, $\theta > 0$ gives $Q'_-(\alpha)=-\infty$, which further implies that there exists some $\gamma^\star\in [\gamma_0,\alpha)$ such that $-\frac{\alpha^2\theta}{(\alpha-\gamma^\star)^2} + \frac{p^2}{4 \gamma^{\star2}}\big(\kappa-H\big(\frac{p}{2 \gamma^\star}\big)\big) + \int_0^{\frac{p}{2z^\star}}w^2{\rm d}H(w) = 0$.
Equation~\eqref{equ:q(z) 2} then yields $\psi^\star_{\gamma^\star} = q^\star_H - \frac{p}{4\gamma^\star}$. This corresponds to case~(\textit{iii}) in the statement.
\hfill $\square$
\vspace{2mm}

\noindent{\bf Proof of \Cref{thm:newsvendor TV}.} $\;$
If $\alpha = 0$, it is easy to see $q_\alpha^\star = 0$, and the result follows. We next focus on $\alpha > 0$. By the dual representation of total variation distance (see, {\it e.g.}, \citealp{muller1997integral}), we have $d_{\rm TV}(F, G) = \max_{\|\chi\|_\infty \leq 1}\{\mathbb{E}_F[\chi(v)] - \mathbb{E}_G[\chi(v)]\}$. Hence, for any $q \geq 0$ and $G \in \mathcal{A}$,
\[
\begin{array}{rl}
\displaystyle \min_{F \in \mathcal{P}}\{\mathbb{E}_F[\pi(q, \tilde{v})] + \alpha\cdot d_{\rm TV}(F, G)\} & \displaystyle = \max_{\|\chi\|_\infty \leq 1}\min_{F \in \mathcal{P}}\big\{\mathbb{E}_F[\pi(q, \tilde{v})] + \alpha\mathbb{E}_F[\chi(\tilde{v})] - \alpha\mathbb{E}_G[\chi(\tilde{v})]\big\}\\
& \displaystyle = \max_{\|\chi\|_\infty \leq 1}\Big\{\min_{v \geq 0}\{\pi(q, v)] + \alpha\chi(v)\} - \alpha\mathbb{E}_G[\chi(\tilde{v})]\Big\},
\end{array}
\]
where the first equality follows from the minimax theorem. Introducing an epigraphical variable $\tau$, the above problem becomes
\[
\begin{array}{rl}
& \displaystyle \max_{\tau, \chi(\cdot)}\Big\{\tau - \alpha \cdot \mathbb{E}_G[\chi(\tilde{v})]~\big|~ \pi(q, v) + \alpha\chi(v) \geq \tau, \;|\chi(v)| \leq 1,\;\forall v \geq 0\Big\}\\
= &  \max\limits_{\tau}\big\{\tau - \alpha\cdot \min\limits_{\chi(\cdot)}\mathbb{E}_G[\chi(\tilde{v})]~\big|~ \chi(v) \geq \frac{\tau - \pi(q, v)}{\alpha},\; -1 \leq \chi(v) \leq 1,\;\forall v \geq 0\big\}.
\end{array}
\]
Given $\tau$, we investigate the minimization of $\chi(\cdot)$. First, to ensure that the feasible set of $\chi(\cdot)$ is nonempty, we must have $\frac{\tau - \pi(q, v)}{\alpha} \leq 1$ for all $v \geq 0$, and hence $\tau \leq \alpha - cq$. Second, it suffices to focus on the pointwise minimum of $\chi(\cdot)$, and it is straightforward to see that the optimum is attained at $\chi^\star(v) = \max\{-1, \frac{\tau - \pi(q, v)}{\alpha}\}$ for any $v \geq 0$. Plugging in the expression of $\chi^\star$ into the objective, problem~\eqref{prob:newsvendor TV} then becomes
\[
\max_{q \geq 0}\;\min_{G \in \mathcal{A}}\;\max_{\tau \leq \alpha - cq}\big\{\tau - \mathbb{E}_G[\max\{-\alpha, \tau - \pi(q, \tilde{v})\}].
\]
We first look at the inner maximization over $\tau$. For any $\tau \geq \alpha + cq$ (resp., $\tau \leq \alpha -(p-c)q$), $\min\{\alpha - \tau + cq, pv, pq\} = 0$ (resp., $= \min\{pv, pq\}$) and hence,, the maximum over $\tau$ must be attained in $[\alpha - (p-c)q, \alpha + cq]$. This, together with the prerequisite $\tau \geq cq - \alpha$, implies that it suffices to focus on $0 \leq \alpha - \tau + cq \leq \min\{pq, 2\alpha\}$. 

We next look at the optimization over $q \geq 0$:
\[
\begin{array}{rl}
\displaystyle\max_{q \geq 0} \; \min_{G \in \mathcal{A}} \; \max_{0 \leq \alpha - \tau + cq \leq  \min\{pq, 2\alpha\}} \; \{\mathbb{E}_{G}[\min\{\alpha - \tau + cq, p\tilde{v}\}]-cq\} =&\displaystyle\max_{q \geq 0} \; \min_{G\in\mathcal{A}} \; \max_{0 \leq \tau \leq \min\{pq, 2\alpha\}} \; \{\mathbb{E}_{G}[\min\{\tau, p\tilde{v}\}]-cq\}\\ [3mm]
=&\displaystyle
\max_{q \geq 0} \; \min_{G\in\mathcal{A}} \; \{\mathbb{E}_{G}[\min\{pq, 2\alpha,  p\tilde{v}\}] - cq\},
\end{array}
\]
where the first equality follows from the variable substitution $\tau \leftarrow \alpha - \tau + cq$ and the second equality follows from the fact that $\mathbb{E}_{G}[\min\{\tau, p\tilde{v}\}]$ is increasing in $\tau$ so its maximum is attained at $\min\{pq, 2\alpha\}$. When $q \geq \frac{2\alpha}{p}$, the objective function $\mathbb{E}_{G}[\min\{pq, 2\alpha,  p\tilde{v}\}] - cq = \mathbb{E}_{G}[\min\{2\alpha,  p\tilde{v}\}] - cq$ is decreasing in $q$ given $G \in \mathcal{A}$. Hence, it is optimal to set $q$ to $\frac{2\alpha}{p}$ whenever $q \geq \frac{2\alpha}{p}$. That is to say, the optimal order quantity of problem~\eqref{prob:newsvendor TV} must reside in $[0,\frac{2\alpha}{p}]$. For this interval, we have
\[
\textstyle\max\limits_{0 \leq q \leq \frac{2\alpha}{p}} \; \min\limits_{G\in\mathcal{A}} \; \{\mathbb{E}_{G}[\min\{pq, 2\alpha,  p\tilde{v}\}] - cq\} = \max\limits_{0 \leq q \leq \frac{2\alpha}{p}} \; \min\limits_{G\in\mathcal{A}} \; \mathbb{E}_{G}[\pi(q, \tilde{v})].
\]
Hence, it suffices to solve the right-hand side problem. We note that $\Phi(q) = \min_{G\in\mathcal{A}} \mathbb{E}_{G}[\pi(q,\tilde{v})]$ is concave in $q$. Hence, $q_\alpha^\star = \min\{\frac{2\alpha}{p}, q_\infty^\star\}$. Combining the expression of $q_\infty^\star$ in \eqref{equ:optimal order ambiguity} then completes the proof.
\hfill $\square$

\section{Equivalence between $\mathcal{V}$ and $\mathcal{A}$ for \ref{prob:misspecification}}\label{appendix:equivalence}
We show that for \ref{prob:misspecification}, it is equivalent to consider $\mathcal{V} = \{G \in \mathcal{P} \mid \mathbb{E}_{G}[\tilde{v}] = \mu,\;\mathbb{E}_{G}[\tilde{v}^2] = \mu^2+\sigma^2 ~~\mbox{for some}~~ \mu \in [\ubar{\mu},\bar{\mu}] ~\mbox{and}~ \sigma^2 \in [\ubar{\sigma}^2,\bar{\sigma}^2]\}$ with uncertain mean and variance and and $\mathcal{A} = \{G \in \mathcal{P} \mid \mathbb{E}_{G}[\tilde{v}] = \ubar{\mu},\;\mathbb{E}_{G}[\tilde{v}^2] = \ubar{\mu}^2+\bar{\sigma}^2\}$ with the lowest mean and the highest variance.

\begin{lemma}
For any $\alpha \geq 0$, it holds that
\[
\max_{q \geq 0}\;\min_{F \in \mathcal{P}} \; \big\{\mathbb{E}_{F}[\pi(q,\tilde{u})] + \alpha \cdot d(F,\mathcal{V})\big\} = \max_{q \geq 0}\;\min_{F \in \mathcal{P}} \; \big\{\mathbb{E}_{F}[\pi(q,\tilde{u})] + \alpha \cdot d(F,\mathcal{A})\big\}.
\]
\end{lemma}

\noindent{\bf Proof.} $\;$
By the definition of $\mathcal{V}$, for any $q \geq 0$, the left-hand side objective for a fixed $q$ can be rewritten as
\[
\min_{\mu \in [\ubar{\mu},\bar{\mu}],\sigma^2 \in [\ubar{\sigma},\bar{\sigma}]}\;\min_{F \in \mathcal{P}} \; \Big\{\mathbb{E}_{F}[\pi(q,\tilde{u})] + \alpha \cdot d(F,\mathcal{A}(\mu, \sigma))\Big\} = \min_{\mu \in [\ubar{\mu},\bar{\mu}],\sigma^2 \in [\ubar{\sigma},\bar{\sigma}]} \; L(\mu, \sigma, q),
\]
where $\mathcal{A}(\mu, \sigma) = \{G \in \mathcal{P} \mid \mathbb{E}_{G}[\tilde{v}] = \mu,\;\mathbb{E}_{G}[\tilde{v}^2] = \mu^2+\sigma^2\}$ and the expression of $L(\mu, \sigma; q)$ has been characterized in \Cref{prop:transformed expectation}. Then it suffices to show that $L(\mu, \sigma, q)$ is always increasing in $\mu$ and decreasing in $\sigma$. Note that the function $L(\mu, \sigma, q)$ consists of two pieces: $L_1(\mu, \sigma, q)=\frac{p}{2}(q+\mu - \frac{p}{4\alpha} - \sqrt{(q-\mu + \frac{p}{4\alpha})^2 + \sigma^2}) - cq$ and $L_2(\mu, \sigma, q)=\frac{\alpha}{2}(\frac{pq}{\alpha}+\mu^2+\sigma^2-\sqrt{(\frac{pq}{\alpha}+\mu^2+\sigma^2)^2-4\mu^2\frac{pq}{\alpha}}) - cq$. On the one hand, it is straightforward to see that $L_1(\mu, \sigma, q)$ is decreasing in $\sigma$, and we also have $\frac{\partial L_1(\mu, \sigma, q)}{\partial \mu} = \frac{p}{2}\big(1 + \frac{q - \mu + \frac{p}{4\alpha}}{\sqrt{(q-\mu + \frac{p}{4\alpha})^2 + \sigma^2}}\big) \geq 0$.
On the other hand, we note that $L_2(\mu, \sigma, p) = \frac{\alpha}{2}\frac{4\mu^2\frac{pq}{\alpha}}{\frac{pq}{\alpha}+\mu^2+\sigma^2 + \sqrt{(\frac{pq}{\alpha}+\mu^2+\sigma^2)^2-4\mu^2\frac{pq}{\alpha}}}$ is decreasing in $\sigma$. Additionally, $\frac{\partial L_2(\mu, \sigma, q)}{\partial \mu} = \alpha\mu\big(1 - \frac{\frac{pq}{\alpha}+\mu^2+\sigma^2 - 2\frac{pq}{\alpha}}{\sqrt{(\frac{pq}{\alpha}+\mu^2+\sigma^2)^2-4\mu^2\frac{pq}{\alpha}}}\big) \geq 0$. Finally, one can easily verify that $L(\mu, \sigma, q)$ is continuous in $\mu$ and $\sigma$, and hence, $L(\mu, \sigma, q)$ is increasing in $\mu$ and decreasing in $\sigma$, concluding the proof. 
\hfill $\square$

\section{Transformed Worst-Case  Distribution}\label{appendix:worst-case transformed distribution}
In this section, we derive the worst-case transformed distribution $T_{\varphi_\alpha}[G_\alpha^\star]$ given an order quantity $q \geq 0$. For ease of notation, let $w(q) = \frac{pq}{\alpha} + \mu^2 + \sigma^2$, $u(q) = q + \frac{p}{4\alpha}$, and define
\[
\left\{
\begin{array}{ll}
v_1 = \frac{1}{2\mu} \big(w(q) - \sqrt{w(q)^2 - 4\mu^2\frac{pq}{\alpha}}\big),\; v_2 = \frac{1}{2\mu} \big(w(q) + \sqrt{w(q)^2 - 4\mu^2\frac{pq}{\alpha}}\big)\\ [2mm]
v_3 = u(q) - \sqrt{(u(q)-\mu)^2+\sigma^2},\; v_4 = u(q) + \sqrt{(u(q)-\mu)^2+\sigma^2}.
\end{array}
\right.
\]

\begin{proposition}[\sc{Worst-Case Transformed Distribution}]\label{prop:worst-case transformed distribution}
Given $q \geq 0$ and $\alpha \geq 0$, the worst-case transformed distribution $T_{\varphi_\alpha}[G_\alpha^\star]$ of \ref{prob:transform} is given by
\[
T_{\varphi_\alpha}[G_\alpha^\star] = \left\{
\begin{array}{ll}
\Big(\frac{1}{2} - \frac{\mu^2 - \sigma^2 - \frac{pq}{\alpha}}{2\sqrt{w(q)^2 - 4\mu^2\frac{pq}{\alpha}}}\Big) \cdot \delta_{\frac{\alpha}{p} v_1^2} + \Big(\frac{1}{2} + \frac{\mu^2 - \sigma^2 - \frac{pq}{\alpha}}{2\sqrt{w(q)^2 - 4\mu^2\frac{pq}{\alpha}}}\Big) \cdot \delta_{\frac{\alpha}{p} v_2^2} & 0 \leq q \leq \frac{p}{4\alpha}\\ [3mm]
\Big(\frac{1}{2} - \frac{\mu^2 - \sigma^2 - \frac{pq}{\alpha}}{2\sqrt{w(q)^2 - 4\mu^2\frac{pq}{\alpha}}}\Big) \cdot \delta_{\frac{\alpha}{p} v_1^2} + \Big(\frac{1}{2} + \frac{\mu^2 - \sigma^2 - \frac{pq}{\alpha}}{2\sqrt{w(q)^2 - 4\mu^2\frac{pq}{\alpha}}}\Big) \cdot \delta_{v_2 - \frac{p}{4\alpha}} & q\geq \frac{p}{4\alpha},(2\mu-\frac{p}{\alpha})q < \mu^2 + \sigma^2 - \frac{p\mu}{2\alpha}\\ [3mm]
\frac{1}{2}\Big(1+\frac{u(q)-\mu}{\sqrt{(u(q)-\mu)^2+\sigma^2}}\Big)\cdot\delta_{v_3 - \frac{p}{4\alpha}} + \frac{1}{2}\Big(1-\frac{u(q)-\mu}{\sqrt{(u(q)-\mu)^2+\sigma^2}}\Big)\cdot\delta_{v_4 - \frac{p}{4\alpha}} & q\geq \frac{p}{4\alpha},(2\mu-\frac{p}{\alpha})q < \mu^2 + \sigma^2 - \frac{p\mu}{2\alpha}.
\end{array}
\right.
\]
\end{proposition}

\noindent{\bf Proof.} $\;$ Given $q \geq 0$, we have derived the expression of the worst-case distribution $G_\alpha^\star$ as in \eqref{equ:primal solution 1} and \eqref{equ:primal solution 2}. Plugging the expression of $\varphi_\alpha$, we then obtain the desired results. 
\hfill $\square$

\begin{figure}[tb]
\begin{subfigure}{0.33\textwidth}
\centering
\includegraphics[width=1.11\linewidth]{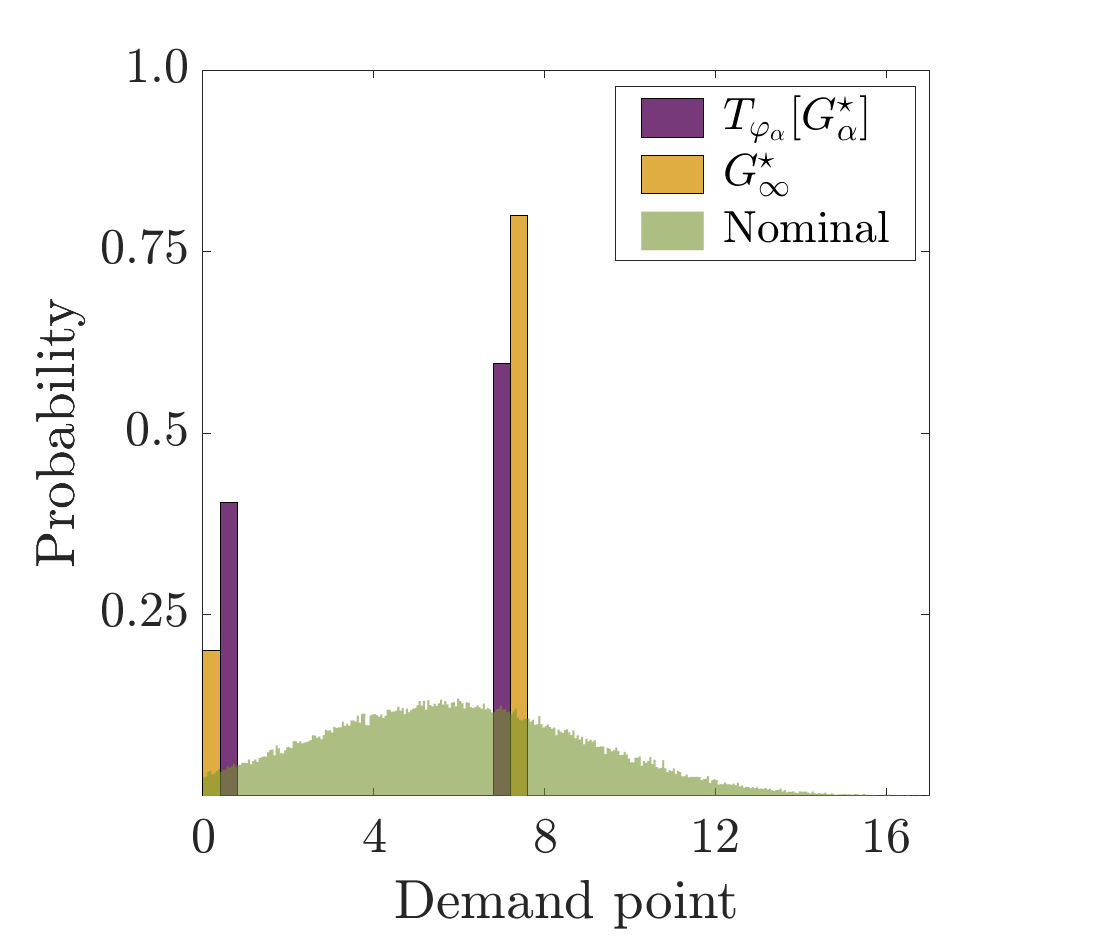}
\end{subfigure}
\begin{subfigure}{0.33\textwidth}
\centering
\includegraphics[width=1.11\linewidth]{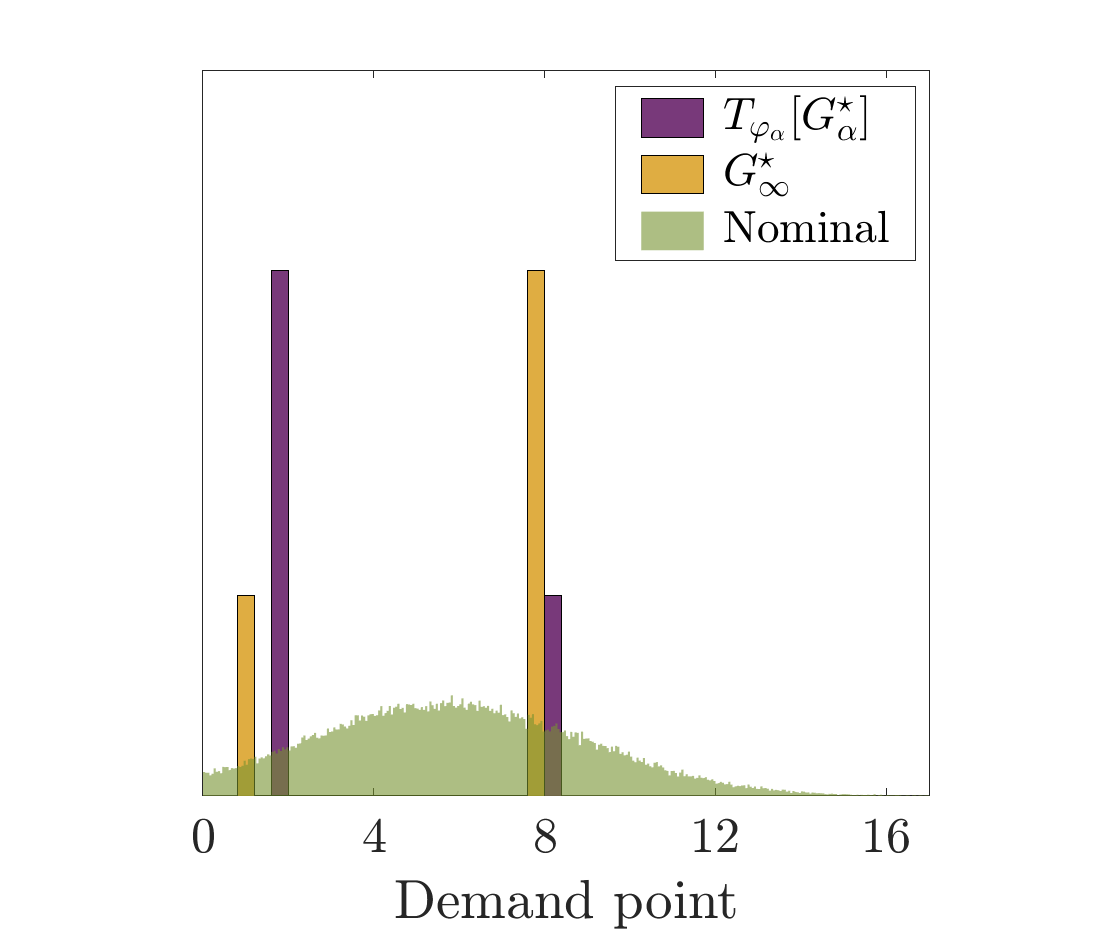}
\end{subfigure}
\begin{subfigure}{0.33\textwidth}
\centering
\includegraphics[width=1.11\linewidth]{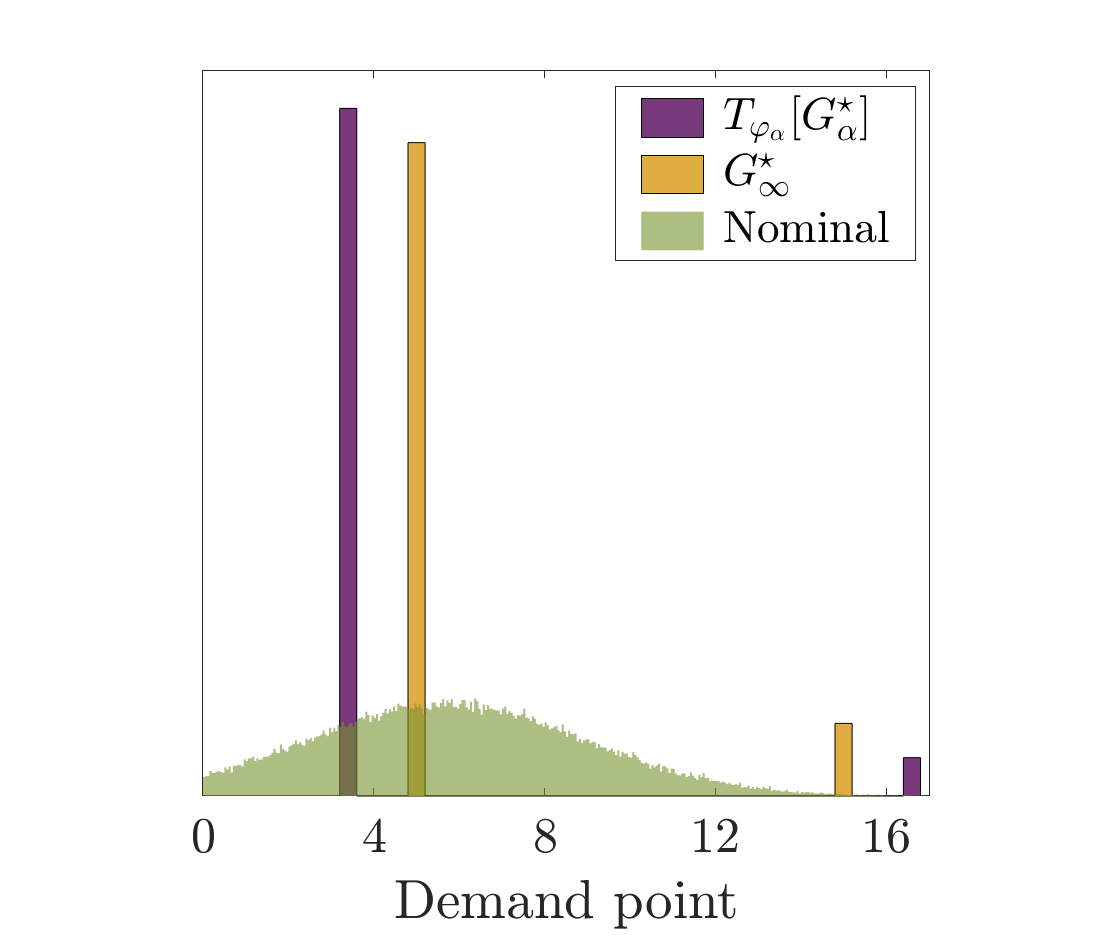}
\end{subfigure}
\caption{\textnormal{Three scenarios depending on the value of $q$: $0 \leq q \leq \frac{p}{4\alpha}$ (left), $q\geq \frac{p}{4\alpha}$ and $(2\mu-\frac{p}{\alpha})q < \mu^2 + \sigma^2 - \frac{p\mu}{2\alpha}$ (middle), and $q\geq \frac{p}{4\alpha}$ and $(2\mu-\frac{p}{\alpha})q \geq \mu^2 + \sigma^2 - \frac{p\mu}{2\alpha}$ (right).}}
\label{fig:transformed dist}
\vspace{-6mm}
\end{figure}

In~\Cref{fig:transformed dist}, we provide a visualization of the transformed worst-case distribution $T_{\varphi_\alpha}[G_\alpha^\star]$, pertaining to the worst-case distribution $G_\infty^\star$ for \ref{prob:ambiguity} and the nominal distribution as truncated normal. Define $v_1^\star = \mu - \sigma\sqrt{\frac{c}{p-c}}$ and $v_2^\star = \mu + \sigma\sqrt{\frac{p-c}{c}}$. Plugging in the expression of $q_\alpha^\star$ in \Cref{thm:newsvendor transport-2}, the worst-case transformed distribution $T_{\varphi_\alpha}[G_\alpha^\star]$ can be characterized as follows.
\begin{enumerate}
\item[($i$)] If $\frac{p - c}{p} < \frac{\sigma^2}{\mu^2 + \sigma^2}$, then $T_{\varphi_\alpha}[G_\alpha^\star] = \big(\frac{\sigma^2}{\mu^2 + \sigma^2}\big)\cdot\delta_0 + \big(\frac{\mu^2}{\mu^2 + \sigma^2}\big)\cdot\delta_{\frac{\mu^2+\sigma^2}{\mu}}$.

\item[($ii$)] If $\frac{p - c}{p} \geq \frac{\sigma^2}{\mu^2 + \sigma^2}$, then
\[
T_{\varphi_\alpha}[G_\alpha^\star] = 
\left\{
\begin{array}{ll} \big(\frac{p-c}{p}\big)\cdot\delta_{\frac{\alpha}{p}v_1^{\star 2}} + \big(\frac{c}{p}\big)\cdot\delta_{\frac{\alpha}{p}v_2^{\star2}} & \alpha < p/(2\sqrt{v_1^\star v_2^\star})\\
\big(\frac{p-c}{p}\big)\cdot\delta_{\frac{\alpha}{p}v_1^{\star2}} + \big(\frac{c}{p}\big)\cdot\delta_{v_2^\star - \frac{p}{4\alpha}} & p/(2\sqrt{v_1^\star v_2^\star}) \leq \alpha < p/(2v_1^\star)\\
\big(\frac{p-c}{p}\big)\cdot\delta_{v_1^\star - \frac{p}{4\alpha}} + \big(\frac{c}{p}\big)\cdot\delta_{v_2^\star - \frac{p}{4\alpha}} & \alpha \geq p/(2v_1^\star).
\end{array}
\right.
\]
\end{enumerate}

Setting $\alpha \to \infty$ then yields the expression of the worst-case distribution of \ref{prob:ambiguity}: $G_\infty^\star = (\frac{\sigma^2}{\mu^2 + \sigma^2})\cdot\delta_0 + (\frac{\mu^2}{\mu^2 + \sigma^2})\cdot\delta_{\frac{\mu^2+\sigma^2}{\mu}}$ if $\frac{p - c}{p} < \frac{\sigma^2}{\mu^2 + \sigma^2}$
and $G_\infty^\star = (\frac{p-c}{p})\cdot\delta_{v_1^\star} + (\frac{c}{p})\cdot\delta_{v_2^\star}$ otherwise.
It is evident that $G_\infty^\star$ and $T_{\varphi_\alpha}[G_\alpha^\star]$ are both two-point distributions with the same probability mass, and the minimal (resp., maximal) support point of $T_{\varphi_\alpha}[G_\alpha^\star]$ is always smaller than that of $G_\infty^\star$.

\section{Generality of \Cref{thm:multi product}}
\label{EC:generalization}

We show how \Cref{thm:multi product} generalizes \Cref{thm:newsvendor transport-2} for a single product to multiple products. To see this, let $K = \mu^2 + \sigma^2$ and note that for the single-product problem, $\lambda^\star$ is the optimal dual variable $r^\star_\alpha$ in the dual reformulation. Plugging the expression of $q_\alpha^\star$ into equations~\eqref{equ:dual solution 1} and \eqref{equ:dual solution 2}, it is immediate to see that
\[
r_\alpha^\star=\left\{
\begin{array}{ll}
\sqrt{\kappa(1-\kappa)}/(2p\sigma) &~~~ \kappa \geq \frac{\sigma^2}{\mu^2+\sigma^2}, \; \alpha \geq p/(2(\mu - \sigma \sqrt{(1-\kappa)/\kappa})) \\ 
\alpha\sqrt{\kappa(1-\kappa)}\big(\mu - \sigma\sqrt{(1-\kappa)/\kappa}\big)/\sigma &~~~ \kappa \geq \frac{\sigma^2}{\mu^2+\sigma^2}, \; \alpha < \frac{p}{2(\mu - \sigma \sqrt{(1-\kappa)/\kappa})}\\ 
0 &~~~\kappa < \frac{\sigma^2}{\mu^2+\sigma^2}.
\end{array}
\right.
\]
If $\kappa \geq \frac{\sigma^2}{\mu^2 + \sigma^2}$ and $\alpha \geq \frac{p}{2(\mu-\sigma\sqrt{(1-\kappa)/\kappa})}$, it holds that $\lambda^\star = \frac{\sqrt{\kappa(1-\kappa)}}{2p\sigma}$ and hence, $\frac{p}{2(\mu - c/(2\lambda^\star))^+} = \frac{p}{2(\mu-\sigma\sqrt{(1-\kappa)/\kappa})}$. By \Cref{thm:multi product}, we have $q_\alpha^\star = \mu + \frac{p-2c}{4\lambda^\star} - \frac{p}{4\alpha} = \mu + \sigma\frac{2\kappa-1}{2\sqrt{\kappa(1-\kappa)}} - \frac{p}{4\alpha}$.
If $\kappa \geq \frac{\sigma^2}{\mu^2 + \sigma^2}$ and $\alpha < \frac{p}{2(\mu-\sigma\sqrt{(1-\kappa)/\kappa})}$, then $\lambda^\star = \frac{\alpha\sqrt{\kappa(1-\kappa)}}{\sigma}(\mu - \sigma\sqrt{(1-\kappa)/\kappa})$. In this case, we can verify that $\alpha < \frac{p}{2(\mu - c/(2\lambda^\star))^+}$, thus, by \Cref{thm:multi product}, $q_\alpha^\star = \frac{\lambda^\star (\lambda^\star + \alpha) p \mu^2}{\alpha(\lambda^\star p/\alpha + c)^2} = (\mu^2 - \sigma^2 + 2 \mu\sigma f( 1-\kappa))\cdot \frac{\alpha}{p}$.
If $\kappa < \frac{\sigma^2}{\mu^2 + \sigma^2}$, we always have $\lambda^\star = 0$, leading to $q_\alpha^\star = 0$.
Consolidating these three scenarios then recovers the optimal solution~\eqref{equ:quantity misspecification} in \Cref{thm:newsvendor transport-2}.

\section{Multiple Products with Complete Covariance Information}\label{appendix:multi product}

In line with the setting of Section \ref{sec:multi product}, we study a misspecification-averse multi-product newsvendor with mean and complete covariance information. Consider $M$ products with random demands $\bmt{u} = (\tilde{u}_1, \dots, \tilde{u}_i) \sim F$ following a multi-dimensional distribution $F$. The misspecification-averse newsvendor then solves
\begin{equation}\label{prob:covariance}
\max_{\bm{q} \geq \bm{0}} \; \min_{F \in \mathcal{P}_M} \; \{\mathbb{E}_{F}[\omega(\bm{q},\bmt{u})] + \alpha \cdot d(F,\mathcal{A})\},
\end{equation}
where the function $\omega(\bm{q}, \bm{u}) = \sum_{i = 1}^M \pi_i(q_i,u_i)$, the optimal-transport cost $d_{\rm OT}(\cdot,\cdot)$ is defined in \eqref{equ:optimal transport 2} with $L_2$-norm, and the ambiguity set is specified by mean and complete covariance information, \textit{i.e.},
\begin{equation}\label{set:covariance}
\mathcal{A} = \big\{G \in\mathcal{P}_M \mid \mathbb{E}_{G}[\bmt{v}]=\bm{\mu},\;\mathbb{E}_{G}[(\bmt{v}-\bm{\mu})(\bmt{v}-\bm{\mu})^\top]= \bm{\Sigma}\big\}.
\end{equation}

\begin{proposition}\label{prop:multi product}
Let $\mathcal{S}^M$ be the space of symmetric matrices with dimension $M$. For the ambiguity set~\eqref{set:covariance} with mean and covariance information, the misspecification-averse multi-product newsvendor problem \eqref{prob:covariance} is equivalent to
\begin{equation}\label{equ:multi product reformulation}
\begin{array}{cll}
\max & t + \bm{\lambda}^\top\bm{\mu}+ \left<\bm{Q}, \bm{\Sigma}+\bm{\mu}\bm{\mu}^\top\right>\\ [1mm]
{\rm s.t.} & t + \bm{\lambda}^\top\bm{v} + \left<\bm{Q}, \bm{v}\bm{v}^\top\right> \leq \sum\limits_{i\in[M]}\Psi_i(\alpha,q_i,v_i) &\forall \bm{v}\geq 0\\ [1mm]
&\bm{q}\geq 0,\;\bm{Q}\in\mathcal{S}^M,\;\bm{\lambda}\in\mathbb{R}^M,\;t\in\mathbb{R},
\end{array}
\end{equation}
where for each $i \in [M]$, $\Psi_i(\alpha,q_i,v_i) = \min\{\alpha v_i^2, pq_i\} - c_iq_i$ if $0\leq q_i \leq \frac{p_i}{4\alpha}$, and if $q_i > \frac{p_i}{4\alpha}$,
\[
\Psi_i(\alpha,q_i,v_i)=\left\{
\begin{array}{ll}
\alpha v_i^2-c_iq_i & 0 \leq v_i \leq \frac{p_i}{2\alpha}\\ [1mm]
p\cdot\min\{v_i - \frac{p_i}{4\alpha}, q_i\} - c_iq_i & v_i > \frac{p_i}{2\alpha}.
\end{array}
\right.
\]
\end{proposition}

\noindent{\bf Proof.} $\;$ For $\bm{q} \geq 0$, by the interchangeability principle in \Cref{lemma:interchangeability}, we have
\[
\max_{\bm{q} \geq 0}\;\min_{G \in \mathcal{C}}\;\min_{F \in \mathcal{P}_M} \big\{\mathbb{E}_F[\omega(\bm{q},\bmt{u})] + \alpha \cdot d_{\rm OT}(F,G)\big\} =  \max_{\bm{q} \geq 0}\;\min_{G \in \mathcal{C}}\;\mathbb{E}_{G}\big[\min_{\bm{u}\geq 0}\{\omega(\bm{q},\bmt{u}) + \alpha\|\bm{u} - \bmt{v}\|_2^2\}\big],
\]
which, by noting that the inner objective is separable, further reduces to
\begin{equation}\label{prob:multi-product}
\textstyle\max\limits_{\bm{q} \geq 0}\;\min\limits_{G \in \mathcal{C}}\;\mathbb{E}_{G}\Big[\sum\limits_{i \in [M]}\min\limits_{u_i\geq 0}\{\pi(q_i, u_i) + \alpha(u_i - \tilde{v}_i)^2\}\Big] = \max\limits_{\bm{q} \geq 0}\;\min\limits_{G \in \mathcal{C}}\;\mathbb{E}_{G}\Big[\sum\limits_{i \in [M]}\Psi_i(\alpha, q_i, \tilde{v}_i)\Big].
\end{equation}
Here, the expression of $\Psi_i(\alpha, q_i, \tilde{v}_i)$ is given in \Cref{lemma:newsvendor type-2}. Applying the standard duality argument for the moment problem, we then obtain the desired result.  
\hfill $\square$

Recall from \cite{natarajan2018asymmetry} that the ambiguity-averse multi-product newsvendor problem
\begin{equation*}
\max_{\bm{q} \geq \bm{0}} \; \min_{G \in \mathcal{A}} \; \mathbb{E}_{G}[\omega(\bm{q},\tilde{\bm{v}})]
\end{equation*}
admits an equivalent dual reformulation
\begin{equation}\label{equ:ambiguity covariance}
\begin{array}{cll}
\max & t + \bm{\lambda}^\top\bm{\mu}+ \left<\bm{Q}, \bm{\Sigma}+\bm{\mu}\bm{\mu}^\top\right>\\ [1mm]
{\rm s.t.} & t + \bm{\lambda}^\top\bm{v} + \big<\bm{Q}, \bm{v}\bm{v}^\top\big> \leq \sum\limits_{i\in[M]}\pi_i(q_i,v_i) ~&\forall \bm{v}\geq 0\\ [1mm]
&\bm{q}\geq 0,\;\bm{Q}\in\mathcal{S}^M,\;\bm{\lambda}\in\mathbb{R}^M,\;t\in\mathbb{R},
\end{array}
\end{equation}
which involves $2^M$ quadratic constraints by a complete expansion of the function $\sum_{i\in[M]}\pi_i(q_i,v_i)$, and is already known to be intractable due to the full covariance structure (\citealp{hanasusanto2015distributionally}). Clearly, the function $\Psi_i(\alpha, q_i, v_i)$ in~\eqref{equ:multi product reformulation} is more complicated than the newsvendor profit function $\pi_i(q_i, v_i)$ as the former depends on the value of $q_i$ (see Figure~\ref{fig:dual tangent} for a visualization). Therefore, the misspecification-averse problem \eqref{equ:multi product reformulation}, compared to \eqref{equ:ambiguity covariance}, is less tractable and more challenging to solve. Specifically, to solve this problem, we first need to partition the decision space into $2^M$ subregions (according to $q_i \leq \frac{p_i}{4\alpha}$ or not), where in each subregion, the function form of $\Psi_i(\alpha,q_i,v_i)$ is fixed. We then solve these $2^M$ subproblems, where for each subproblem, there are at least $2^M$ quadratic constraints in the dual reformulation by expanding $\sum_{i\in[M]}\Psi_i(\alpha, q_i, v_i)$. Therefore, solving a multi-product misspecification-averse problem is much harder than its ambiguity-averse counterpart.

To resolve this issue, we can identify a lower bound of~\eqref{prob:multi-product}, by noting that for each $i\in[M]$,
\[
\Psi_i(\alpha, q_i, v_i) \geq p_i\cdot\min\Big\{q_i, v_i - \frac{p_i}{4\alpha}\Big\} - c_iq_i = L_i(\alpha, q_i, v_i),~\forall q_i \geq 0,\;v_i \geq 0.
\]
Replacing $\Psi_i(\alpha, q_i, v_i)$ with $L_i(\alpha, q_i, v_i)$ in problem \eqref{prob:multi-product}, we then obtain a lower bound. The function $L_i(\alpha, q_i, v_i)$ inherits a similar structure to that of $\pi(q_i, v_i)$. Hence, the resulting problem can also be approximately solved via the quadratic decision rules or semi-definite programming relaxations as discussed in \cite{hanasusanto2015distributionally1, natarajan2017reduced} and \cite{natarajan2018asymmetry1}.

\section{Calibrating $\alpha$ with Limited Knowledge of Distribution Shift}\label{appendix:calibration}

Assume the access to training samples $\hat{v}_1,\ldots,\hat{v}_N$ drawn from the data-generating distribution $D$ and a very small number of testing samples $\hat{u}_1,\ldots,\hat{u}_M$ drawn from the out-of-sample distribution $F$, with $M \ll N$ (in our experiments, $N=30$ and $M=5$). The corresponding empirical distributions are $\hat{D} = \frac{1}{N}\sum_{n=1}^N \delta_{\hat{v}_n}$ and $\hat{F} = \frac{1}{M}\sum_{m=1}^M \delta_{\hat{u}_m}$. This allows us to compute the empirical distribution shift $d_{\rm OT}(\hat{F}, \hat{D})$ used to estimate the underlying $d_{\rm OT}(F,D)$. 

We next describe in detail the two strategies for calibrating $\alpha$ based on the testing samples.

\begin{itemize}

\item {\bf Formula-based approach}. This approach explores formula~\eqref{equ:alpha-N} derived in \Cref{thm:guarantee}, \textit{i.e.}, $\alpha_N = \frac{1}{2}\sqrt{\frac{p(p-c)}{\epsilon_N + d_{\rm OT}(F,D)}}$,
to determine $\alpha$. To implement the formula in practice, it remains to estimate the unknown quantities $d_{\rm OT}(F,D)$ and $\epsilon_N$, which correspond to the following two steps accordingly. 

\textit{Step~1}. We first employ the term $\beta \cdot d_{\rm OT}(\hat{F}, \hat{D})$ as an estimate for the 
distance term $d_{\rm OT}(F,D)$, where  $d_{\rm OT}(\hat{F}, \hat{D})$ is the empirical distribution shift as computed before and $\beta$ is randomly drawn from $\mathbb{U}[0.5, 1]$ as a randomized discount factor.

\textit{Step~2}.  The quantity $\epsilon_N$ involves unknown problem-specific constants in its explicit form given in Theorem~3 and is thus not directly computable. Therefore, to select $\alpha_N$, we treat $\epsilon_N$ as a tuning parameter and perform a $5$-fold cross-validation for $\epsilon_N$ over an interval $[\underline{\epsilon}, \bar{\epsilon}]$ using the training distribution $\hat{D}$. This is actually equivalent to cross-validation of $\alpha_N$ over the interval $\big[\frac{1}{2}\sqrt{\frac{p(p-c)}{\bar{\epsilon} + d_{\rm OT}(F,D)}}, \frac{1}{2}\sqrt{\frac{p(p-c)}{\underline{\epsilon} + d_{\rm OT}(F,D)}}\big]$.

\item {\bf Stress-testing approach}. This approach can be viewed as a \textit{controlled stress test} (\citealp{kupiec2002stress}) that proceeds in two steps: first, constructing a stress-testing distribution; and second, validating against the constructed stress-testing distribution. 

\textit{Step~1}. To mimic the testing environment with distribution shift, we generate a synthetic stress-testing distribution $F_{\rm stress}$ supported on $N$ points that satisfies
$d_{\rm OT}(F_{\rm stress}, \hat{D}) = \beta\cdot d_{\rm OT}(\hat{F}, \hat{D})$,
where $\beta \in [0.5, 1]$ is a discount factor that downscales the empirical distribution shift to account for potential overestimation arising from the imbalance in sample sizes between $\hat{D}$ and $\hat{F}$. Here, we randomly choose the discount factor $\beta$ from a uniform distribution $\mathbb{U}[0.5, 1]$. To achieve this, one admissible construction is $F_{\rm stress} = \frac{1}{N}\sum_{n=1}^N \delta_{\hat{w}_n} 
~~{\rm with}~~
\hat{w}_n = (1-\rho)\hat{v}_n + \rho \hat{v}_*$,
where $\hat{v}_* = \min_{n\in[N]} \hat{v}_n$ and 
$
\rho = \sqrt{\frac{N \beta d_{\rm OT}(\hat{F}, \hat{D})}{\sum_{n\in[N]} (\hat{v}_n - \hat{v}_*)^2}}.
$
Here, $\rho$ is the stress multiplier that determines how much mass is shifted downward so that $F_{\rm stress}$ exactly matches the empirical distance $\beta d_{\rm OT}(\hat{F}, \hat{D})$ from $\hat{D}$. By construction, one can easily verify that the constructed $F_{\rm stress}$ satisfies  
$
d_{\rm OT}(F_{\rm stress}, \hat{D}) = \frac{1}{N}\sum_{n\in[N]}(\hat{v}_n - \hat{w}_n)^2 = \beta d_{\rm OT}(\hat{F}, \hat{D}),
$ 
and the resultant samples of $F_{\rm stress}$ are shifted downward relative to those of the training distribution $\hat{D}$, which is then used for validating  
$\alpha$ in a distribution shifted environment.

\textit{Step~2}. Furthermore, based on the empirical mean $\hat{\mu}$ and standard deviation $\hat{\sigma}$ of $\hat{D}$, the hyperparameter $\alpha$ is then cross-validated over an interval $[\underline{\alpha}, \bar{\alpha}]$ using the generated testing distribution $F_{\rm stress}$.

\end{itemize}


}


\end{document}